\def\ersteseite{1} 
\def\letzteseite{54}
\documentclass[editor]{artjlt-29}
\volume{30}
\annum{2020}

\usepackage{amsmath}
\usepackage{latexsym}
\usepackage{amssymb}
\usepackage{hyperref}

\font\tengoth=eufm10 at 10pt
\font\sevengoth=eufm7 at 6pt
\newfam\gothfam
\textfont\gothfam=\tengoth
\scriptfont\gothfam=\sevengoth

\newcommand{\fS}{{\mathfrak S}}

\newcommand{\g}{{\mathfrak g}}

\newcommand{\z}{{\mathfrak z}}

\newcommand{\fg}{{\mathfrak g}}
\newcommand{\fh}{{\mathfrak h}}

\newcommand{\fk}{{\mathfrak k}}

\newcommand{\fo}{{\mathfrak o}}

\newcommand{\fp}{{\mathfrak p}}

\newcommand{\ft}{{\mathfrak t}}
\newcommand{\fu}{{\mathfrak u}}

\newcommand{\fz}{{\mathfrak z}}

\newcommand\heis{\mathfrak {heis}} 
 
\newcommand\pu{\mathfrak{pu}}

\renewcommand{\:}{\colon}
\newcommand{\1}{\mathbf{1}}

\newcommand{\cA}{\mathcal{A}}

\newcommand{\cD}{\mathcal{D}}
\newcommand{\cE}{\mathcal{E}}

\newcommand{\cH}{\mathcal{H}}

\newcommand{\cL}{\mathcal{L}}
\newcommand{\cM}{\mathcal{M}}

\newcommand{\cO}{\mathcal{O}}
\newcommand{\cP}{\mathcal{P}}
\newcommand{\cQ}{\mathcal{Q}}

\newcommand{\cT}{\mathcal{T}}

\newcommand{\cV}{\mathcal{V}}
\newcommand{\cW}{\mathcal{W}}

\newcommand\bx{{\bf{x}}}

\newcommand{\eset}{\emptyset}

\newcommand{\derat}[1]{\frac{d}{dt}\big\vert_{t = #1}}

\newcommand{\dd}{{\tt d}}

\newcommand{\subeq}{\subseteq}
\newcommand{\supeq}{\supseteq}

\newcommand{\into}{\hookrightarrow}
\newcommand{\eps}{\varepsilon}

\newcommand{\shalf}{{\textstyle{\frac{1}{2}}}}

\newcommand{\N}{{\mathbb N}}
\newcommand{\Z}{{\mathbb Z}}
\newcommand{\R}{{\mathbb R}}
\newcommand{\C}{{\mathbb C}}

\newcommand{\K}{{\mathbb K}}

\newcommand{\bP}{{\mathbb P}}

\renewcommand{\H}{{\mathbb H}}
\newcommand{\T}{{\mathbb T}}

\newcommand{\bS}{{\mathbb S}}

\renewcommand{\hat}{\widehat}
\newcommand{\hats}[1]{\hat{\hat{\hbox{$#1$}}}}

\renewcommand{\tilde}{\widetilde}

\renewcommand{\L}{\mathop{\bf L{}}\nolimits}
\newcommand{\Lie}{\mathop{\bf L{}}\nolimits}

\newcommand{\Meas}{\mathop{{\rm Meas}}\nolimits}

\newcommand{\GL}{\mathop{{\rm GL}}\nolimits}

\newcommand{\AU}{\mathop{{\rm AU}}\nolimits}

\newcommand{\OO}{\mathop{\rm O{}}\nolimits}

\newcommand{\U}{\mathop{\rm U{}}\nolimits}

\newcommand{\aut} {\mathop{{\mathfrak{aut}}}\nolimits}

\newcommand{\gau} {\mathop{{\mathfrak{gau}}}\nolimits}
\newcommand{\gl}  {\mathop{{\mathfrak{gl} }}\nolimits}

\newcommand{\vir} {\mathop{{\mathfrak{vir}}}\nolimits}
\newcommand{\fsl} {\mathop{{\mathfrak{sl} }}\nolimits}
\newcommand{\fsp} {\mathop{{\mathfrak{sp} }}\nolimits}
\newcommand{\su}  {\mathop{{\mathfrak{su} }}\nolimits}

\newcommand{\ad}{\mathop{{\rm ad}}\nolimits}
\newcommand{\Ad}{\mathop{{\rm Ad}}\nolimits}

\newcommand{\epi}{\mathop{{\rm epi}}\nolimits}
\newcommand{\tr}{\mathop{{\rm tr}}\nolimits}

\newcommand{\Hom}{\mathop{{\rm Hom}}\nolimits}

\newcommand{\Ext}{\mathop{{\rm Ext}}\nolimits}

\newcommand{\Isom}{\mathop{{\rm Isom}}\nolimits}

\newcommand{\Herm}{\mathop{{\rm Herm}}\nolimits}

\newcommand{\Aut}{\mathop{{\rm Aut}}\nolimits}

\newcommand{\diag}{\mathop{{\rm diag}}\nolimits}
\newcommand{\End}{\mathop{{\rm End}}\nolimits}
\newcommand{\id}{\mathop{{\rm id}}\nolimits}

\renewcommand{\dim}{\mathop{{\rm dim}}\nolimits}

\newcommand{\im}{\mathop{{\rm im}}\nolimits}

\newcommand{\supp}{\mathop{{\rm supp}}\nolimits}

\newcommand{\cone}{\mathop{{\rm cone}}\nolimits}
\newcommand{\conv}{\mathop{{\rm conv}}\nolimits}

\newcommand{\Spann}{\mathop{{\rm span}}\nolimits}

\newcommand{\der}{\mathop{{\rm der}}\nolimits}

\newcommand{\indlim}{{\displaystyle \lim_{\longrightarrow}}\ }
\newcommand{\prolim}{{\displaystyle\lim_{\longleftarrow}}\ }

\newcommand{\nin}{\noindent} 
\newcommand{\oline}{\overline}

\newcommand{\la}{\langle}
\newcommand{\ra}{\rangle}

\newcommand{\res}{\vert}

\newcommand{\Spec}{{\rm Spec}}

\newcommand{\ssssarr}{\hbox to 15pt{\rightarrowfill}}
\newcommand{\sssarr}{\hbox to 20pt{\rightarrowfill}}
\newcommand{\ssarr}{\hbox to 30pt{\rightarrowfill}}
\newcommand{\sarr}{\hbox to 40pt{\rightarrowfill}}
\newcommand{\arr}{\hbox to 60pt{\rightarrowfill}}
\newcommand{\larr}{\hbox to 60pt{\leftarrowfill}}
\newcommand{\Arr}{\hbox to 80pt{\rightarrowfill}}

\newcommand{\pmat}[1]{\begin{pmatrix} #1 \end{pmatrix}}

\newcommand{\one}{\ensuremath{\mathbf{1}}}

\renewcommand{\phi}{\varphi} 
\renewcommand{\derat}[1]{\frac{d}{dt}\Big\vert_{t=#1}}

\newcommand{\cont}{\mathfrak{cont}}
\newcommand{\Inj}{\mathop{\rm Inj}\nolimits}
\newcommand{\Bij}{\mathop{\rm Bij}\nolimits}
\newcommand\co{\mathop{\rm co}\nolimits}

\newcommand{\bc}{{\bf{c}}}
\newcommand{\bd}{{\bf{d}}}

\title{A Survey on Invariant Cones\\ 
 in Infinite Dimensional Lie Algebras}  
\author{Karl-Hermann Neeb}
\lastname{Neeb} 



\address{%
Karl-Hermann Neeb \\ 
Department Mathematik \\ 
Friedrich-Alexander-Universit\"at \\ 
Erlangen-N\"urnberg \\ 
Cauerstrasse 11 \\ 
91058 Erlangen, Germany \\ 
neeb@math.fau.de} 

\editor{G. \'Olafsson} 
\received{October 28, 2019} 

\begin{document}

\maketitle


\begin{abstract} For the Lie algebra $\g$ of a connected 
infinite-dimensional Lie group~$G$, there is a natural duality 
between so-called semi-equicontinuous weak-$*$-closed convex 
$\Ad^*(G)$-invariant subsets of the dual space $\g'$ and 
$\Ad(G)$-invariant lower semicontinuous positively homogeneous convex functions 
on open convex cones in $\g$. In this survey, we discuss various 
aspects of this duality and some of its applications to a 
more systematic understanding of open invariant cones and convexity 
properties of coadjoint orbits. In particular, we show that 
root decompositions with respect to elliptic Cartan subalgebras 
provide powerful tools for important classes 
of infinite Lie algebras, such as completions
of locally finite Lie algebras, Kac--Moody algebras 
and twisted loop algebras with infinite-dimensional range spaces. 
We also formulate various open problems. 
\end{abstract}

\tableofcontents

\section{Introduction} 
\label{sec:1}

Let $G$ be a connected (possibly infinite-dimensional) Lie 
group  with Lie algebra $\g$ 
(see \cite{Ne06} for a survey on infinite-dimensional Lie groups).  
The conjugation action of $G$ on itself induces on $\g$ the 
{\it adjoint action} 
$\Ad \: G \to \Aut(\g)$ and by dualization we obtain the 
{\it coadjoint action} 
\[ \Ad^* \: G\to \GL(\g')\quad \mbox{  with  } \quad 
\Ad^*_g \lambda := \lambda \circ \Ad_{g^{-1}}.\] 
Here $\g'$ is the space of continuous linear functionals on $\g$. 
We call a subset of $\g$, resp., $\g'$ {\it invariant} if it is invariant 
under $\Ad(G)$, resp., $\Ad^*(G)$. 
This survey is primarily concerned with 
open convex cones $\Omega \subeq \g$ which are invariant 
under the adjoint action of $G$ on $\g$. 

In the context of finite-dimensional groups, invariant cones play a central role 
in the theory of Lie semigroups \cite{HHL89}, 
the theory of causal symmetric spaces \cite{HO97}, 
the complex geometry of semigroup complexifications of Lie groups 
of the type $S = G\exp(i\Omega)$ (Olshanski semigroups) 
\cite{Ne98b}, and,  last but not least, 
in the unitary representation theory of Lie groups, 
where one studies holomorphic extensions to such semigroups (\cite{Ne99}).

The theory of invariant cones in finite-dimensional Lie algebras 
was initiated by Kostant, I.E.~Segal and Vinberg (\cite{Se76}, 
\cite{Vin80}). A classification of simple Lie algebras 
with invariant cones was 
obtained independently by Ol'shanski\u\i{}, Paneitz, and Kumaresan--Ranjan 
(cf.\ \cite{Ols82}, \cite{Pa84}, \cite{KR82}). 
A structure theory of invariant convex cones in general finite-dimensional 
Lie algebras was developed by Hilgert and Hofmann 
in \cite{HH89} and \cite{HHL89}. 
The characterization of finite-dimensional Lie algebras 
containing pointed open invariant cones was obtained 
in \cite{Ne94} in terms of symplectic modules 
of convex type, and these have been classified in 
\cite{Neu00}; see \cite{Ne99} for a rather self-contained 
exposition of this theory.  
Typical Lie algebras containing open invariant cones 
are compact ones with non-trivial center 
(such as $\fu_n(\C)$), 
hermitian Lie algebras 
corresponding to automorphism groups 
of bounded symmetric domains (such as $\fsp_{2n}(\R)$), and semidirect 
products, such as the Jacobi Lie algebra $\heis(\R^{2n}) \rtimes \fsp_{2n}(\R)$ 
(polynomials of degree $\leq 2$ on $\R^{2n}$ with the Poisson bracket). 

In this survey we review some results extending the 
theory of invariant cones to infinite-dimensional Lie algebras. 
Due to the lack of a general structure theory, any effective theory 
has to build on a well-developed  toolbox that applies to important 
specific classes of Lie algebras.   
So one of our goals is to contribute to this toolbox  
by showing how general results on convex set, invariance under 
group actions, fixed point theorems 
and Lie theory can be combined to derive information on invariant cones. 

Some types of invariant cones in infinite-dimensional 
Lie algebras have been discussed in \cite{Ne10}, 
including a classification of invariant cones in the Virasoro algebra 
and typical relations with representation theory.  
Invariant cones also arise naturally in differential geometry: 
In the Poisson--Lie algebra $(C^\infty(M), \{\cdot,\cdot\})$ of a Poisson 
manifold~$M$, we have the cone of positive functions, 
and in the Lie algebra $\cont(M,\eta)$ of contact vector fields 
on the contact manifold $(M,\eta)$, the cone of those vector fields $X$ 
for which $\eta(X) \geq 0$ (\cite{EP00}). 

Any theory of invariant open cones implicitly contains information 
on invariant convex functions: 
If $\cD \subeq \g$ is an open domain and $f \: \cD \to \R$ is a lower-semicontinuous 
$\Ad(G)$-invariant convex function, then its open epigraph 
\[ \Omega_f := \{ (x,t) \in \g \times \R \: f(x) < t \} \] 
is an open invariant cone in the direct Lie algebra sum $\g \oplus \R$. 

Invariant convex functions arise naturally in the unitary representation theory 
of Lie groups, and this connection is an important motivation 
for the interest in invariant cones. 
For a unitary representation $(\pi,\cH)$ of $G$ and $x \in \g$, we write 
$\partial \pi(x) := \frac{d}{dt}\big|_{t=0} \pi(\exp tx)$ 
for the infinitesimal generator of the unitary 
one-parameter group 
$(\pi(\exp tx))_{t \in \R}$ (Stone's Theorem). 
If the representation is {\it smooth} in the 
sense that the subspace $\cH^\infty$ 
of elements with smooth orbit maps is dense in $\cH$, 
then the {\it support functional} 
\begin{equation}
  \label{eq:spi}
s_\pi \: \g \to \R \cup \{\infty\}, \quad 
s_\pi(x) :=\sup\big(\Spec(-i\partial \pi(x))\big)  
\end{equation}
is a lower semicontinous positively homogeneous convex function. 
It encodes the spectral bounds of the selfadjoint 
operators $i\partial\pi(x)$. 
The relation $\partial \pi(\Ad(g)x) = \pi(g) \partial \pi(x) \pi(g)^{-1}$ 
implies that it is $\Ad(G)$-invariant. The representation $\pi$ is 
said to be {\it semibounded} if $s_\pi$ is bounded on some non-empty open 
subset of $\g$. 
We thus obtain a natural connection between 
unitary representations and invariant convexity in the Lie algebra~$\g$. 
For a survey on recent  
progress concerning semibounded representations, 
we refer to \cite{Ne16}. Particularly striking results 
concerning so-called smoothing operators and decomposition theory 
can be found in \cite{NSZ15}. 

A systematic understanding of open invariant cones in particular 
required to verify that unitary representations are semibounded and to 
estimate their support functional~$s_\pi$. It can be used 
to show that, for certain Lie groups, all semibounded representations 
are trivial or that all semibounded representations are bounded, 
which is the case if all non-empty open invariant cones coincide with $\g$. 
This is the case for compact semisimple Lie algebras 
(Remark~\ref{rem:center-compact}), 
the unitary Schatten Lie algebras $\fu_p(\cH)$ for $p > 1$ 
(Proposition~\ref{prop:schatten}), and the 
projective unitary Lie algebra $\pu(\cH) = \fu(\cH)/\R i \one$ of 
a Hilbert space $\cH$ (\cite[Thm.~5.6]{Ne12}). 
This property is also related to interesting 
fixed point properties of the corresponding groups, and 
this kind of information is essential for the method of holomorphic 
induction (see \cite{Ne12, Ne13, Ne14}), 
where one requires that the representation one induces from is bounded. 

We start in Section~\ref{sec:2} with a discussion of 
convex subsets of a locally convex space $E$ and its topological dual~$E'$. 
For our purposes, the most important class of subsets 
of $C \subeq E'$ are those which are semi-equicontinuous in the sense that 
their support functional $s_C(v)  = \sup \la C,  v\ra$ is bounded on some 
non-empty open subset of $E$. This leads us to the Duality 
Theorem~\ref{thm:2.3} between weak-$*$-closed convex subsets 
and positively homogeneous lower semicontinuous convex functions 
on open convex cones. As we want to apply this duality to 
invariant subsets of Lie algebras and their duals, we discuss 
in Subsection~\ref{subsec:2.2} various fixed point theorems, 
based on compactness, equicontinuity and the Bruhat--Tits Theorem. 
In Subsection~\ref{subsec:cox} we provide a description 
of closed convex hulls of orbits of Coxeter groups 
that builds on \cite{HN14} and makes some results 
more  easily applicable in the Lie algebra context. 

In Section~\ref{sec:3} we then turn to invariant convexity 
in Lie algebras, formulate the specific problems (P1-6) one would 
like to answer for infinite-dimensional Lie algebras and 
explain how these relate to unitary representations 
and spectral conditions for corresponding selfadjoint operators 
(Subsection~\ref{subsec:3.0}). Generalizing root decompositions 
of finite-dimensional Lie algebras with respect to a compactly 
embedded Cartan algebra, we introduce in Subsection~\ref{subsec:3.1} 
the notion of an elliptic Cartan subalgebra~$\ft \subeq \g$. Here an important point 
is that it suffices that the sum of the root spaces is dense in $\g_\C$, 
which is more natural in the context of Banach and Fr\'echet--Lie algebras. 
In this rather general context, the Reduction Theorem~\ref{thm:redux} 
tells us that the projection $p_\ft \: \g \to \ft$, resp., the restriction 
map $p_{\ft'} \: \g' \to \ft'$ leave open, resp., 
semi-equicontinuous weak-$*$-closed  invariant convex subsets invariant. 
So these sets can be studied to a large extent in 
terms of their intersections with $\ft$ and $\ft'$ respectively. 
This permits us to generalize various results from 
the finite-dimensional case. 
We conclude Section~\ref{sec:3} with an explanation of 
how the algebraic theory of unitary 
highest weight representations can be used to obtain convexity 
theorems for coadjoint orbits. 

In the brief Section~\ref{sec:4} we recall the convexity 
theorems for adjoint and coadjoint orbits for finite-dimensional Lie algebras 
and refer to \cite{Ne99} for details. 
Several specific classes of infinite-dimensional Lie algebras 
are discussed in Section~\ref{sec:5}. 
For nilpotent and $2$-step solvable algebras everything reduces to the 
abelian case (Theorem~\ref{thm:nilp}), but the class of 
$3$-step solvable algebras contains the oscillator algebras 
which display various types of non-trivial behavior 
(Subsection~\ref{subsubsec:osci}). 
For Kac--Moody Lie algebras we formulate the Kac--Peterson Convexity Theorem and 
we briefly mention some results on Lie algebras of vector fields. 
Since they display the differences between finite and infinite-dimensional theory 
so nicely, we take a closer look at infinite-dimensional 
versions of unitary groups in Subsection~\ref{subsec:unitliealg}. 
We conclude Section~\ref{sec:5} with the observation that projective 
limits of Lie algebras do not lead to new phenomena. 

In Section~\ref{sec:6} we turn to Lie algebras which 
contain a dense directed union 
of finite-dimensional Lie algebras. Typical examples 
are the Lie algebras $\g = \fu(J,\K)$ of skew-hermitian  $J \times J$-matrices 
with finitely many non-zero entries in $\K \in \{\R,\C,\H\}$ 
and their Banach completions. The Lie algebras $\fu(J,\K)$ 
have  the  interesting property that  all semi-equicontinuous orbits 
in $\g'$ and all invariant norms extend to the minimal Banach completion, 
the Lie algebea $\fu_1(\ell^2(J,\K))$ of skew-hermitian trace class operators. 
We also discuss the existence of invariant scalar products on general direct 
limits of compact Lie algebras and show that open invariant cones 
in Hilbert--Lie algebras intersect the center. 

Motivated by direct limits of affine Kac--Moody algebras 
and double extensions of twisted loop algebras, we turn in Section~\ref{sec:7} 
to various aspects of double extensions. The most well-behaved class 
are double extensions of Lie algebras with an invariant scalar product 
(euclidean Lie algebras) because they carry an invariant Lorentzian form 
and this immediately leads to open invariant cones and many semi-equicontinuous 
orbits. However, a finer analysis of the convexity properties of double 
extensions turns out to be quite delicate because of the huge variety 
of double extensions and because many natural representation theoretic 
constructions lead to non-Lorentzian double extension. Here we briefly 
explain the source of these difficulties. 

We collect some auxiliary results in appendices on constructing open cones, 
Lorentzian geometry and the Bruaht--Tits Theorem. \\

{\bf Notation:} For a set $J$, we write $S_J$ for the group of 
all bijections of $J$, and $S_{(J)}$ for the 
subgroup of those elements moving only finitely many elements. 
These groups act on the space $\R^J$ of functions $J \to \R$ 
and preserve the subspace $\R^{(J)}$ of finitely supported functions. 

For $\K \in \{\R,\C,\H\}$, we write $\gl(J,\K) \subeq \End(\K^{(J)})$ 
for the Lie algebras of all $J \times J$-matrices 
with finitely many non-zero entries, 
$\fsl(J,\C) := \{ x \in \gl(J,\C) \:  \tr x = 0\}$, and 
\[  \fu(J,\K) := \{ x \in \gl(J,\K) \: x^* = -x\}.\] 

For topological vector spaces $E,F$, we write $\cL(E,F)$ for the space 
of continuous linear maps $E \to F$ and $\cL(E) = \cL(E,E)$. 
For a subset $C$ of a real linear space $E$ we write 
$\conv(C)$ for its convex hull and $\cone(C)$ for the convex cone 
generated by~$C$.

\section{Open cones and semi-equicontinuous subsets} 
\label{sec:2}

This section is independent of Lie algebraic structures. 
We first recall some basic facts on convex subsets 
of locally convex spaces and the duality between positively homogeneous 
lower semicontinuous convex functions on open cones 
and semi-equicontinuous weak-$*$-closed convex subsets of the topological dual space  
(Subsection~\ref{subsec:2.1}). We discuss some 
fixed point results for group actions on convex sets 
in Subsection~\ref{subsec:2.2}, 
and in Subsection~\ref{subsec:cox} we provide some specific new results 
on convex hulls of orbits of Coxeter group which can be used 
to study adjoint and coadjoint orbits in Lie theory. 

\subsection{Basic facts and concepts} 
\label{subsec:2.1}

Let $E$ be a real locally convex space and $E'$ be its topological dual, i.e., 
the space of continuous linear functionals on~$E$. The space 
of continuous linear endomorphisms of $E$ is denoted $\cL(E)$. 
We write $\la \alpha, v \ra = \alpha(v)$ for the natural pairing 
$E' \times E \to \R$ and endow $E'$ with the weak-$*$-topology, i.e., 
the coarsest topology for which all evaluation maps 
$\eta_v \: E' \to \R,  \eta_v(\alpha) := \alpha(v)$ 
are continuous. 

For a subset $C \subeq E'$, the sets 
\[ C^\star := \{ v \in E \: (\forall \alpha \in C)\ \alpha(v) \geq 0\} \quad 
\mbox{ and }  \quad  B(C) := \{ v \in E \: \inf\la C, v \ra > - \infty\}  \] 
are convex cones. 

A function $f \: E \to \R \cup \{\infty\}$ is said to be {\it convex} 
({\it lower semicontinuous}) 
if its epigraph $\epi(f) = \{(x,t) \in E \times \R\: f(x) \leq t \}$ 
is convex (closed). 
The {\it support functional of a subset $C\subeq E'$} 
\[ s_C \: E \to \R \cup \{ \infty\}, \quad s_C(v)  := \sup \la C,  v\ra,\] 
takes finite values on $- B(C) = B(-C)$. As a sup 
of a family of continuous linear functionals, 
$s_C$ is convex, lower semicontinuous and positively homogeneous. 
Note that $s_C$ does not change if we replace $C$ by its weak-$*$-closed convex 
hull. It is easy to see that 
\begin{equation}
  \label{eq:episupp}
  \epi(s_C) = ((-C) \times \{1\})^\star
\end{equation}

\begin{Definition} We call a subset $C \subeq E'$ {\it semi-equicontinuous} if 
$s_C$ is bounded on some non-empty open subset of $E$ (cf.\ \cite{Ne09}). 
This implies that the cone $B(C)$ has interior points 
and that $s_C$ is continuous on $B(C)^0$ 
(\cite[Prop.~6.8]{Ne08}). 
\end{Definition} 

\begin{Remark} \label{rem:2.3old} 
(a) If the space $E$ is barrelled, which includes in particular 
Banach and Fr\'echet spaces and their locally convex direct limits 
(\cite{GN}), 
then $C \subeq E'$ is semi-equicontinuous 
if and only if $B(C)$ has interior points 
(apply \cite[Thm.~6.10]{Ne08} to the lower semicontinuous function~$s_C$). 

\nin (b) If $C \subeq E'$ is semi-equicontinuous, then the lower semicontinuity 
of $s_C$ implies that, for each $c > 0$, the subset 
$\{ x \in E \: s_C(x) \leq c\}$ is closed with dense interior. 
\end{Remark} 

\begin{Example} \label{ex:2.2} 
(a) A subset $C \subeq E'$ is equicontinuous  
if $\{v \in E \: s_C(v)  \leq 1\}$ is a {$0$-neighborhood}. 
This is equivalent to $B(C) = E$ and $s_C$ being 
continuous. This clearly implies that 
$C$ is  semi-equicontinuous.

\nin (b) If $\Omega \subeq E$ is an open convex cone, then its dual cone 
$\Omega^\star 
= \{ \alpha \in E' \: \alpha\res_{\Omega} \geq 0\}$ 
is semi-equicontinuous because its support functional $s_{\Omega^\star}$ 
is constant $0$ on the open subset~$-\Omega$. 

If, conversely, $C \subeq E'$ is a convex cone, then 
$B(C) = C^\star$,  so that 
$C$ is semi-equicontinuous if and only if $C^\star$ has interior points.
\end{Example}

\begin{Theorem} \label{thm:2.3} {\rm(Duality Theorem)} 
The assignment $C \mapsto s_C \res_{B(-C)^0}$ defines a 
bijection from semi-equicontinuous weak-$*$-closed convex subsets 
$C \subeq E'$ to continuous positively homogeneous convex functions 
defined on non-empty open convex cones. 
\end{Theorem}

\begin{proof} Since $C$ can be reconstructed from $s_C\res_{B(-C)^0}$ by 
  \begin{equation}
    \label{eq:recov}
C = \{ \alpha \in E' \: (\forall v \in B(-C)^0)\, 
\alpha(v) \leq s_C(v) \} 
  \end{equation}
(\cite[Prop.~6.4]{Ne08}), 
the assignment is injective. It remains to show that it is surjective. 

So let $\eset\not = 
\Omega \subeq E$ be an open convex cone and let $f \: \Omega \to \R$ 
be a continuous, positively homogeneous convex function. We claim that 
the weak-$*$-closed convex subset 
\begin{equation}
  \label{eq:scequal}
C_f := \{ \alpha \in E' \: \alpha\res_{\Omega} \leq f \}  
\quad \mbox{ is equicontinuous with } \quad 
 s_{C_f} \res_{\Omega} = f.
\end{equation}
Here $s_{C_f} \leq f$ follows from the definition and since $\Omega$ is open, 
$C_f$ is semi-equi\-con\-tinuous. 
To verify~\eqref{eq:scequal}, 
let $x \in \Omega$. Then the pair $(x,f(x)) \in E \times \R$ lies 
in the boundary of the open convex cone 
$\epi(f)^0 = \{(x,t) \in \Omega \times \R \: f(x) < t\}$. 
By the Hahn--Banach Separation Theorem, there exists a 
non-zero continuous linear functional $(-\alpha, c) \in \epi(f)^\star$ 
with $-\alpha(x) + c f(x) = 0$. Then $c \geq 0$. 
To see that we actually have $c > 0$, assume that $c = 0$. 
Then $\alpha(x) = 0$ and $\alpha\not=0$, so that 
$\alpha(\Omega)$ is a neighborhood of $0$, contradicting 
$-\alpha\res_\Omega \geq 0$. We conclude that $c > 0$,  and 
$c^{-1} \alpha \in C_f$ satisfies $c^{-1}\alpha(x) = f(x)$. 
This proves~\eqref{eq:scequal}. 

The dual  cone $\epi(f)^\star$ is contained in the closed half space 
$E' \times [0,\infty)$. Since it contains all elements of the form 
$(-\alpha,1)$, $\alpha \in C_f$, the intersection of  
$\epi(f)^\star$ with the open half space $E' \times (0,\infty)$ is 
weak-$*$-dense (because open line segments are dense in their closure). 
Hence $\epi(f)^\star$ is generated by $(-C_f) \times \{1\}$ 
as a weak-$*$-closed convex cone. We thus obtain 
with Proposition~\ref{prop:bidual} and \eqref{eq:episupp} 
\begin{equation}
  \label{eq:dualgraph}
\oline{\epi(f)} = (\epi(f)^\star)^\star 
= ((-C_f) \times \{1\})^\star = \epi(s_{C_f}). 
\end{equation}
As a consequence, $B(-C_f) \subeq \oline\Omega$, 
so that 
$B(-C_f)^0 \subeq \Omega$ by Lemma~\ref{lem:bou}(ii). 
Since $\Omega \subeq B(-C_f)^0$ by definition, we have equality. 
This proves that $f = s_{C_f}\res_{B(-C_f)^0}$. 
\end{proof}

\begin{Remark}
Passing to the larger space $E \times \R$, 
any continuous, positively homogeneous convex function 
$f \: \Omega \to \R$ is determined by 
its open epigraph $\epi(f)^0 = \{ (x,t) \: t > f(x)\}$ which is an open 
convex cone. This reduces the analysis of such functions largely 
to the analysis of open convex cones. As follows from 
\eqref{eq:dualgraph}, the open cone $\epi(f)^0$ is the interior 
of the dual cone of the semi-equicontinuous subset 
$(-C_f) \times \{1\} \subeq E' \times \R$. 
\end{Remark}

\begin{Remark} \label{rem:2.5} If $\iota \: F \to E$ is an injective, continuous linear 
map with dense range, then its adjoint defines an injection 
$\iota' \: E' \into F'$. Then any open convex subset $\Omega \subeq E$ 
is determined by its intersection $\Omega \cap  F = \iota^{-1}(\Omega)$, 
which is an open subset of~$F$. Accordingly, any semi-equicontinuous 
subset $C \subeq E'$ defines a 
semi-equicontinuous subset $\iota'(C) \subeq F'$. This observation 
can often be used to study open convex subsets of $E$ through~$F$, 
which may be more accessible. 
\end{Remark}

\begin{Lemma} \label{lem:endo-strong-closed} 
For a weak-$*$-closed convex subset 
$C \subeq E'$, we have: 
\begin{itemize}
\setlength\itemsep{0em}
\item[\rm(a)] If $A \in \cL(E)$, then the  
adjoint $A' \: E' \to E', \lambda \mapsto \lambda \circ A$ 
satisfies $A'C \subeq C$ if and only if 
$s_C(Av) \leq s_C(v)$ for every $v \in E$. 
\item[\rm(b)] $\cL(E)_C := \{ A \in \cL(E) \: A'C \subeq C\}$ is closed 
in the topology of pointwise convergence. 
\end{itemize}
\end{Lemma}

\begin{proof} (a)  By the Hahn--Banach Separation Theorem, 
\[ C = \{ \alpha \in E' \: (\forall v \in E)\, \alpha(v) \leq s_C(v) \}.\] 
Therefore $A'C \subeq C$ is equivalent to 
$s_C \circ A = s_{A'C} \leq s_C$. 

\nin (b) Since $s_C$ is lower semicontinuous, 
the subset $\{ w \in E \: s_C(w) \leq s_C(v)\}$ is closed 
for every $v \in E$. Hence 
(b) follows from (a). 
\end{proof}

\begin{Definition} 
For a convex subset $C$ in the real linear space $E$, we define 
its {\it recession cone} 
\[  \lim(C) := \{ x \in E \: C + x \subeq C \}  
\quad \mbox{ and } \quad H(C) := \lim(C) \cap -\lim(C).\] 
Then $\lim(C)$ is a convex cone and $H(C)$ a linear subspace. 
\end{Definition} 

\begin{Lemma} \label{lem:limcone} {\rm(\cite[Lemma~2.9]{Ne10})} 
If $\eset\not=C \subeq E$ 
is an open or closed convex subset, then the following assertions hold: 
\begin{description}
\setlength\itemsep{0em}
  \item[\rm(i)] $\lim(C) = \lim(\oline C)$ is a closed convex cone. 
  \item[\rm(ii)] $v \in \lim(C)$ if and only if 
$t_j c_j \to v$ for a net with 
$t_j \geq 0$, $t_j \to 0$ and $c_j \in C$. 
  \item[\rm(iii)] If $c \in C$ and $d \in E$ satisfy 
$c + \R_+ d \subeq C$, then $d \in \lim(C)$. 
\item[\rm(iv)] $H(C) = \{0\}$ if and only if $C$ contains no affine lines. 
\item[\rm(v)] $B(C)^\star = \lim(C)$ and $B(C)^\bot = H(C)$. 
\end{description}
\end{Lemma}

\begin{Remark} \label{rem:2.8} 
(a) If $\dim E < \infty$, then a closed convex subset $C\subeq E'$ 
is semi-equicontinuous if and only if it contains no affine lines, i.e., 
if the cone $\lim(C)$ is pointed. This follows from 
$\lim(C) = B(C)^\star$  (Lemma~\ref{lem:limcone}(v)). 

\nin (b) If $E$ is a locally convex space and 
$C \subeq E'$ is equicontinuous, then the recession cone
$\lim(C)$ is trivial because half-lines are not equicontinuous. 
\end{Remark}

\begin{Example} We consider the real Hilbert space 
\[ E := \ell^2 = \Big\{ (a_n)_{n \in \N} \: a_n \in\R, 
\sum_{n=1}^\infty a_n^2 < \infty\Big\}\] 
and identify $E'$ with $\ell^2$. Then 
$C := \{  x \in E' \: \|x\|_\infty \leq 1 \}$ 
is a weak-$*$-closed convex subset.  
The cone $B(C)$ coincides with the subspace $\ell^1 \subeq \ell^2$, hence 
has no interior points. In particular $C$ is not semi-equicontinuous, 
although $\lim(C) = \{0\}$ and $C$ contains no affine lines 
(Remark~\ref{rem:2.8}(a)).  

This example is also interesting in the context of Remark~\ref{rem:2.5}. 
The continuous inclusion $\iota \:  \ell^1 \to \ell^2$ induces 
an inclusion $\iota' \: \ell^2 \into \ell^\infty \cong (\ell^1)'$. 
Now $\iota'(C)$ is contained in the unit ball of $\ell^\infty$, 
hence is equicontinuous on $\ell^1$, 
although $C$ itself is not equicontinuous on~$\ell^2$.
\end{Example}

\begin{Example} \label{ex:2.8} 
(a) If $\cE$ is a real Hilbert space, then the closed unit ball 
$C := \{ \xi \in \cE \: \|\xi\| \leq 1\}$ 
is a weakly compact convex set with $\lim(C) = \{0\}$.

\nin (b) If $\cE$ is a real Hilbert space, then $\cL := 
\R \times \cE \times \R$ carries the Lorentzian bilinear form 
\[ \beta((z,x,t), (z',x',t')) := zt' + z't - \la \xi, \eta \ra,\] 
and 
\[ \cL_+ := \{ \bx := (z,x,t) \: z \geq 0,  \beta(\bx,\bx) \geq 0\}  \] 
is a self-dual closed convex cone. 
For every $m > 0$, the hyperboloid 
\[ H := \{ \bx \in \cL_+ \: \beta(\bx,\bx) \geq m^2 \} \] 
is a weakly closed convex subset with 
$\lim(H) = \cL_+$ 
(cf.~Proposition~\ref{prop:c.2}). 
The closed convex subset 
\[ P :=\cL_+ \cap (\R \times \cE \times \{1\}) 
=  \{ \bx = (z,x,1) \: z \geq \shalf \|x\|^2 \}  \] 
is a paraboloid.  It is the epigraph of the function $\|x\|^2/2$ on $\cE$, 
and its recession cone satisfies 
\[ \lim(P) \subeq \cL_+ \cap (\R \times \cE \times \{0\}) 
= \R_+ (1,0,0),  \] 
hence coincides with the ray $\R_+ (1,0,0)$.\begin{footnote}
{In view of the preceding examples, we may call a weak-$*$-closed convex subset 
$C \subeq E'$ 
{\it elliptic} if ${\lim(C) = \{0\}}$, 
{\it parabolic} if $\lim(C) = \R_+ \alpha$ for some $\alpha \not=0$, and 
{\it hyperbolic} if $\lim(C)$ separates the points of $E$. 
}\end{footnote}
\end{Example}

\begin{Definition} \label{def:subtang} 
For a convex subset $C \subeq E$ and $x \in C$, we write 
\[ L_x(C) := \oline{\R_+ (C - x)} \] 
for the {\it subtangent cone of $C$ in $x$}. Then the restriction 
of a linear functional 
$\alpha \in E'$ to $C$ takes a minimal value in $x$ if and only if 
$\alpha \in L_x(C)^\star$. In particular, $L_x(C)^\star \subeq B(C)$. 
\end{Definition}

\subsection{Fixed points} 
\label{subsec:2.2}

In this subsection we collect several results on the existence 
of fixed points for affine group actions. 
We start with the result for the most general environment, 
where the group has to be compact 
(cf.~\cite[Prop.~2.11]{Ne10}). 

\begin{Proposition}\label{prop:project}
Let $K$ be a compact group acting continuously on the 
complete locally convex space $E$ by the representation 
$\pi \: K \to \GL(E)$. 

\nin{\rm(i)} If $\Omega \subeq E$ is an open or closed  
$K$-invariant convex subset, then $\Omega$ is invariant under the 
fixed point projection 
$p(v) := \int_K \pi(k)v\, d\mu_K(k),$
where $\mu_K$ is a normalized Haar measure on~$K$. In particular, 
$p(\Omega) = \Omega \cap E^K$. 

\nin {\rm(ii)} If $C \subeq E'$ is a weak-$*$-closed convex 
$K$-invariant subset, then 
$C$ is invariant under the adjoint 
$p'(\lambda)v := \lambda(p(v)) = \int_K \lambda(\pi(k)v)\, d\mu_K(k).$
In particular, $p'(C) = C\cap (E')^K$. 
\end{Proposition} 

\begin{Remark} \label{rem:2.15} 
(a) The preceding proposition can also be used to 
study fixed point projections for non-compact groups. 
Let us assume that $G$ contains a directed family $(K_j)_{j \in J}$ 
of compact subgroups for which the union $\bigcup_{j \in J} K_j$ 
is dense in~$G$. If $G$ acts continuously on $E$ and 
$C \subeq E'$ is a weak-$*$-compact $G$-invariant subset, 
then we obtain with Proposition~\ref{prop:project}(ii)  
for each $j \in J$ and $\lambda \in C$ that 
$p_j'(\lambda) := \int_{K_j} \lambda \circ \pi(k)\, d\mu_{K_j}(k) \in C$.
As $C$ is weak-$*$-compact, the net $(p_j'(\lambda))_{j \in J}$ has a cluster 
point in $C$ which is invariant under every $K_j$, hence also under~$G$.

\nin (b) Proposition~\ref{prop:project}(i)  does not extend to the situation 
under (a). Here a typical example arises from the action 
of the group  $S_\infty = S_{(\N)}$ of all finite permutations on 
$E = \R^{(\N)}$. Then $\Omega := \{ x \in \R^{(\N)} \: \sum_j x_j > 0\}$ 
is an open invariant convex cone which contains no fixed point, 
whereas its dual cone $\Omega^\star = (\R_+)^\N$ contains constant functions 
which are fixed points. 
\end{Remark}

Kakutani's Fixed Point Theorem  below  (\cite[Thm.~5.11]{Ru73}) 
weakens the assumption on the group but strengthens 
the assumptions on the invariant subset. 

\begin{Theorem} \label{thm:kakutani} {\rm(Kakutani's Fixed Point Theorem)} 
Let $E$ be a locally convex space and let $G \subeq \GL(E)$ be an 
equicontinuous group. Then each $G$-invariant compact convex 
subset of $E$ contains a $G$-fixed point. 
\end{Theorem}

The following proposition is a bridge between equicontinuity 
and compactness. 
\begin{Proposition} \label{prop:2.17}
Let $E$ be a locally convex space and 
$G \subeq \GL(E)$ an equicontinuous group for which 
$E$ is the union of finite-dimensional 
$G$-invariant subspaces. Then the following assertions hold: \\[-7mm]
\begin{itemize} 
\setlength\itemsep{0em}
\item[\rm(i)] $G$ is  equicontinuous on the dual space 
$E'$ with respect to the weak-$*$-topology.   
\item[\rm(ii)] There exists a continuous representation 
$\rho \: K \to \GL(E)$ of a compact group $K$ for which $G \subeq \rho(K)$ 
is dense in the topology of pointwise convergence.
\end{itemize}
\end{Proposition}

Kakutani's Theorem applies to $G$-invariant 
weak-$*$-compact convex subsets $C \subeq E'$ (cf.\ also Remark~\ref{rem:2.15}(a)), 
but, in view of (ii), one can also apply 
Proposition~\ref{prop:project} to the compact group~$K$.
  
\begin{proof} (i) On $E'$ the weak-$*$-topology coincides with the topology 
of uniform convergence on compact absolutely convex subsets of $E$ 
which are contained in a finite-dimensional subspace. 
Our assumptions imply that $E$ is the union of $G$-invariant 
compact absolutely convex subsets and the polars of these sets form 
a fundamental systems of $G$-invariant $0$-neighborhoods in~$E'$. 

\nin (ii) Consider $E$ as a direct limit $E = \indlim E_j$ 
of finite-dimensional $G$-invariant 
subspaces. Then the representations $(\rho_j,E_j)$ of $G$ on the $E_j$ 
yield a homomorphism 
$G \to \prod_{j \in I} \GL(E_j), g \mapsto (\rho_j(g))_{j \in I},$ 
and since each group $\rho_i(G)$ is relatively compact, we obtain a 
homomorphism of $G$ into the compact group 
$\prod_{i \in I} \oline{\rho_i(G)}$. 
Let $K := \oline{G} \subeq \GL(E)$ denote the closure with respect to the topology 
of pointwise convergence. Then $K$ also preserves all subspaces 
$E_i$, so that it may be considered as a subgroup of the compact group 
$\prod_{i \in I} \oline{\rho_i(G)}$, in which it 
is the closure of the image of $G$ in this group, hence compact. 
That the action of $K$ on $E$ is continuous 
follows from the equicontinuity of $K$ and the fact that 
all its orbit maps are continuous. 
\end{proof}

The following proposition 
is an example of a result not covered by these general 
tools. It can be proved by direct arguments: 
\begin{Proposition} \label{prop:3.2}
Let $J$ be a set. Then the following assertions hold: \\[-8mm]
\begin{itemize}
\setlength\itemsep{0em}
\item[\rm(a)] Every open convex cone $\Omega \subeq \ell^\infty(J)$ 
invariant under the natural action of the full permutation group 
$S_J$ contains a fixed point, i.e., a constant function. 
\item[\rm(b)] For a non-empty open $S_{(J)}$-invariant cone 
$\Omega \subeq \R^{(J)}$, we have:  \\[-7mm]
\begin{itemize}
\setlength\itemsep{0em}
\item[\rm(i)]  There exist disjoint 
finite subsets $\eset\not=F_1, F_2 \subeq J$ and $a,b> 0$ 
with $a e_{F_1} - b e_{F_2} \in \Omega$, where $e_F := \sum_{j \in F} e_j$. 
\item[\rm(ii)]  If $\Omega$ is proper, then the 
summation functional $\chi(x) := \sum_{j \in J} x_j$ satisfies 
$\chi \in \Omega^\star \cup - \Omega^\star$. 
\item[\rm(iii)] $\Omega^\star \subeq \ell^\infty(J,\R)$. 
\end{itemize}
\item[\rm(c)] $\lambda \in \R^J \cong (\R^{(J)})'$ 
has a semi-equicontinuous orbit under $S_{(J)}$ if and only if 
it is bounded, and then its orbit is equicontinuous. 
\end{itemize}
\end{Proposition}

\begin{proof} (a) is \cite[Lemma~3.5]{Ne12}. 

\nin (b) (i): As $\Omega$ is open, there exists an element 
$x \in \Omega$ with $\chi(x) \not=0$. If $F \supeq \supp(x)$ is finite, 
then averaging over $S_{F}$ shows that some multiple of $e_F$ is contained 
in $\Omega$. If $ e_F \in \Omega$, then we put $F_1 := F$ 
and observe that, for any sufficiently small $b > 0$ and any finite subset 
$F_2\subeq J$ disjoint from $F_1$, the element $e_{F_1} - b e_{F_2}$ is contained 
in $\Omega$.

\nin (b)(ii): If $\chi$ is not contained in $\Omega^\star$ or $-\Omega^\star$, 
then there exist $x,y \in \Omega$ with $\chi(x) < 0 < \chi(y)$. Then the argument 
in (i) provides a finite non-empty subset $F \subeq J$ 
with $\pm e_F \in \Omega$. Hence $0 = e_F - e_F \in \Omega$ implies that 
$\Omega = \R^{(J)}$ is not proper. 

\nin (b)(iii): Let $\lambda \in \Omega^\star$ and choose $F_1, F_2$ and $a,b > 0$ 
as in (i) with $x := a e_{F_1} - b e_{F_2} \in \Omega$. Applying 
$\lambda$ to the $S_{(J\setminus F_2)}$-orbit of~$x$, it follows that $\lambda$ 
is bounded from below. We likewise see that $\lambda$ is bounded from above 
by evaluating on the $S_{(J\setminus F_1)}$-orbit of~$x$. 

\nin (c): Clearly, the orbit of a bounded function $\lambda \: J \to \R$ 
is equicontinuous, hence in particular semi-equicontinuous. 
If, conversely, $S_{(J)}\lambda$ is semi-equicontinuous, 
then $B(S_{(J)}\lambda)^0 \subeq \R^{(J)}$ is a non-empty open invariant 
cone. This cone does not change if we add a constant function to $\lambda$. 
Adding a suitable constant $c \in \R$, we see that the cone 
$(\cO_{\lambda + c})^\star$ has interior points, and thus (b)(iii) implies 
that $\lambda$ is bounded. 
\end{proof}

From the Bruhat--Tits Fixed Point Theorem (Theorem~\ref{thm:BT}), we obtain: 

\begin{Proposition} \label{prop:2.16} 
Let $\sigma \: G \to \cH \rtimes \OO(E)$ define an affine 
isometric action of the group $G$ on the euclidean space $E$ 
and $\eset\not=\Omega \subeq E$ be an invariant open or closed convex subset. 
Then the following assertions hold: \\[-7mm]
\begin{itemize}
\setlength\itemsep{0em}
\item[\rm(i)] If $E$ is complete and $\Omega$ contains a bounded $G$-orbit, 
then it contains a fixed point. 
\item[\rm(ii)] Suppose that the affine subspace $E^G$ of fixed points is 
non-empty and complete. Then the orthogonal 
projection $p \: E \to E^G$ satisfies $p(\Omega) = \Omega \cap E^G.$ 
\end{itemize}
\end{Proposition}

\begin{proof} (i) By Lemma~\ref{lem:1.3}, every proper open invariant convex 
subset of $\cH$ is exhausted by closed invariant convex subsets. 
Hence it suffices to prove the assertion if $\Omega$ is closed. 
Then $\Omega$ is a Bruhat--Tits space with respect to the induced 
metric (Example~\ref{ex:BTexample}), so that 
the assertion follows from the Bruhat--Tits Theorem~\ref{thm:BT}. 

\nin (ii) As in (i), it suffices to consider the case where $\Omega$ is closed. 
As $E^G$ is complete, we have $E = E^G \oplus (E^G)^\bot$. 
Let $\hat E$ denote completion of $E$, extend the $G$-action
 to $\hat E$, and write $\hat \Omega$ for the closure of $\Omega$ in $\hat E$. 
Then $\hat E = E^G \oplus \hat{(E^G)^\bot}$ shows that 
$\hat E^G = E^G$. 

We have to show that $c \in \Omega$ implies $p(c) \in \Omega$. 
The subset $C := p^{-1}(p(c)) \cap \hat\Omega \subeq \hat E$ is a 
closed convex $G$-invariant subset containing $c$. 
As $G$ fixes $p(c)$, the $G$-orbit of $c$ is bounded, so that 
$C$ contains a fixed point by (i). Since $p(c)$ is the only $G$-fixed 
point in $p^{-1}(p(c))$, it follows that $p(c) \in C \subeq \hat\Omega$. 
The closedness of $\Omega$ in $E$ implies that 
$\hat \Omega \cap E = \Omega$, so that we obtain 
$p(c) \in E^G \cap \hat \Omega = E^G \cap \Omega.$ 
\end{proof}

\subsection{Coxeter groups} 
\label{subsec:cox}

In \cite{HN14} we studied convex hulls of orbits of Coxeter groups 
to apply these results to orbit projections in infinite-dimensional Lie algebras 
(see~\cite{Ne14, MN17} for some concrete applications). 
Here we  extend some of these, so that they fit better 
the typical requirements in the Lie algebra context. 

\begin{Definition} (a) 
Let $V$ be a finite-dimensional real vector space. 
A {\it reflection data on $V$} consists of a finite family 
$(\alpha_s)_{s \in S}$ of linear functionals on $V$ 
and a family $(\alpha^\vee_s)_{s \in S}$ of elements of $V$ 
satisfying 
\[\alpha_s(\alpha^\vee_s) = 2 \quad \mbox{ for } \quad s \in S.\] 
Then $r_s(v) := v - \alpha_s(v) \alpha^\vee_s$ 
is a reflection on $V$. We write $\cW := \la r_s \: s \in S \ra 
\subeq \GL(V)$ for the subgroup generated by these reflections, 
which is a Coxeter group (\cite{Vin71}, \cite{Hu92}). 
The cone 
\[ K := \{ v \in V \: (\forall s \in S)\, \alpha_s(v) \geq 0 \}  \] 
is called the {\it fundamental chamber} 

\nin (b) A reflection data of finite type is called a  {\it linear Coxeter system} 
(cf.\ \cite{Vin71}) if \\[-7mm]
\begin{description}
\setlength\itemsep{0em}
\item[\rm(LCS1)] $K$ has interior points. 
\item[\rm(LCS2)] $(\forall s \in S)\, \alpha_s \not\in \cone(\{ \alpha_t \: 
t\not=s\})$. 
\item[\rm(LCS3)] $(\forall w \in \cW\setminus \{e\})
\, wK^0 \cap K^0 =\eset$. 
\end{description}
\end{Definition}

\begin{Theorem} \label{thm:1.6} {\rm(\cite[Thm.~2]{Vin71})} 
For a linear Coxeter system 
$(V, (\alpha_s)_{s \in S}, (\alpha^\vee_s)_{s \in S})$, the 
following assertions hold: 
    \begin{itemize}
\setlength\itemsep{0em}
    \item[\rm(i)] The subset $\cT := \cW K$  is a convex cone. 
    \item[\rm(ii)] $\cW$ acts discretely on the interior $\cT^0$ of $\cT$, i.e., 
point stabilizers are finite, $\cT^0/\cW$ is Hausdorff, and 
 every $x \in \cT^0$ has a neighborhood $U$ with 
$wU \cap U =\eset$ for $w \in \cW \setminus \{e\}$.
    \item[\rm(iii)] An element $x \in K$ is contained 
in $\cT^0$ if and only if the stabilizer $\cW_x$ is finite. 
    \end{itemize}
\end{Theorem}

The cone $\cT$ is called the {\it Tits cone} of the linear Coxeter system.
We write 
\begin{equation}
  \label{eq:co}
 \co(v) :=  \conv(\cW v)\quad \mbox{ and } \quad 
\oline\co(v) := \oline\conv(\cW v) 
\end{equation}
for the (closed) convex hull of a $\cW$-orbit. 
We define two cones 
\[ C_S := \cone\{ \alpha^\vee_s \: s \in S \} \subeq V, \qquad 
C_S^\vee := \cone\{ \alpha_s \: s \in S \} \subeq V^*. \] 

A {\it reflection in $\cW$} is an element 
conjugate to some $r_s$, $s \in S$. Any reflection 
can be written as $r_\alpha(v) = v - \alpha(v) \alpha^\vee$ 
with $\alpha \in V^*$ and $\alpha^\vee \in V$, where 
$\alpha$ belongs to the set 
$\Delta := \cW \{ \alpha_s \: s \in S \}$ of {\it roots} and 
$\alpha^\vee \in \Delta^\vee := \cW \{ \alpha^\vee_s \: s \in S \}$ is a 
{\it coroot}. To $v \in E$, we associate  the cone \\[-5mm]
\[  C_v := \cone \{ \alpha^\vee \: \alpha(v) >  0 \}.\] 
We also recall 
from \cite[Rem.~2.5]{HN14} that the  subset $\Delta^+ := \Delta \cap K^\star$ 
of positive roots satisfies $\Delta = \Delta^+ \dot\cup - \Delta^+$ 
and 
\begin{equation}
  \label{eq:dag}
\Delta^+ \subeq  C_S^\vee.
\end{equation}

\begin{Lemma} \label{lem:orb-cont} 
Let $C\not=\eset$ be an open convex cone in the finite-dimensional 
real vector space~$V$ and $\Gamma \subeq \GL(V)$ be a subgroup such that 
$|\det(\gamma)| = 1$ for all $\gamma \in \Gamma$. Then 
\[ \oline\conv(\Gamma.v) \subeq C \quad \mbox{ for all } \quad v \in C.\] 
\end{Lemma}

\begin{proof} From \cite[Thm.~V.5.4]{Ne99} we know that there exists a 
smooth {$\Gamma$-invariant}  convex
function $\phi \: C \to (0,\infty)$ such that 
$x_n \to x \in \partial C$ and $x_n \in C$ implies 
$\phi(x_n) \to  \infty$. 
Therefore all sublevel sets $\{\phi \leq c\}$, $c > 0$, are convex subsets of $C$ 
which are closed in $V$. As $\phi$ is $\Gamma$-invariant, the assertion follows from 
$\conv(\Gamma v) \subeq \{ \phi \leq \phi(v)\}$. 
\end{proof}

\begin{Lemma}\label{lem:2.2} If $v \in \cT$ and 
$w \in C_v$, 
then there exists an $\eps > 0$ such that 
$v - \eps w \in \co(v)$. 
\end{Lemma}

\begin{proof} We write $w = \sum_{j = 1}^k c_j \alpha^\vee_j$ 
with $c_j > 0$ and $\alpha_j(v) > 0$. 
Then $r_{\alpha_j}(v) = v - \alpha_j(v) \alpha^\vee_j \in \co(v)$ 
implies that 
\[ v - \sum_{j = 1}^k d_j  \alpha_j(v) \alpha^\vee_j \in \co(v) 
\quad \mbox{ for } \quad d_j \geq 0, \sum_j d_j \leq 1.\] 
For $t > 0$ we compare this expression with 
\[ v_t := v - t w 
=  v - \sum_{j = 1}^k t c_j \alpha^\vee_j 
=  v - \sum_{j = 1}^k t \frac{c_j}{\alpha_j(v)} \alpha_j(v) \alpha^\vee_j,\] 
and find that $v_t \in \co(v)$ if 
$t \sum_{j = 1}^k \frac{c_j}{\alpha_j(v)} \leq 1.$ 
This prove the lemma. 
\end{proof}

\begin{Lemma} \label{lem:2.3} 
For $v \in K$, the following assertions hold: 
\begin{itemize}
\setlength\itemsep{0em}
\item[\rm(i)] $C_v \subeq C_S$ with equality for $v \in K^0$.   
\item[\rm(ii)]  If $v \in  \cT^0$, then the stabilizer $\cW_v$ is finite and 
$\bigcap_{w \in \cW_v} wC_S = C_v.$ 
This implies in particular that $C_v$ is closed and that 
\begin{equation}
  \label{eq:interchar}
\bigcap_{w \in \cW_v} w(v - C_S) = v - C_v
\quad \mbox{ and } \quad 
\bigcap_{w \in \cW} w(v - C_v) =  \bigcap_{w \in \cW} w(v - C_S).
\end{equation}
\end{itemize}
\end{Lemma}

\begin{proof} (i) If $\alpha(v) > 0$ for some $v \in K$ and $\alpha \in \Delta$, then 
$\alpha \in \Delta^+ \subeq C_S^\vee$ by \eqref{eq:dag}. 
Write $\alpha = w \alpha_s \in \Delta$ for some 
$w \in \cW$ and $s \in S$. Then 
$\alpha^\vee = w\alpha_s^\vee \in C_S$ follows from 
$\alpha = w\alpha_s \in C_S^\vee$ by \cite[Thm.~1.10]{HN14}. This proves $C_v \subeq C_S$. 
If $v \in K^0$, then $\alpha_s(v) > 0$ for every $s \in S$, so that 
$C_S \subeq C_v$ and thus $C_S = C_v$. 

\nin (ii) As $C_v \subeq C_S$ for $v \in K$, and $C_v$ is $\cW_v$-invariant,
it follows that 
$C_v \subeq \bigcap_{w \in \cW_v} wC_S.$ 

Let $S_v := \{ s \in S \: \alpha_s(v) = 0\}$ and observe that 
$(V, (\alpha_s)_{s \in S_v}, (\alpha^\vee_s)_{s \in S_v})$ 
also is a linear Coxeter system (\cite[Rem~1.5]{HN14}). 
We know from \cite[Prop.~1.12]{HN14} that $\cW_v = \la r_s \: s \in S_v\ra$. 
If $\alpha \in \Delta^+$ vanishes on $v$, then $r_\alpha\in \cW_v$, and 
thus $\alpha \in \cW_v \{ \alpha_s \: s \in S_v\}$. We conclude that 
$\alpha^\vee \in C_{S_v}$. This shows that 
$C_S = C_v + C_{S_v},$ 
and this leads to 
\[ \bigcap_{w \in \cW_v} wC_S 
\subeq \bigcap_{w \in \cW_v} w(C_v + C_{S_v}) 
= C_v,\]
where the last equality follows from \cite[Cor.~V.2.10]{Ne99}, applied with 
$x = 0$. 
The first equality in \eqref{eq:interchar} now follows 
from $w(v - C_S) = v - wC_S$ for $w \in \cW_v$. 

Finally, the second equality in \eqref{eq:interchar} follows from 
\[ \bigcap_{w \in \cW} w(v - C_v) 
=  \bigcap_{w \in \cW} \bigcap_{w' \in \cW_v} ww'(v - C_S) 
=  \bigcap_{w \in \cW} w(v - C_S).\qedhere\] 
\end{proof}

\begin{Theorem} \label{thm:coxconvtheo}{\rm(Convexity Theorem 
for Coxeter groups)} 
For $v \in \cT$, we have 
  \begin{equation}
    \label{eq:inc1}
\oline\co(v) \cap \cT = \cT \cap \bigcap_{w \in \cW} w(v - \oline{C_v}).
  \end{equation}
If $v \in \cT^0$, then 
$C_v$ is closed, $\oline\co(v) \subeq \cT^0$, and 
  \begin{equation}
    \label{eq:inc2}
\oline\co(v)  
= \cT \cap \bigcap_{w \in \cW} w(v - C_v) 
= \cT \cap \bigcap_{w \in \cW} w(v - C_S).
  \end{equation}
\end{Theorem}

For the special case where $v \in K^0$, this theorem is contained 
in \cite[Satz~3.21]{Ho99}, which was never published. 
For the case where $\cW$ is finite, it is contained in \cite[Prop.~V.2.9]{Ne99}. 

\begin{proof} From \cite[Thm.~2.7]{HN14} we obtain 
$\co(v) \subeq v - C_v,$ hence 
$\oline\co(v) \subeq v - \oline{C_v}$. 
Since $\co(v)$ is $\cW$-invariant, 
it follows that 
$\oline\co(v) \subeq \bigcap_{w \in \cW} w(v - \oline{C_v}).$ 
Intersecting with $\cT$ thus proves ``$\subeq$'' in \eqref{eq:inc1}. 

To verify ``$\supeq$'',
we use $\cT = \cW K$ to see that 
it suffices to show that, for $v \in K$, any 
$u \in K \cap (v- \oline{C_v})$ is contained in $\oline\co(v)$. 

\nin {\bf Case 1:} $u \in v - C_v$. 
Consider the line segment 
\[ \gamma \: [0,1] \to K, \qquad \gamma(t) := v + t(u-v).\] 
Then $\gamma^{-1}(\oline\co(v)) = [0,c]$ for some $c \in [0,1]$ 
and we have to show that $c = 1$. 

As $u,v \in K$, we have $\gamma(t) \in K$ for $0 \leq t \leq 1$. 
Hence $\alpha(u), \alpha(v) \geq 0$ for $\alpha \in \Delta \cap K^\star$, and if 
$\alpha(v) > 0$, then $\alpha(\gamma(t)) > 0$ for $0 \leq t < 1$. 
This implies that $C_v \subeq C_{\gamma(t)}$ for $t < 1$. 

Suppose that $c < 1$. For $v' := \gamma(c) \in \oline\co(v)$ we then obtain 
$u - v' = (1-c)(u-v) \in - C_v \subeq - C_{v'}$. 
By Lemma~\ref{lem:2.2}, there exists an $\eps \in (0,1)$ such that 
$v' + \eps (u-v') \in \oline\co(v') \subeq \oline\co(v)$. As
\[ v' +\eps (u-v') 
= v +  c(u-v) + \eps(1-c)(u-v) 
= \gamma(c + \eps(1-c)),\] 
this contradicts the maximality of $c$, so that we have shown that~$c = 1$ 
and hence that $(v - C_v) \cap K \subeq \oline\co(v)$.

\nin {\bf Case 2:} $u \in v - (\oline C_v \setminus C_v)$. \\
\nin {\bf Case a:} If $v - \oline{C_v}$ intersects $K^0$, then 
$(v - \oline{C_v}) \cap K^0$ is dense in $(v - \oline{C_v}) \cap K$, 
and since $v - C_v$ is dense in $v - \oline{C_v}$, the set
$(v - C_v) \cap K^0 \subeq (v - C_v)\cap K$ 
is also dense in $(v - \oline{C_v}) \cap K$. 
In Case 1 we have seen that 
$(v - C_v)\cap K \subeq \oline\co(v)$, so that we also obtain 
$(v - \oline{C_v})\cap K \subeq \oline\co(v)$. 

\nin {\bf Case b:} $(v - \oline{C_v}) \cap K^0 = \eset$: 
By definition of $K^0$, 
this is equivalent to the existence of an $s \in S$ with 
$\alpha_s(v - \oline{C_v}) \subeq (-\infty,0]$. 
As $\co(v) \subeq v - C_v$ (\cite[Thm.~2.17]{HN14}), this entails 
$\alpha_s \in -\co(v)^\star$, so that $r_{\alpha_s}^* \alpha_s = - \alpha_s$ 
further leads to $\alpha_s \in \co(v)^\bot$. 

We now reduce this case to the first one, for a smaller linear 
Coxeter system. We write 
\[ S = S' \dot\cup S'' \quad \mbox{ with } \quad 
S' := \{ s \in S \: \alpha_s(C_v) \not=\{0\}\}, \qquad 
S'' := \{ s \in S \: \alpha_s(C_v) =\{0\}\}.\] 
By \cite[Rem~1.5]{HN14}, 
$(V, (\alpha_s)_{s \in S'}, (\alpha^\vee_s)_{s \in S'})$ 
also is a linear Coxeter system. 
For $t \in S'$ and $s \in S''$ we have 
$\alpha_s(\alpha^\vee_t) = 0$ and by \cite[Prop.~1.4]{HN14} also 
$\alpha_t(\alpha^\vee_s) = 0$. This implies that 
$r_s$ commutes with $r_t$. For 
\[ \cW' := \la r_s \: s \in S' \ra \quad \mbox{ and }  \quad 
\cW'' := \la r_s \: s \in S'' \ra \] 
we therefore obtain $\cW = \cW' \cW''$, where both subgroups commute,  
and $\cW''$ fixes all points in $\co(v)$. In particular, we obtain 
\[ \co(v) = \co'(v) := \conv(\cW'v).\] 
Next we observe that 
$\oline{C_v} = \oline\cone(v-\co(v)) = \oline{C_v'}$ 
and 
\[ K' := \{ v \in V \: (\forall s \in S')\, \alpha_s(v) \geq 0 \}  
\supeq K,\] 
so that it suffices to show that 
$(v - \oline{C_v'}) \cap K' \subeq \oline\co'(v).$ 
By definition of $S'$, no $\alpha_s$, $s \in S'$, vanishes on $C_v'$, 
so that this follows from Case 2a. 

We thus obtain $(v - \oline{C_v}) \cap K \subeq \oline\co(v)$ 
in all cases, and this proves \eqref{eq:inc1}. 

If, in addition, $v \in \cT^0$, then 
Lemma~\ref{lem:orb-cont} and $|\det(r_s)| = 1$ for $s \in S$ imply 
that $\oline\co(v) \subeq \cT^0$, 
so that \eqref{eq:inc2} follows from Lemma~\ref{lem:2.3}(ii). 
\end{proof}

\begin{Examples} \label{ex:cox} For more details on typical 
classes of Coxeter groups, we refer to \cite{Hu92}. 

\nin (a) For finite Coxeter groups we have 
$\cT = E$. They arise in particular from finite-dimensional Lie groups 
as Weyl  groups (Definition~\ref{def:weylgrp} 
and \cite[Prop.~VII.2.10]{Ne99}). 
A typical example is the action of the symmetric 
group $S_n$ on the euclidean space $\R^n$, for 
$\alpha_j(x) = x_{j+1} - x_j$ and $\alpha_j^\vee = e_{j+1} - e_j$, 
$j = 1,\ldots, n-1$. 

\nin (b) For affine Coxeter groups, the space  
$E'$ carries a positive semidefinite $\cW$-invariant 
form degenerate along a one-dimensional subspace $\R \lambda$, 
so that $\cW$ fixes $\lambda$. Then the 
hyperplane $H = \lambda^{-1}(1) \subeq E$ is $\cW$-invariant 
and the interior $\cT^0 = \lambda^{-1}(\R^\times_+)$ of the Tits cone 
is an open half space. In this case $\cW$ can be viewed as a 
group generated by affine euclidean reflections on~$H$. 

\nin (c) Other important classes of Coxeter groups 
preserve a Lorentzian form and can be viewed as isometry 
groups of hyperbolic spaces (\cite{Vin71}, \cite{Ka90}). 
\end{Examples}

\section{The Lie algebra context} 
\label{sec:3}

\subsection{The problems} 
\label{subsec:3.0}

Let $G$ be a connected Lie group with Lie algebra $\g$. 
We are aiming at a systematic analysis of 
$\Ad^*(G)$-invariant semi-equicontinuous subsets $C \subeq \g'$, resp., 
$\Ad(G)$-invariant continuous positively homogeneous convex functions 
on open cones in $\g$. By the Duality Theorem~\ref{thm:2.3} 
and the subsequent remark, 
it basically suffices to study semi-equicontinuous 
invariant subsets $C \subeq \g'$ 
and open invariant cones in $\g$, resp., $\g \oplus \R$. 
This leads to the following problems: 
\begin{itemize}
\setlength\itemsep{0em}
\item[\rm(P1)] Classify the semi-equicontinuous 
coadjoint orbits $\cO_\lambda := \Ad^*(G)\lambda \subeq \g'$. 
\item[\rm(P2)] Classify open invariant convex cones in $\g$. 
\end{itemize}

Semi-equicontinuous weak-$*$-closed convex subsets $C \subeq \g'$ 
are locally compact (with respect to the weak-$*$-topology) and
every $x \in B(C)^0$ defines a proper evaluation function  
$\eta_x \: C \to \R, \alpha \mapsto \alpha(x)$. 
In particular, there exists an extreme point $\lambda$ 
minimizing $\eta_x$ (\cite[Prop.~6.13]{Ne08}). We then have
\begin{align}
  \label{eq:min}
\lambda(x) = \min \cO_\lambda(x) = \min \lambda(\cO_x)
\quad \mbox{ for } \quad 
\cO_x = \Ad(G)x.  \tag{{\rm min}}\end{align} 
For $y \in \g$, the smooth function $f(t) := \lambda(e^{t \ad y}x)$ then 
satisfies $0 = f'(0) \leq f''(0)$, which leads to the necessary conditions 
\begin{align}
  \label{eq:nec-cond}
  \lambda([x,y]) = 0 \quad \mbox{ and } \quad 
\lambda([y,[y,x]])\geq 0 \quad \mbox{ for } \quad y \in \g. \tag{{\rm pos}}
\end{align}
Here elements $x \in \g$ which are extreme points in $\conv(\cO_x)$ are most 
interesting. 

\begin{itemize}
\item[\rm(P3)] Determine natural conditions under which 
\eqref{eq:nec-cond} implies \eqref{eq:min}. 
\end{itemize}

{\bf Unitary representations:} 
If $(\pi, \cH)$ is a smooth unitary representation of $G$, 
then its support functional $s_\pi$ 
(see \eqref{eq:spi} in the introduction), coincides with the 
support functional of the range of the momentum map 
on the projective space of $\cH^\infty$: 
\[ \Phi_\pi \: \bP({\cal H}^\infty)\to \g' \quad \hbox{ with } \quad 
\Phi_\pi([v])(x) 
= \frac{\la v,  -i\partial\pi(x)v \ra}{\la v, v \ra}\quad \mbox{ for } 
[v] = \C v. \]
The representation is semibounded if and only if 
the range of $\Phi_\pi$ is semi-equicontinuous. 
Then the {\it momentum set} $I_\pi := \oline\conv(\im \Phi_\pi)$ 
also has this property. 

\begin{itemize}
\item[\rm(P4)] Calculate the momentum set $I_\pi$ for 
irreducible semibounded representations $\pi$ as concretely as possible. 
It is enough to do this for irreducible ones 
because we have direct integral 
decompositions for semibounded representations (cf.~\cite{NSZ15}), 
so that \cite[Thm.~X.6.16]{Ne99} on the momentum set of direct integrals 
reduces the problem to the irreducible case. As we shall see below, 
one may expect that 
\begin{equation}
  \label{eq:momset}
I_\pi = \oline\conv(\cO_\lambda) \quad \mbox{ for some } \quad 
\lambda \in \g' 
\end{equation}
in all cases where the representation has an extremal weight $\lambda$ 
(see \cite[\S 5.2]{Ne10} and Examples~\ref{ex:3.18a} and \ref{ex:3.18}). 
\end{itemize}

\nin {\bf Double extensions 
and projective unitary representations}: 
In the context of unitary representations, the examples arising 
in physics \cite{Ot95, PS86, SG81, Ne10} often lead to 
projective unitary representations $(\dd\pi, \cH^\infty)$ 
of a semidirect product Lie algebra $\g^\sharp = \g \rtimes_D \R$, where 
$D \in \der(\g)$ and $\bd := (0,1)$ 
satisfies the {\it positive energy condition} 
$-i\dd\pi(\bd) \geq 0$ (see also \cite{Ne14}, \cite[\S 8]{JN19}). 
Here $\dd\pi(x) := \partial\pi(x)\res_{\cH^\infty}$ and 
$\cH^\infty$ is the space of smooth vectors.

Lifting to a unitary representation of a central extension, we thus 
obtain a Lie algebra of the form  
\begin{equation}
  \label{eq:doub-ext}
 \hat\g = (\R \oplus_\omega \g) \rtimes \R \bd, \qquad \bd = (0,0,1),   
\end{equation}
with Lie bracket 
\begin{equation}
  \label{eq:doublextbrack}
[(z,x,t), (z',x',t')] = (\omega(x,x') + t\delta(x') -t' \delta(x), [x,x'] + tDx' - t' Dx, 0),
\end{equation}
where $\omega \: \g\times \g \to \R$ is a Lie algebra $2$-cocycle and 
$[\bd, (z,x)] = (\delta(x), Dx)$ 
is a lift of the derivation $D$ to the central 
extension $\tilde\g = \R \oplus_\omega \g = \R \bc \oplus \g$. 
Here $\delta$, $D$ and $\omega$ 
are related by the equation 
\begin{equation}
  \label{eq:domegarel}
\omega(Dx,x') + \omega(x,Dx') = \delta([x,x']) \quad \mbox{ for } \quad 
x,x' \in \g.
\end{equation}

\begin{Definition} \label{def:doubleext} 
We call $\hat\g$ the 
{\it double extension} defined by the pair $(\omega, D)$ (if $\delta =0$), 
resp., by $(\omega, D, \delta)$. 
\end{Definition}

As we may assume that $\lambda :=(1,0,0)\in\hat\g'$ 
is contained in the momentum set (cf.\ \cite[Prop.~15]{JN16}), we then obtain 
\begin{equation}
  \label{eq:1b}
 \inf \lambda(\cO_\bd)=  \inf \cO_\lambda(\bd)  > - \infty. 
\end{equation}
A particularly interesting case arises if 
there exists an eigenvector $v_0$ of $-i\dd\pi(\bd)$ corresponding to the 
minimal spectral value $0$. Then $\lambda := \Phi_\pi([v_0])$ 
satisfies the minimality condition 
$\inf \cO_\lambda(\bd) = \lambda(\bd)=0$, and 
the necessary conditions from \eqref{eq:nec-cond} take the form 
\begin{equation}
  \label{eq:doub-nec-con}
\delta = 0 \quad \mbox{ and } \quad 
 \omega(Dx,x) \geq 0 \quad \mbox{ for }\quad x \in \g.
\end{equation}
In particular the symmetric bilinear form 
$\omega(Dx,y)$ is positive semidefinite. 
For oscillator algebras ($\g$ abelian), these conditions 
were the starting point for the classification results in \cite{NZ13} 
(see Subsection~\ref{subsubsec:osci}
for details). 
Accordingly, one has to explore the structural implications of \eqref{eq:doub-nec-con} 
for other classes of Lie algebras. This leads to our last general problem: 
\begin{itemize}
\item[\rm(P5)] Classify the pairs $(\omega, D)$ satisfying 
\eqref{eq:domegarel} with $\delta =0$ and the {\it positive energy condition} (PEC): 
\begin{align}
\label{eq:pec}
 \inf \cO_\lambda(\bd) = \lambda(\bd) = 0.\tag{PEC}
\end{align}
Determine for which of these pairs $(\omega,D)$ 
the coadjoint orbit $\cO_\lambda$ is semi-equicontinuous. 
\end{itemize}
\nin The positive energy condition is very natural because it 
is often satisfied in examples related to physics, where $-i\dd \pi(\bd)$ 
corresponds to a Hamilton operator (\cite{Ne10}, \cite{Ot95}). 
For some concrete results concerning the case where $\hat\g$ is locally 
finite, we refer to \cite{MN16}. 

\begin{Remark} \label{rem:dcone} 
If the convex hull of the adjoint orbit $\cO_\bd$ has 
interior points in the hyperplane $\R \times \g \times \{1\}$, 
then the cone generated by $\cO_\bd$ has interior points, and therefore 
$\cO_\bd^\star$ consists of semi-equicontinuous orbits 
(cf.~Lemma~\ref{lem:doubcone} for such situations). 
\end{Remark}

\subsection{Root decompositions} 
\label{subsec:3.1}

\begin{Definition} \label{def:basic} (a) We call an abelian subalgebra 
$\ft$ of the locally convex Lie algebra $\g = \L(G)$ 
an  {\it elliptic Cartan subalgebra} 
if $\ft$ is maximal abelian, the root spaces 
\[ \g_\C^\alpha := \{ x \in \g_\C \: (\forall h \in \ft_\C) [h,x]= \alpha(h)
x\} \] 
span a dense subspace $\g_\C^{\rm alg} \subeq \g_\C$, 
and 
\[ \alpha(\ft) \subeq i \R \quad \mbox{ for every  root } 
\quad \alpha \in \Delta :=  
\{ \alpha \in \ft_\C^* \setminus \{0\} \: 
\g_\C^\alpha\not= \{0\}\}.\]

\nin (b) Since $[\g_\C^\alpha, \g_\alpha^\beta] \subeq \g_\C^{\alpha+\beta}$, the 
subspace $\g^{\rm alg}_\C$ is a Lie subalgebra of $\g_\C$ and it intersects 
$\g$ in the dense subalgebra $\g^{\rm alg} := \g_\C^{\rm alg}\cap \g$. 
We shall always {\bf assume} that the canonical projection 
$\g^{\rm alg} \to \ft$ extends to a continuous projection 
$p_\ft \: \g \to \ft$, so that we can identify $\ft'$ with the subspace 
$(\ker p_\ft)^\bot \subeq \g'$ of all functionals vanishing on all root spaces. 
We write $p_{\ft'} \: \g' \to \ft'$ for the restriction map.

\nin (c) If $\sigma \: \g_\C \to \g_\C$ denotes the complex conjugation with 
respect to $\g$, we write $x^* := -\sigma(x)$ for $x \in \g_\C$, so that 
$\g = \{ x \in \g_\C \: x^* = -x\}$. 
We then have $x_\alpha^* \in \g_\C^{-\alpha}$ for 
$x_\alpha \in \g_\C^\alpha$. 
\end{Definition}

\begin{Remark} 
If $\g$ carries a positive definite $e^{\ad \ft}$-invariant form $\kappa$ and 
$\ft$ is complete with respect to the induced scalar product, then 
$\g = \ft \oplus \ft^\bot$ (even if $\g$ itself 
is not complete). Therefore we have a natural orthogonal projection 
$p_\ft \: \g \to \ft$ and for every open or closed 
invariant subset $\Omega \subeq \g$, we have 
 $p_\ft(\Omega) = \Omega \cap \ft$ by Proposition~\ref{prop:2.16}(ii). 
The completeness requirement for $\ft$ is automatic if $\dim \ft < \infty$, 
but there are also interesting situations where it is satisfied for 
infinite-dimensional $\ft$ but $\g$ itself is not complete. 
Typical examples are the twisted loop algebras with values in Hilbert--Lie 
algebras (see Subsection~\ref{subsec:twistloop}). 
\end{Remark}

In the present context it is of particular interest to understand adjoint 
orbits of elements $x \in \ft$ and coadjoint orbits of 
elements $\lambda \in \ft' \cong  [\ft,\g]^\bot$. 
This can be done in terms of their projections onto $\ft$, resp., $\ft'$. 

The finite-dimensional 
Lie subalgebras arising from root vectors come in four types 
(\cite[Lemma~C.2]{Ne10}): 

\begin{Lemma}
  \label{lem:e.1} For $0 \not= x_{\alpha} \in \g_\C^{\alpha}$ the subalgebra 
$\g_\C(x_\alpha) := \Spann\{x_\alpha, x_{\alpha}^*,
[x_\alpha, x_{\alpha}^*]\}$
is $\sigma$-invariant and of one of the following types: 
\begin{itemize}
\setlength\itemsep{0em}
\item[\rm(A)] The abelian type: $[x_\alpha, x_{\alpha}^*] = 0$, i.e., 
$\g_\C(x_\alpha)$ is two-dimensional abelian. 
\item[\rm(N)] The nilpotent type: $[x_\alpha, x_{\alpha}^*] \not= 0$ 
and $\alpha([x_\alpha, x_{\alpha}^*]) = 0$, i.e., 
$\g_\C(x_\alpha)$ is a three-dimensional Heisenberg algebra. 
\item[\rm(S)] The simple type: $\alpha([x_\alpha, x_{\alpha}^*]) \not= 0$,
i.e., $\g_\C(x_\alpha) \cong \fsl_2(\C)$. In this case we distinguish the
two cases: 
\item[\rm(CS)] $\alpha([x_\alpha, x_{\alpha}^*]) > 0$, i.e., 
$\g_\C(x_\alpha) \cap \g \cong \su_2(\C)$, and 
\item[\rm(NS)] $\alpha([x_\alpha, x_{\alpha}^*]) < 0$, i.e., 
$\g_\C(x_\alpha) \cap \g \cong \su_{1,1}(\C) \cong \fsl_2(\R)$. 
\end{itemize}
\end{Lemma} 

\begin{Definition} \label{def:weylgrp} 
We call a root $\alpha \in \Delta$ {\it compact} 
if there exists an element $x_\alpha \in \g_\C^\alpha$ 
with $\alpha([x_\alpha, x_{\alpha}^*]) > 0$ 
and if the corresponding $\fsl_2(\C)$-subalgebra 
$\g_\C(x_\alpha)$ acts in a locally finite way on $\g_\C^{\rm alg}$. 
Then $\dim \g_\C^\alpha = 1$ (\cite[Prop.~I.6]{Ne00}) 
and there exists a unique element 
$\alpha^\vee \in i \ft \cap [\g_\C^\alpha, \g_\C^{-\alpha}]$ 
with $\alpha(\alpha^\vee) = 2$. We write 
$\Delta_k \subeq \Delta$ for the subset of compact roots. 
The linear endomorphism 
\[ r_\alpha \: \ft \to \ft,\quad  r_\alpha(x) 
:= x - \alpha(x)  \alpha^\vee 
= x  + (i \alpha)(x) i\alpha^\vee \] 
is called the corresponding reflection and 
\[ \cW := \cW(\fg,\ft) := \la r_\alpha \: \alpha \in \Delta_k \ra 
\subeq \GL(\ft) \] 
is called the {\it Weyl group}. 
\end{Definition}

The following proposition 
provides useful information for the 
analysis of invariant cones and orbit projections. 

\begin{Proposition}
  \label{prop:orb-pro} 
For $x \in \ft$, $x_\alpha \in \g_\C^\alpha$ and 
$\lambda \in \ft'$,  the following assertions hold: 
\begin{description}
\item[\rm(i)] 
$p_\ft(e^{\R\ad (x_\alpha - x_\alpha^*)}x) 
= x  + 
\begin{cases}
\R^+ \alpha(x) [x_\alpha^*, x_\alpha]  & \text{ for }\  
\alpha([x_\alpha, x_\alpha^*]) \leq 0 \\ 
[-1,0] \alpha(x) \alpha^\vee 
& \text{ for }\ \alpha([x_\alpha, x_\alpha^*]) > 0. 
\end{cases}$
\item[\rm(ii)] 
$p_{\ft'}(e^{\R\ad^* (x_\alpha - x_\alpha^*)}\lambda) 
= \lambda  + 
\begin{cases}
\R^+ \lambda([x_\alpha^*, x_\alpha]) \alpha   & \text{ for }\  
\alpha([x_\alpha, x_\alpha^*]) \leq 0 \\ 
[-1,0] \lambda(\alpha^\vee) \alpha 
& \text{ for }\ \alpha([x_\alpha, x_\alpha^*]) > 0. 
\end{cases}$
\end{description}
\end{Proposition}

\begin{proof} (i) is an immediate consequence of \cite[Lemma~VII.2.9]{Ne99}, 
and (ii) follows from (i) and $p_{\ft'}(\lambda)(e^{\ad y}x)= \lambda(p_\ft(e^{\ad y}x)).$
\end{proof}

Proposition~\ref{prop:orb-pro} easily shows that 
root vectors of type (A) 
lead to degeneracies with respect to invariant convex sets: 

\begin{Proposition} \label{prop:3.8} Let $0 \not= x_\alpha \in \g_\C^\alpha$ 
with $[x_\alpha, x_\alpha^*] =0$. 
\begin{itemize}
\setlength\itemsep{0em}
\item[\rm(i)]  If $C \subeq  \g$ is an open invariant convex subset 
intersecting $\ft$, then $i(x_\alpha + x_\alpha^*) \in \lim(C)$. 
\item[\rm(ii)] If $\cO_\lambda \subeq \g'$ is semi-equicontinuous 
with $B(\cO_\lambda)^0 \cap \ft \not= \eset$, 
then $i(x_\alpha + x_\alpha^*) \in \cO_\lambda^\bot$. 
\end{itemize}
\end{Proposition} 

Note that (ii) implies that semi-equicontinuous coadjoint 
orbits lie in a proper weak-$*$-closed hyperplane. 

\begin{proof} (i) Let $x \in C \cap \ft$. Then $y := x_\alpha- x_\alpha^*$ and 
$z := i(x_\alpha + x_\alpha^*)$ satisfy 
$(\ad y)^2 x  =0$, so that 
\[ e^{\R \ad y} x = x + \R [x,y] = x + \R i \alpha(x) z.\] 
Since $\alpha$ does not vanish on the open subset $C \cap \ft$, it follows that 
$\pm z \in \lim(C)$. 

\nin (ii) Since $\ft \cap B(\cO_\lambda)^0$ is open, it contains 
an element $x$ with $\alpha(x) \not=0$.  
For $y,z$ as above and $\mu \in \cO_\lambda$, we now have 
\[ \mu(e^{\R \ad y}x) 
= \mu(x)  + \R i \alpha(x) \mu(z).\] 
As $\cO_\lambda(x)$ is bounded from below, it follows that $\mu(z) =0$, 
and thus $z \in \cO_\lambda^\bot$. 
\end{proof}

To exclude the pathologies described in the preceding proposition, 
we introduce the following concept: 
\begin{Definition} \label{def:cone-pot} 
(a) We say that $\g$ has {\it cone potential} if type (A) in 
Lemma~\ref{lem:e.1} does not occur, i.e., $0 \not= x_\alpha\in \g_\C^\alpha$ 
implies $[x_\alpha, x_\alpha^*] \not=0$ 
(\cite[Def.~VII.2.22]{Ne99}). 

\nin (b) We call a root $\alpha \in \Delta$ {\it non-compact} if 
there exists a non-zero $x_\alpha \in \g_\C^\alpha$ with 
$\alpha([x_\alpha, x_\alpha^*]) \leq 0$. We write 
$\Delta_p\subeq \Delta$ for the set of non-compact roots. 
\end{Definition}

The following theorem shows that, in central aspects, the 
action of the exponential $T = \exp \ft$ of an elliptic Cartan subalgebra $\ft$ 
behaves very much like a compact group. We may therefore expect that
much of the finite-dimensional machinery developed in \cite{Ne99} 
carries over.  

\begin{Theorem}
  \label{thm:redux} 
Let $\ft \subeq \g$  be an elliptic Cartan subalgebra for which 
there exists a continuous $\ft$-equivariant 
linear projection $p_\ft \: \g \to \ft$. \\[-9mm]
\begin{itemize}
\setlength\itemsep{0em}
\item[\rm(i)] If $C \subeq \g$ is open invariant and convex, 
then $p_\ft(C) = C \cap \ft$. 
\item[\rm(ii)] If $C \subeq \g'$ is semi-equicontinuous weak-$*$-closed 
and convex, then $p_{\ft'}(C) = C \cap \ft'$. 
\end{itemize}
\end{Theorem}

\begin{proof} (i) Let $T := \exp(\ft) \subeq G$ be the connected 
subgroup corresponding to $\ft$. 
If $x \in \g^{\rm alg}$, then the orbit $\Ad(T)x$ is finite-dimensional 
and its closure coincides with the orbit of a compact 
torus group~$T_1$. Hence Proposition~\ref{prop:project}, applied to the $T$-action 
on finite-dimensional  invariant subspaces, implies that 
$p_\ft(C \cap \g^{\rm alg}) \subeq C$. 
As $C \cap \g^{\rm alg}$ is dense in $C$, 
we obtain $p_\ft(C) \subeq \oline{C \cap \ft}$ from the continuity of~$p_\ft$. 
Since $p_\ft \: \g \to \ft$ is a projection, 
the summation map $\ft \oplus \ker(p_\ft) \to \g$ 
is a topological isomorphism, so that $p_\ft$ is also  open. 
Hence $p_\ft(C)$ is an open subset of $\oline{C \cap \ft}$. 
Lemma~\ref{lem:bou}(ii) now show that $p_\ft(C) \subeq C \cap \ft$. 
The inclusion $C \cap \ft \subeq p_\ft(C)$ is trivial. 

\nin (ii) Applying (i) to the open invariant cone $B(C)^0$, we find 
an element $x_0 \in B(C)^0 \cap \ft$. Then all subsets 
\[ C_r := \{ \lambda \in C \: \lambda(x_0) \leq r \}, \quad r \in \R,  \] 
are weak-$*$-compact (\cite[Prop.~6.13]{Ne08}) and $\Ad(T)$-invariant. 

Let $R \: \g' \to (\g^{\rm alg})'$ denote the restriction map. 
As $C_r$ is weak-$*$-compact, its image $R(C_r)$ is also 
weak-$*$-compact. The action of $\Ad(T)$ on $E := \g^{\rm alg}$ 
satisfies the assumptions of Proposition~\ref{prop:2.17}, 
so that Kakutani's Theorem implies that, for every 
$\lambda \in C_r$, the subset 
$\oline\conv(\Ad(T)'\lambda)$ contains a fixed point, hence 
an element of $\ft'$. This implies that 
$p_{\ft'}(\lambda) = {\lambda\res_{\ft} \in C}$. 
\end{proof}

Combining Theorem~\ref{thm:redux}  with 
Proposition~\ref{prop:orb-pro}, we obtain the following two corollaries: 

\begin{Corollary} \label{cor:3.7} 
If $C \subeq \g$ is an open invariant convex subset,  
then $C_\ft := C \cap \ft = p_\ft(C)$ is $\cW$-invariant and satisfies 
\[ i \alpha(C_\ft) \cdot i [x_\alpha, x_\alpha^*] \subeq  \lim(C) 
\quad \mbox{ for } \quad \alpha([x_\alpha, x_\alpha^*]) \leq 0.\] 
In particular, if $\lim(C)$ is pointed and $\g$ has cone potential, then  
$i\Delta_p \subeq C_\ft^\star \cup - C_\ft^\star$. 
\end{Corollary}

\begin{Corollary} \label{cor:3.8} 
If $C \subeq \g'$ is a weak-$*$-closed convex subset,  
then $C_{\ft'} := C \cap \ft' = p_{\ft'}(C)$ is $\cW$-invariant with 
\[ \la C, i [x_\alpha, x_\alpha^*] \ra \cdot i \alpha \subeq  \lim(C_{\ft'}) 
\quad \mbox{ for } \quad \alpha([x_\alpha, x_\alpha^*]) \leq 0.\] 
If $\lim(C_{\ft'})$ is pointed, then 
$i[x_\alpha, x_\alpha^*] \in C^\star \cup - C^\star$, 
and $i[x_\alpha, x_\alpha^*] \in C^\star$ implies $i\alpha \in \lim(C_{\ft'})$. 
\end{Corollary}

\begin{Remark} \label{rem:maxcone} Assume that $\g$ has cone potential. 

\nin (a) Any open invariant convex cone $C \subeq \g$ intersects $\ft$ by 
Theorem~\ref{thm:redux}. If 
$\oline C$ is pointed, then Corollary~\ref{cor:3.7} implies that 
$C$ specifies a $\cW$-invariant {\it positive system of non-compact 
roots} 
\[ \Delta_p^+ := \{ \alpha \in \Delta_p \: i \alpha \in C_\ft^\star \} 
= \Delta_p \cap - i C_\ft^\star
\quad \mbox{ with }\quad 
\Delta_p = \Delta_p^+ \dot\cup -\Delta_p^+.\] 
This in turn defines two cones: 
\begin{equation}
  \label{eq:maxcon}
 C_{\rm min}(\Delta_p^+) 
:= \oline\cone(\{ i [x_\alpha, x_\alpha^*] \: 
x_\alpha \in \g_\C^\alpha, \alpha \in \Delta_p^+\}), \qquad 
C_{\rm max}(\Delta_p^+) := (i\Delta_p^+)^\star.
\end{equation} 
Corollary~\ref{cor:3.7} now entails that 
\begin{equation}
  \label{eq:sandwich}
C_{\rm min}(\Delta_p^+)  \subeq \oline{C_\ft} 
\subeq C_{\rm max}(\Delta_p^+).
\end{equation}

\nin (b) Suppose that $\lambda \in \ft'$ is such that 
$\cO_\lambda =\Ad^*(G)\lambda \subeq \g'$ is semi-equicontinuous 
and put $C_\lambda := \oline\conv(\cO_\lambda)$ (weak-$*$ closure). 
We further assume that $\cO_\lambda^\bot = \{0\}$. 

Then $B(\cO_\lambda)^0 = B(C_\lambda)^0 \subeq \g$ is an open invariant convex cone 
with $B(\cO_\lambda)^\star = \lim(C_\lambda)$ 
(Lemma~\ref{lem:limcone}). As in (a), we choose the positive system $\Delta_p^+$ 
such that 
\[i \Delta_p^+ 
= (i \Delta_p) \cap B(\cO_\lambda)^\star 
= (i \Delta_p) \cap \lim(C_\lambda)\] 
and obtain the relation 
\begin{equation}
  \label{eq:convrel1}
\conv(\cW\lambda) + \cone(i\Delta_p^+) 
\subeq p_{\ft'}(\oline\conv(\cO_\lambda)).
\end{equation}
In Remark~\ref{rem:convth1} below we shall see  
how representation theoretic arguments can be 
used to obtain the converse inclusion. In any case we see that 
\begin{equation}
  \label{eq:bolambda}
 B(\cO_\lambda) \cap \ft 
\subeq B(\cW \lambda) \cap (i\Delta_p^+)^\star 
= B(\cW \lambda) \cap C_{\rm max}(\Delta_p^+).
\end{equation}
In particular, the cone 
$C_{\rm max}(\Delta_p^+)$ needs to have interior points if $\cO_\lambda$ 
is semi-equicontinuous, 
and the Weyl group orbit $\cW\lambda \subeq \ft$ has to be 
semi-equicontinuous as well. 
\end{Remark}

\begin{Remark} Let $\ft \subeq \g$ be an elliptic Cartan subalgebra and 
$T := \exp \ft$. 
For $\alpha \in \Delta$, we consider the character 
$\chi_\alpha(\exp x) := e^{i\alpha(x)}$ of $T$ and 
the subgroup $\cQ \subeq \hat T := \Hom(T,\T) \into \ft'$, 
generated by these characters. 
We consider $\cQ$ as a discrete group, so that 
\[ S := \hat \cQ = \Hom(\cQ,\T) \] 
is a compact abelian group. The adjoint representation defines a 
homomorphism 
\[ A \: T\to S, \quad A(\exp x)_\beta := e^{i\beta(x)} \] 
whose range separates the character group $\hat S \cong \cQ$ 
and therefore is dense in $S$ (\cite[Thm.~7.64]{HM13}). 

Each $s = (s_\beta)_{\beta \in \Delta} \in S$ defines an automorphism 
$\psi(s) \in \Aut(\g^{\rm alg}) \cong \Aut(\g^{\rm alg}_\C)^\sigma$ by 
$\psi(s)x_\alpha = \chi_\alpha(s) x_\alpha$ for $x_\alpha \in \g_\C$

\nin (a) If $\g = \g^{\rm alg}$, then the topology of pointwise 
convergence coincides on $A(T)$ with the product topology induced 
from~$S$, so that $S$ is the closure of $A(T)$ in the topology of pointwise 
convergence. 
This group preserves 
all weak-$*$-closed $\Ad(T)$-invariant convex subsets of 
the algebraic dual $(\g^{\rm alg})^*$ (Lemma~\ref{lem:endo-strong-closed}). 

\nin (b) If $\g$ is strictly larger than $\g^{\rm alg}$, 
then it is not clear that all elements of $S$ act by automorphisms on $\g$. 
However, in all concrete examples we are aware of, this is the case. 
For example if $\g = \fu_p(\cH)$ and $(e_j)_{j \in J}$ is an orthonormal 
basis, then $\ft \cong i \ell^p(J,\R)$ is an elliptic 
Cartan subalgebra of $\g$, and $\g^{\rm alg} = \fu(J,\C)$. 
Here $S \cong \T^J/\T$, where 
$\T \subeq \T^J$ represents the constant functions $J \to \T$. 
The compact group $S$ acts continuously on $\fu_p(\cH)$ 
(cf.\ Remark~\ref{rem:extend-action}). 
\end{Remark}

\subsection{The connection with unitary representations} 
\label{subsec:3.2}

Unitary representations can provide information 
on invariant convex subsets through their momentum sets and support functionals. 

\begin{Proposition} If $(\pi, \cH)$ is a unitary representation of $G$, 
then 
\[ C_\pi := \{ x \in \g \: -i \partial \pi(x) \geq 0\} = I_\pi^\star \] 
is a closed convex invariant cone in $\g$. 
If there exists an element $\bc \in \g$ with 
$\partial \pi(\bc) = i \one$ and $\pi$ is semibounded, then $C_\pi$ has interior points. 
\end{Proposition}

\begin{proof} The definition of the momentum set (Subsection~\ref{subsec:3.0}) 
implies that 
$C_\pi = I_\pi^\star$, so that the $\Ad(G)$-invariance of $I_\pi$ implies that 
the closed convex cone $C_\pi$ is also invariant. 

For the second assertion, suppose that $s_\pi \leq m$ on 
an open subset $U \subeq \g$. 
Then $-m \bc + U \subeq C_\pi$ shows that $C_\pi$ has interior points. 
\end{proof}

\begin{Remark} (Connections with Lie supergroups) 
For a unitary representation of a (possible infinite-dimensional) Lie superalgebra 
$\g = \g_{\oline 0} \oplus \g_{\oline 1}$, 
we have ${-i \dd\pi(x)^2 \geq 0}$ for odd elements $x \in \g_{\oline 1}$ 
(\cite{NS11}). 
Therefore 
\[ C_\pi := \{ x \in \g_{\oline 0} \:  -i \partial\pi(x) \geq 0\} \]
is a closed convex invariant cone in the even part $\g_{\oline 0}$ which contains 
all brackets $[x,x]$, $x \in \g_{\oline 1}$. In many situations this cone has interior 
points, so that a classification of open invariant cones in $\g_{\oline 0}$ 
is of central importance to understand the unitary representations of 
the Lie superalgebra~$\g$ (\cite{AN15}, \cite{NY18}). 
\end{Remark}

\begin{Remark} \label{rem:convth1}
(a) Let $\iota \: H \to G$ be an injective morphism of connected Lie 
groups and identify $\fh$ by the tangent map 
$\L(\iota)$ with a Lie subalgebra of $\g$. 
If $(\pi, \cH)$ is  a semibounded unitary representation for 
which there exists an $x_0 \in \fh \cap B(I_\pi)^0$, 
then $\pi\res_H$ is a semibounded representation of $H$ with momentum set 
\[ I_{\pi\res_H} = I_\pi\res_{\fh}.\]
If the spectrum of the restriction $\pi\res_H$ 
is known, this representation theoretic information 
can be used to derive information on the momentum set $I_\pi$ 
and on coadjoint $G$-orbits (see~\cite{Ne12, Ne14} for examples).  
In particular, we obtain for every $\lambda \in I_\pi$ the relation 
\begin{equation}
  \label{eq:momsetrest}
\cO_\lambda\res_{\fh} \subeq I_{\pi\res_H}.
\end{equation}

\nin (b) If $H$ is abelian, then every semibounded representation $(\pi,\cH)$ of $H$ 
is given by a spectral measure $P$ on $\fh'$ (\cite[\S 7]{Ne08}) and 
$I_\pi = \oline\conv(\supp(P))$, so that, in the context of (a), 
\begin{equation}
  \label{eq:pro-inc}
\cO_\lambda\res_\fh \subeq \oline\conv(\supp(\pi\res_H)).
\end{equation}
This information is particularly useful if one has explicit information 
on the spectrum of $\pi\res_H$. 
\end{Remark} 

\begin{Examples} \label{ex:3.18a}
An important special case arises for 
$H = T = \exp(\ft)$, $\ft$ an elliptic Cartan subalgebra, 
if $\pi\res_T$ is a direct sum of weight spaces 
\[ \cH_{i\mu}(\ft) = \{ v \in \cH \: (\forall x \in \ft)\, 
\dd\pi(x)v = i \mu(x) v \}, \qquad 
\mu \in \cP_\pi := \{ \nu \in \ft' \: \cH_{i\nu}(\ft) \not=\{0\}\}.\] 
If 
\begin{equation}
  \label{eq:wconvlambda}
\oline\conv(\supp(\pi\res_T)) = \oline\conv(\cW \lambda) 
\quad \mbox{ for a weight } \quad 
\lambda \in \ft',
\end{equation} 
then we obtain from \eqref{eq:pro-inc} 
\begin{equation}
  \label{eq:pro-incb}
 \cO_\lambda\res_\ft \subeq \oline\conv(\cW \lambda).
\end{equation}

This situation arises in particular if: 
\begin{itemize}
\setlength\itemsep{0em}
\item[\rm(a)] $\Delta = \Delta_k$, $\g_\C = \g_\C^{\rm alg}$ is a locally finite Lie algebra, 
and $\pi_\lambda$ is a unitary representation with extremal weight 
$i\lambda$ (\cite{Ne98}, \cite[\S 3]{Ne04}). 
\item[\rm(b)] $\g = \g^{\rm alg}$ is a unitary real form of the Kac--Moody Lie algebra $\g_\C$, 
i.e., all root vectors are of type (N) or (CS), 
and $\pi_\lambda$ is a unitary representation with 
highest weight $i\lambda$ with respect to some positive system $\Delta^+$ 
(\cite[Ch.~11]{Ka90}). These representations exists if and only if 
$i\lambda(\check \alpha) \in \N_0$ for all simple roots $\alpha \in \Delta_+$. 
\item[\rm(c)] $\g = \g^{\rm alg}$ is a unitary form of a locally affine 
Lie algebra $\g_\C$ and $\pi_\lambda$ is a unitary representation with 
extremal weight $i\lambda$ not vanishing in the central element~$\bc$ 
(see \cite[Thms.~4.10/11]{Ne10b} for details). 
These Lie algebras arise in particular 
as dense subalgebras of twisted loop algebras 
$\cL_\phi^\tau(\fk_\C)$ (see Subsection~\ref{subsec:twistloop}). 
\end{itemize}
In all these cases, combining \eqref{eq:pro-incb} with 
the trivial inclusion $\cW\lambda \subeq \cO_\lambda \cap \ft' \subeq  
p_{\ft'}(\cO_\lambda)$ leads to 
\begin{equation}
  \label{eq:repconvtheo}
p_{\ft'}(\oline\conv(\cO_\lambda)) = \oline\conv(\cW\lambda) 
\end{equation}
for the extremal weight~$i\lambda$ of the unitary representation~$\pi$.
\end{Examples}

\begin{Example} \label{ex:3.18} 
If $\g = \g^{\rm alg}$ is locally finite and not all roots 
are compact, then the situation is more complicated. 
For the corresponding unitary highest weight representations $\pi_\lambda$ 
(see \cite{Ne01} for a classification) we obtain under the assumption 
that $\pi_\lambda$ has discrete kernel that 
\begin{equation}
  \label{eq:contheoweights}
\conv(\cP_{\pi_\lambda}) = \conv(\cW\lambda) +i \cone(\Delta_p^+),
\end{equation}
which leads to 
\begin{equation}
  \label{eq:repconvtheo2}
p_{\ft'}(\oline\conv(\cO_\lambda)) = \oline{\conv(\cW\lambda) +i \cone(\Delta_p^+)}.
\end{equation}
\end{Example}

\begin{Problem}
We define a quasi-order on $\g'$ by 
\[ \mu \prec \lambda \quad \mbox{ if } \quad 
\oline\conv(\cO_\mu) \subeq    \oline\conv(\cO_\lambda)\] 
and say that  $\lambda$ {\it majorizes} $\mu$ if $\mu \prec \lambda$. \\[-7mm]
\begin{itemize}
\item[\rm(P6)] Let $\ft \subeq \g$ be an elliptic Cartan subalgebra. 
Give an explicit description of the quasi-order $\prec$ on $\ft' \subeq \g'$. 
\end{itemize}
If all roots are compact, this this is closely related to the order relation 
\[ \mu \prec_\cW \lambda \quad \mbox{ if } \quad 
\oline\conv(\cW\mu) \subeq    \oline\conv(\cW\lambda).\] 
See Example~\ref{ex:schur}, Remark~\ref{rem:5.4}, 
Remark~\ref{rem:extend-action} and Problem~\ref{prob:2} for related issues. 

\end{Problem}

\section{Finite-dimensional Lie algebras} 
\label{sec:4}

In this section we recall the key results on invariant convex cones in 
finite-dimensional Lie algebras.

\subsection{Compact Lie algebras} 
\label{subsec:4.1}

The best behaved examples 
are {\it compact} Lie algebras $\g$, i.e., Lie algebras of a compact Lie group~$G$. 
Then $\g$ is finite-dimensional and carries an $\Ad(G)$-invariant scalar product 
$\kappa$. The map $\g \to \g', x \mapsto x^* := \kappa(x,\cdot)$ 
is a $G$-equivariant equivalence between 
adjoint and coadjoint representation. Moreover, 
every $x$ is an extreme point of $\conv(\cO_x)$ because the orbit lies 
in a euclidean sphere. 

We collect the main results in the following theorem: 
\begin{Theorem} \label{thm:compact-lie} 
Let $\g$ be a compact Lie algebra, $\ft \subeq \g$ be maximal abelian and 
$p_\ft \:  \g \to \ft$ be the projection along $[\ft,\g]$. 
Then the following assertions hold: \\[-7mm]
\begin{itemize}
\setlength\itemsep{0em}
\item[\rm(i)] Every adjoint orbit $\cO_x = \Ad(G)x \subeq \g$ intersects 
$\ft$ in an orbit of the Weyl group $\cW$, which is a finite Coxeter group. 
\item[\rm(ii)] $p_\ft(\cO_x) = \conv(\cW x)$ for $x \in \ft$  
{\rm(Kostant's Convexity Theorem)}. 
\item[\rm(iii)] The map $C \mapsto C \cap \ft$ defines a bijection 
from $\Ad(G)$-invariant open/closed convex subsets of $\g$ to 
$\cW$-invariant open/closed convex subsets of $\ft$. 
\item[\rm(iv)] For $x \in \ft$, the cone $L_x(\conv \cW x)$ 
{\rm(Definition~\ref{def:subtang})} 
is generated by the elements $-\alpha(x) \alpha^\vee$, 
$\alpha \in \Delta(\g,\ft)$.  
\item[\rm(v)] 
If $\lambda([x,y]) = 0$ and 
$\lambda([y,[y,x]])\geq 0$ hold for all $y \in \g$, 
then $\lambda(x) = \min \lambda(\cO_x)$. 
\end{itemize}
\end{Theorem}

\begin{proof} (i) \cite[Thm.~12.2.2, Lemma 12.2.16]{HN12}; 
(ii) is contained in \cite{Ko73}, 
and (iii) in \cite[Thm.~III.17]{Ne96}. 

\nin (iv) follows from the fact that $\cW$ is a finite Coxeter group 
with $\cT = \ft$ and \cite[Thm.~2.7]{HN14} (see also 
Theorem~\ref{thm:coxconvtheo}).

\nin (v) For $z^* = \kappa(z,\cdot)$, 
the first part of condition \eqref{eq:nec-cond} 
leads to $[x,z] = 0$ so that both lie in a maximal abelian 
subalgebra $\ft$, and the second part of \eqref{eq:nec-cond} 
can be analyzed in terms of the root decomposition 
(Definition~\ref{def:basic}) which implies~(v). 
\end{proof}

Assertions (iii) and (iv) in Theorem~\ref{thm:compact-lie} 
are the key to a complete classification of 
convex $\cW$-invariant subsets of $\ft$ and $\Ad(G)$-invariant convex subsets 
of $\g$, carried out in \cite{Ne96, Ne99}. 

\begin{Example} \label{ex:schur-horn}
For $G = \U_n(\C)$ and $\g = \fu_n(\C)$, we chose the subspace 
$\ft\cong i \R^n$ of diagonal matrices as an elliptic Cartan subalgebra 
and 
\[ p_\ft \:\fu_n(\C) \to i\R^n, \quad 
p_\ft(x) = \diag(x_{11}, \ldots, x_{nn}) \] 
is the corresponding projection. 
If $x = \diag(\lambda_1, \ldots, \lambda_n)\in \ft$, 
then the adjoint orbit 
 $\cO_x := \{ gxg^{-1} \: g \in \U_n(\C)\}$ is the set of all skew-hermitian 
matrices with the eigenvalues $\{\lambda_1,\ldots, \lambda_n\}$, and 
the Schur--Horn Convexity Theorem asserts that 
\begin{equation}
  \label{eq:schurhorn}
p_\ft(\cO_{x}) 
= \conv(\{(\lambda_{\sigma(1)}, \ldots, \lambda_{\sigma(n)}) \: \sigma \in S_n\}) 
= \conv(S_n.\lambda)  
\end{equation}
is the convex hull of the orbit of $x$ under the symmetric group $S_n\cong \cW$. 
\end{Example}

\begin{Remark} \label{rem:center-compact}
For a compact Lie algebra $\g$, 
Proposition~\ref{prop:project} implies in particular that 
every open invariant cone $\Omega \subeq \g$ intersects 
the center $\fz(\g)$. In particular, $\g$ contains no non-trivial open 
invariant cones if $\fz(\g) = \{0\}$, i.e., if $\g$ is semisimple. 
If, conversely, $0 \not= z \in \fz(\g)$, then there exist open 
convex $\Ad(G)$-invariant neighborhoods $U$ of $z$ not containing $0$, 
and then $\Omega := \bigcup_{t > 0} t U$ is a proper open invariant cone 
(cf.~Lemma~\ref{lem:a.6}(b)). 
\end{Remark}

\subsection{Non-compact finite-dimensional Lie algebras}

In this subsection we assume that $\g$ is finite-dimensional 
and that $\ft \subeq \g$ is a compactly embedded (=elliptic) 
Cartan subalgebra. The existence of $\ft$ 
is ensured by the assumption that the subset 
$\g'_{\rm seq}\subeq \g'$ of functionals $\lambda$ with semi-equicontinuous 
orbit $\cO_\lambda$ separates the points of $\g$ 
(\cite[Thm.~VII.3.28]{Ne99}). 
If $G$ is a connected Lie group with $\g = \L(G)$, then we 
write $T = \exp(\ft)$ for the subgroup corresponding to~$\ft$.

For $\lambda \in \ft'$, the orbit $\cO_\lambda \subeq \g'$ 
is semi-equicontinuous if and only if it is admissible 
in the sense of \cite{Ne99} because it is automatically closed 
(\cite[Thm.~VIII.1.8]{Ne99}). The Lie algebra 
$\g$ is said to be {\it quasihermitian} if 
$\fz_\g(\fz(\fk)) = \fk$ holds for the compact Lie subalgebra 
$\fk$ with $\fk_\C = \ft_\C + \sum_{\alpha \in \Delta_k} \g_\C^\alpha$. 
This condition is always satisfied if $\g'$ is generated by  
semi-equicontinuous orbits (\cite[Thm.~VII.3.28]{Ne99}).

\begin{Theorem} {\rm(Convexity theorem for semi-equicont.\ orbits---
finite-dim.\  case)}  
Suppose that $\g$ is quasihermitian. 
For  $\lambda \in \ft'$, the coadjoint orbit 
$\cO_\lambda$ is semi-equicontinuous if and only if 
there exists a $\cW$-invariant positive system $\Delta_p^+$ 
with $\lambda \in C_{\rm min}(\Delta_p^+)^\star$, and then 
$p_{\ft'}(\cO_\lambda)$ is convex with 
\begin{equation}
  \label{eq:findimconvtheo}
 p_{\ft'}(\cO_\lambda) \subeq \conv(\cW\lambda) + \cone(i\Delta_p^+).
\end{equation}
If, in addition, $\cO_\lambda^\bot=\{0\}$, then equality holds in \eqref{eq:findimconvtheo}. 
\end{Theorem}

\begin{proof} The first part and \eqref{eq:findimconvtheo} follow from 
\cite[Def.~VII.2.6, Thm.~VIII.1.19]{Ne99}. 
If $\cO_\lambda$ spans~$\g'$, then 
\eqref{eq:convrel1} in Remark~\ref{rem:maxcone}(b) 
implies equality in \eqref{eq:findimconvtheo}. 
\end{proof}

For adjoint orbits we have (\cite[Thm.~VIII.1.36]{Ne99}): 
\begin{Theorem} {\rm(Convexity theorem for adjoint orbits---
finite-dimensional case)}
Suppose that $\g$ is finite-dimensional and $\ft \subeq \g$ is an 
elliptic Cartan subalgebra. 
If $\Delta_p^+$ is a $\cW$-invariant positive system of non-compact roots, 
then 
\[ p_\ft(\cO_x) \subeq \conv(\cW x) + C_{\rm min}(\Delta_p^+) 
\quad \mbox{ for } \quad x \in C_{\rm max}.\] 
\end{Theorem}

These convexity theorems 
lead to complete information on invariant convex subsets 
(cf.~Remark~\ref{rem:convth1}(c)). 
The main tool is the action of the finite Coxeter group $\cW = \cW(\g,\ft)$, 
and the determination of the cones 
$C_{\rm max}(\Delta_p^+)$ for the different $\cW$-invariant positive 
systems in $\Delta_p$. 

\begin{Remark} (a) If $\lambda \in \g'$ is contained in the algebraic interior of an 
invariant semi-equicontinuous subset $C \subeq \g'$, i.e., 
$C - \lambda$ is absorbing in the linear space $C - C$  
($\lambda$ is strictly admissible in the terminology of \cite{Ne99}), 
then $\cO_\lambda$ intersects $\ft'$ (\cite[Thm.~VIII.1.8, Lemma~VIII.1.27]{Ne99}), 
$\conv(\cO_\lambda)$ 
is closed and $\cO_\lambda = \Ext(\conv(\cO_\lambda))$ is the set of its 
extreme points (\cite[Prop.~VIII.1.30]{Ne99}). 
Not every coadjoint orbit in $\g'_{\rm seq}$ intersects 
$\ft' \cong [\ft,\g]^\bot$, as the nilpotent orbits in $\g = \fsl_2(\R)$ show. 

\nin (b) For every 
irreducible semibounded representation $(\pi,\cH)$, there exists a 
$\lambda \in \ft'\cong [\ft,\g]^\bot$ such that the set $\Ext(I_\pi)$ 
of extreme points of $I_\pi$ is a single coadjoint orbit $\cO_\lambda$ 
and $I_{\pi} = \conv(\cO_\lambda)$ (\cite[Thm.~X.4.1]{Ne99}). 
Further 
$\cH$ decomposes into weight spaces $\cH_{i\alpha}(T)$, $\alpha \in \ft'$. 
The corresponding weight set $\cP_\pi \subeq \ft'$ has the property that 
\[ \cO_\lambda\res_{\ft} = \conv(\cP_\pi) \quad \mbox{ and } \quad 
\lambda \in \Ext(\conv(\cP_\pi)),\] 
i.e., {\it $\lambda$ is an extremal weight} (cf.\ Example~\ref{ex:3.18}). 
\end{Remark}

In the finite-dimensional context, many convexity theorems have been proved 
by symplectic techniques using convexity properties of momentum maps. 
We do not expect this to work in the infinite-dimensional context, where 
we are mostly interested in inclusions such as in \eqref{eq:pro-incb} or 
\eqref{eq:findimconvtheo} and not in equalities (see in particular \cite{BFR93, 
BEFR97}). 
We therefore put a stronger emphasis on 
functional analytic arguments using convex sets. 

\section{Infinite-dimensional Lie algebras} 
\label{sec:5}

For infinite-dimensional Lie algebras, only very particular results concerning 
problems (P1-6) are known. They show certain common patterns, but so far no systematic 
theory has been developed to create a unifying picture. 
We now briefly discuss several classes of infinite-dimensional Lie algebras 
and their invariant convex cones. 

\subsection{Nilpotent and $2$-step solvable Lie algebras} 

If $x \in \g$ satisfies $(\ad x)^2 = 0$, then 
$\Ad^*(\exp tx)\lambda = \lambda \circ e^{-t \ad x}
= \lambda - t (\lambda \circ \ad x)$ for 
$\lambda \in \g'$ shows that 
the orbits of the corresponding one-parameter group 
in $\g'$ are either trivial of affine lines. 
If $\cO_\lambda$ is semi-equicontinuous, the orbit must be trivial 
(cf.\ Proposition~\ref{prop:3.8}). 
A closer inspection of this simple observation leads to:
\begin{Theorem} \label{thm:nilp} 
Suppose that $\g$ is either nilpotent or $2$-step solvable, i.e., 
$[\g,\g]$ is abelian. 
Then $\g'_{\rm seq} \subeq [\g,\g]^\bot$, i.e., 
all semi-equicontinuous coadjoint orbits are trivial. 
If there exists a pointed invariant cone $W \subeq \g$ with non-empty 
interior, then $\g$ is abelian. 
\end{Theorem}

\begin{proof} The first assertion follows from \cite[Thm.~1.5]{NZ13}. 
For the second assertion, we note that $W^\star \subeq \g'_{\rm seq}$, 
and since $W$ is pointed, $W^\star$ separates the points of $\g$. 
As the coadjoint action is trivial on $\g'_{\rm seq}$, the adjoint action 
on $\g$ must be trivial as well. This means that $\g$ is abelian. 
\end{proof}

\subsection{Oscillator algebras} 
\label{subsubsec:osci} 

In {\rm(\cite[Thms.~2.8, 3.2, Prop.~3.4]{NZ13})}, 
we determine all semi-equicontinuous coadjoint 
orbits in {\it oscillator algebras}, i.e.,  
double extensions  $\g = \g(V,\omega, D) = (\R \oplus_\omega V) \rtimes_D \R$ 
of an abelian Lie algebra $V$, where $(V,\omega)$ is a symplectic vector space 
and $D \in \fsp(V,\omega)$ (Definition~\ref{def:doubleext}):  

\begin{Theorem}
 \label{thm:osci} 
For an oscillator algebra, the following  are equivalent: \\[-7mm]
\begin{itemize}
\setlength\itemsep{0em}
\item[\rm(i)] $\g'_{\rm seq} \not=\eset$. 
\item[\rm(ii)] $q(x,y) := \omega(Dx,y)$ or $\omega(x,Dy)$ is positive definite, 
all functionals $i_x \omega = \omega(x,\cdot)$ are 
$q$-continuous, and the corresponding quadratic form 
$x \mapsto \|i_x\omega\|_q^2$ is continuous. 
\item[\rm(iii)] $\g$ is a Lorentzian double extension $\g(V,\kappa,D)$ 
of a locally convex 
euclidean vector space $(V,\kappa)$ on which $D$ is a $\kappa$-skew-symmetric 
derivation and $\omega(v,w) = \kappa(v,Dw)$. 
\end{itemize}
If these conditions are satisfied, then 
$\lambda = (z^*, \alpha, t^*) \in \g'_{\rm seq}$ if and only if 
$z^*\not=0$ and $\alpha\res_{D(V)}$ is $\kappa$-bounded. 
If, in addition, $D(V)$ is dense in $V$, then the $\kappa$-boundedness of 
$\alpha$ on $D(V)$ implies its $\kappa$-boundedness. 
\end{Theorem}

\pagebreak 
\subsection{Kac--Moody Lie algebras} 
\label{subsec:kacmoo}

In \cite[Thm.~2(b)]{KP84} Kac and Peterson generalized 
Kostant's Convexity 
Theorem to symmetrizable Kac--Moody algebras, 
and this implies the Atiyah--Pressley Convexity Theorem for loop groups 
\cite[Thm.~1]{AP83}. 
Concretely, let $\ft \subeq \g$ be an elliptic Cartan subalgebra 
such that $\g = \g^{\rm alg}$, 
$\g_\C$ is a symmetrizable Kac--Moody 
algebra, and all root vectors are of type (N) or (CS) (cf.\ Lemma~\ref{lem:e.1}), 
i.e., $\g$ is a unitary real form of $\g_\C$. 
For every element $x$ in the Tits cone $\cT \subeq \ft$ 
(Subsection~\ref{subsec:cox}), it asserts that 
\[ p_\ft(\cO_x)= \conv(\cW x) \quad \mbox{ for } \quad x \in \cT,\] 
where $\cW = \cW(\fg,\ft)$ is a finitely generated Coxeter group 
(Definition~\ref{def:weylgrp}, \cite[Ch.~6]{Ka90}). 
This result applies in particular to 
twisted loop algebras with finite-dimensional compact target groups 
(cf.~Subsection~\ref{subsec:twistloop}). 

\nin Since $\g$ carries a non-degenerate invariant symmetric bilinear form, 
we likewise obtain for coadjoint orbits $\cO_\lambda$, 
$\lambda \in \cT'$ (the Tits cone in $\ft'$), that 
\[ p_{\ft'}(\cO_\lambda)= \conv(\cW \lambda) \quad \mbox{ for } \quad 
\lambda  \in \cT'.\] 

We now sketch a rather general representation theoretic argument that 
provides the inclusion ``$\subeq$'' of this convexity theorem. 
If $\lambda \in \ft'$ is dominant integral, i.e., 
$i\lambda(\alpha^\vee) \in \N_0$ for all simple roots $\alpha\in \Pi$, 
then it is a highest weight of a unitary representation 
$(\pi_\lambda, \cH_\lambda)$ (see \cite{KP84}), so that 
the discussion in Example~\ref{ex:3.18a}(b) implies that 
\begin{equation}
  \label{eq:dom-conv}
 p_{\ft'}(\cO_\lambda) \subeq \oline\conv(\cW\lambda) \subeq 
\lambda + i \cone(\Pi).  
\end{equation}
Approximating general elements $\lambda \in \cT'$ by 
positive multiples of integral ones, we see that 
\eqref{eq:dom-conv} holds for every $\lambda \in (\cT')^0$. 
The Convexity Theorem for Coxeter groups, combined with 
Lemma~\ref{lem:orb-cont}, now yields 
\begin{equation}
  \label{eq:convtheokacmo}
p_{\ft'}(\cO_\lambda) 
\subeq \cT' \cap \bigcap_{w \in \cW} 
w(\lambda + i \cone(\Pi)) = \oline\conv(\cW\lambda).
\end{equation}

\subsection{Lie algebras of vector fields} 
\label{subsec:5.3}

In this  subsection we briefly comment on two results on Lie algebras 
of vector fields. 

\nin {\bf Vector fields on the circle.} 
Let $\g = \cV(\bS^1) = C^\infty(\bS^1) \partial_\theta$ 
be the Lie algebra of smooth vector fields 
on the circle $\bS^1 \cong \R/\Z$, where $\bd := \frac{\partial}
{\partial \theta}$ denotes the generator of the right rotations. Then 
$$ W := \Big\{ f \frac{\partial}{\partial \theta} \: f > 0 \Big\} $$
is an open invariant cone and $\pm W$ are the only open invariant 
cones in this Lie algebra. 
Based on this observation and a specific Convexity Theorem, 
we classify in \cite[\S 8]{Ne10} all semi-equicontinuous 
coadjoint orbits of the Virasoro algebra 
and the Lie algebra $\cV(\bS^1)$ of smooth 
vector fields on the circle. Here Lemma~\ref{lem:doubcone} 
applies with $\bd$ as above. It implies that $\pm W$ 
are the only proper open invariant cones in 
$\cV(\bS^1)$ (\cite[Thm.~8.3]{Ne10}) and $W^\star \cup - W^\star$ is 
the set of semi-equicontinuous coadjoint orbits 
(\cite[Prop.~8.4]{Ne10}). 

In the Virasoro algebra $\vir$, a central extension 
$\vir = \R \bc \oplus_\omega \cV(\bS^1)$, the situation is more 
complicated. Here the open invariant cones are classified by a sign 
and an angle between $\pi/2$ and $\pi$ 
(\cite[Thm.~8.15]{Ne10}). 

\nin {\bf The diffeomorphism group of an annulus.} 
In \cite{BFR93, BEFR97} a convexity theorem for the diffeomorphism group 
of the annulus $M := [0,1] \times \bS^1$ is obtained. 
Here $\g$ is the corresponding Lie algebra of smooth functions 
with the Poisson bracket and $\ft \subeq \g$ is the abelian subalgebra of 
radial functions $f(r,\theta) = f(r)$. 
In this context several completion procedures are required to 
identify the natural generalization of the Weyl group, which in this context 
leads to the semigroup $\oline\cW$ of measure preserving maps on $[0,1]$ 
and the closed convex hulls of its orbits in $L^1([0,1])$ 
(\cite{Br66}, \cite{Ry65, Ry67}); 
see also Problem~\ref{prob:2}). 
This exhibits an interesting analogy with the situation 
discussed in Remark~\ref{rem:extend-action} below, where 
we show that the monoid $\Isom(\cH)$ acts on orbit closures in 
the duals spaces $\fu_p(\cH)'$ and $\fu(\cH)'$.

\subsection{Unitary Lie algebras} 
\label{subsec:unitliealg}

{\bf The full unitary Lie algebra of a Hilbert space.} 
Let $\cH$ be a complex Hilbert space 
with an orthonormal basis $(e_j)_{j \in J}$, so that 
$\cH \cong \ell^2(J,\C)$.  
For the full unitary group $G = \U(\cH)$, 
a generalization of the 
Schur--Horn Theorem (Example~\ref{ex:schur-horn})  
was obtained by A.~Neumann \cite{Neu99, Neu02}. 
We write $\ft \cong i \ell^\infty(J,\R) \subeq \g = \fu(\cH)$ 
for the subalgebra of diagonal operators with respect to the orthonormal 
basis. Although it is maximal abelian, it is not an elliptic 
Cartan subalgebra in the sense of Definition~\ref{def:basic} 
because the subalgebra $\g^{\rm alg}$ is not dense in $\g$ 
with respect to the operator norm. Nevertheless, we have a natural 
projection on diagonal matrices: 
\[ p_\ft \:  \fu(\cH) \to \ft, \qquad p_\ft(x)_j = \la e_j, x e_j \ra.\] 
The permutation group 
$S_J$ naturally acts on $\cH$ by $\sigma e_j := e_{\sigma(j)}$, 
normalizing $\ft$, which gives 
rise to a {\it big Weyl group} $\oline\cW \cong S_J \subeq \GL(\ft)$. 
Neumann's result asserts that 
\begin{itemize}
\setlength\itemsep{0em}
\item[\rm(a)] $ \oline{p_\ft(\cO_x)} = \oline\conv(S_J x)$ for $x\in\ft$. 
\item[\rm(b)] Every open or closed invariant subset $\Omega \subeq \fu(\cH)$ is 
determined by the $S_J$-invariant open subset $\Omega \cap \ft$ 
(see \cite[Thm.~4.33]{Neu00b} for closed subsets and use 
Lemma~\ref{lem:1.3} for open ones). 
\end{itemize}

If $x \in \fu(\cH)$ is not diagonalizable, then Neumann 
also describes $\oline{p_\ft(\cO_x)}$ in terms of the essential spectrum 
of $x$. The norm continuous unitary representations 
$(\pi_k)_{k \in \N}$ of $\U(\cH)$ on the spaces $\Lambda^k(\cH)$ 
provide a sequence of invariant convex functions 
by their support functionals $s_{\pi_k}$ (see \eqref{eq:spi} in the introduction). 
To evaluate these functionals on diagonal operators, we observe that 
the restrictions of $\pi_k$ to the diagonal subgroup 
$T$ is diagonalizable with the weights 
\[ x_F := \sum_{j \in F} x_j, \quad \mbox{ where } \quad 
F \subeq J, |F| = k, \qquad 
x_j := \la e_j, x e_j \ra.\] 
For $x \in \ell^\infty(J,\R)$ and the corresponding element 
$\diag(x) \in i\ft$, we then have 
\begin{equation}
  \label{eq:sk-funct}
s_k(x) := s_{\pi_k}(-i\diag(x))  = \sup \{ x_F \: 
F \subeq J, |F| = k \}.
\end{equation}
Neumann then shows that, if $J$ is countable, then 
the support functionals $s_{\pi_k}$ 
suffice to determine closed convex hulls of adjoint orbits 
(\cite[Thm.~4.37]{Neu00b}, \cite[Prop.~2.8]{Neu99}): 
\begin{align} 
\oline\conv(\cO_x) 
&= \{ y \in \fu(\cH) \:  (\forall k \in \N)
\ s_{\pi_k}(y) \leq s_{\pi_k}(x), \ s_{\pi_k}(-y) \leq s_{\pi_k}(-x) \} 
\label{neu:facts1}\\
\oline\conv(S_J x) &= \{ y \in \ell^\infty(J,\R)  \:  
(\forall k \in \N)\ s_k(y) \leq s_k(x), s_k(-y) \leq s_k(-x)\}.
\label{neu:facts2} 
\end{align}

\begin{Example} \label{ex:schur}
If $\cH = \C^n$, then we likewise obtain finitely many invariant 
convex functions $s_k, k =1,\ldots, n,$ on $\R^n$, such that 
\[  \conv(S_n x) 
= \{ y \in \R^n\: s_k(y) \leq s_k(x), k =1,\ldots, n-1; 
s_n(y) = s_n(x)\}.\] 
If $x$ is decreasing, then $s_k(x) = x_1 + \cdots +x_k$, 
so that one obtains the classical 
Hardy--Littlewoord--P\'olya inequalities describing convex hulls 
of $S_n$-orbits in $\R^n$ (cf.~\eqref{eq:schurhorn}, \cite{Ry65}). 
This idea is fundamental in Neumann's generalization of the 
Schur--Horn--Kostant Theorem to the pair 
$(\g,\ft) = (\fu(\cH), i\ell^\infty(J,\R))$. 
\end{Example}

\begin{Remark} \label{rem:5.4} 
Each $s_k$ extends by the same formula 
to a lower semicontinuous function  
${s_k \: \R^J \to (-\infty, \infty]}$ 
on the full space of real-valued functions on~$J$. 
These functions are 
invariant under the full permutation group~$S_J$. 
If $F \subeq J$ is a $k$-element subset, then 
\[ s_k(\lambda) 
= \sup \la \lambda, S_J e_F \ra 
= \sup \la \lambda, S_{(J)} e_F \ra 
\quad \mbox{ for } \quad 
e_F := \sum_{j \in F} e_F \in \R^{(J)} \cong (\R^J)'.\] 
It is easy to see that $s_k(\lambda) < \infty$ if and only if $\lambda$ is bounded from 
above. Accordingly, $s_k(\pm \lambda) < \infty$ is equivalent to the boundedness 
of $\lambda$, and by Proposition~\ref{prop:3.2} this is equivalent to the 
(semi-)equicontinuity of $S_{(J)}\lambda$, resp., $S_J\lambda$. 
\end{Remark}


Neumann also generalized Kostant's Convexity Theorem  
to the full orthogonal and symplectic Lie algebras 
of bounded operators. We refer to \cite{Neu02} for details. 

From Proposition~\ref{prop:3.2}(a) one derives in particular 
that every open invariant cone in $\fu(\cH)$ intersects 
the center $i\one$ (\cite[Thm.~5.6]{Ne12}, see also 
Remark~\ref{rem:center-compact}) and that, 
for a real Hilbert space $\cH$ of dimension $>2$ 
and a quaternionic  Hilbert space $\cH$, all open invariant 
cones in $\fo(\cH)$, resp., $\fu_\H(\cH)$ are trivial. 
These are interesting analogs of the corresponding observations 
concerning the finite-dimensional compact Lie algebras~$\fu_n(\K)$. 

{\bf Unitary Lie algebras of compact operators.} 
One may also consider the the unitary Banach--Lie groups  
$\U_p(\cH)$, $1 \leq p \leq \infty$, whose Lie algebras 
are the Banach spaces $\fu_p(\cH)$ of skew-hermitian 
operators of Schatten class~$p$ with 
$\|x\|_p = \tr(|x|^p)^{1/p}$. For these Lie algebras we have: 
\begin{Proposition} \label{prop:schatten} 
If $\cH$ is an infinite-dimensional complex Hilbert space, 
then $\fz(\fu_p(\cH)) =\{0\}$, and $\fu_p(\cH)$ contains non-trivial open 
invariant cones if and only if $p = 1$. 
\end{Proposition}
\begin{proof}  That $\fz(\fu_p(\cH)) = \{0\}$ follows from the 
fact that any operator commuting with all rank-one operators is a 
multiple of $\one$, which is not a compact operator. 

For $p = 1$, the trace functional is continuous, hence defines an 
invariant open half space. To obtain a pointed cone, consider the subset 
\[ \Delta := \{ x \in \Herm_1(\cH) \: x \geq 0, \tr x = 1\} \] 
which is bounded, $\|\cdot\|_1$-closed, convex and $\U(\cH)$-invariant,
and its distance from $0$ is~$1$. Hence Lemma~\ref{lem:a.6}(b) applies to 
any open convex set of the form $\Omega := \Delta + B_r(0)$, $0 < r < 1$, 
and shows that $\R^\times_+ \Omega$ is a pointed open invariant cone.  

For $1 < p < \infty$, let $\Omega \subeq \fu_p(\cH)$ be an open invariant cone.
Intersecting with diagonal operators with respect to an orthonormal 
basis, we obtain an open cone $W \subeq \ell^p(J,\R)$, 
invariant under the action of the group $S_{(J)}$ of finite permutations. 
To see that $0 \in W$, it suffices to assume that $J =\N$. 
As $W$ is open, it contains an element $x$ with finite support in~$\N$. 
For every $y \in \ell^p(\N)$ with finite support disjoint from 
$\supp(x)$, there exists an $\eps > 0$ with $x \pm \eps y \in W$. 
For $k \in \N$, pick permutations $\sigma_1, \ldots, \sigma_k \in S_{(\N)}$ 
fixing all elements in $\supp(y)$, and 
for which the subsets $\sigma_j(\supp x)$ are mutually disjoint. Then 
\[\pm \eps y + \frac{1}{k} \sum_{j = 1}^k \sigma_j x  \in W 
\quad \mbox{ and } \quad 
 \frac{1}{k}\Big\|\sum_{j = 1}^k \sigma_j x\Big\|_p  
= \frac{1}{k} (k \|x\|_p^p)^{1/p} 
= k^{-1 + \frac{1}{p}} \|x\|_p,\] 
and this expression converges to $0$ for $k \to \infty$. 
This implies that $\pm y \in \oline W$. By permutation invariance, 
it follows that $H(\oline W)$ contains all elements of finite support, 
and hence that $\oline W = \ell^p(\N)$. This in turn implies 
$W = \ell^p(\N)$ (Lemma~\ref{lem:bou}). 

For $p = \infty$ we argue similarly with 
$\frac{1}{k}\big\|\sum_{j = 1}^k \sigma_j x\big\|_\infty \leq \frac{1}{k}$. 
\end{proof}

\begin{Remark} \label{rem:extend-action}
(a) 
The Banach--Lie algebras $\fu_p(\cH)$ exhibit an 
interesting infinite-dimensional feature, namely that the unitary group 
$\U(\cH)_s$, i.e., the group $\U(\cH)$, endowed with the strong operator topology, 
 acts on these Lie algebras by conjugation 
as automorphisms: 
\[ \Ad \: \U(\cH)_s \times \fu_p(\cH) \to \fu_p(\cH), \quad 
(g,X) \mapsto gXg^{-1},\] 
and this action by isometries is continuous. 
As the subgroup $\U_p(\cH)$ is dense in $\U(\cH)_s$, 
the subgroup $\Ad(\U_p(\cH)) \subeq \cL(\fu_p(\cH))$ is dense 
in $\Ad(\U(\cH))$. By Lemma~\ref{lem:endo-strong-closed}, 
it follows that every weak-$*$-closed convex subset 
$C \subeq \fu_p(\cH)'$ is $\Ad(\U(\cH))'$-invariant and 
thus, by the Duality Theorem~\ref{thm:2.3}, also every 
open invariant cone, and further every 
continuous positively homogeneous convex function on such a cone. 

\nin (b) One can even go one step further by observing that 
the closure of $\U(\cH)$ in the space $\cL(\cH)$ with respect to 
the strong operator topology is the monoid 
\[ \Isom(\cH) = \{ A \in \cL(\cH) \: A^* A = \1 \} \] 
of isometries.\begin{footnote}{Since each $A \in \Isom(\cH)$ is 
determined by the orthonormal family $f_j := Ae_j$ and we may w.l.o.g.\ 
assume that $f_j = e_{\sigma(j)}$ for an injection 
$\sigma \: J \to J$, the density of $\U(\cH)$ in $\Isom(\cH)$ 
follows from the density of the monoid $\Inj(J)$ of injections 
in the group $S_J = \Bij(J)$ of all bijections with respect to the 
topology of pointwise convergence.   
}\end{footnote}
This monoid also acts continuously on $\fu_p(\cH)$ 
by Lie algebra endomorphisms via 
\[ \Ad(g)X := g X g^*,\] 
and this action by contractions is continuous with respect to the 
strong operator topology on $\Isom(\cH)$. Here we use that 
\[ gXYg^* = gXg^*gYg \quad \mbox{ for } \quad g \in \Isom(\cH).\] 
The same argument as under (a) now shows that every $\Ad(\U_p(\cH))$-invariant 
weak-$*$-closed convex subset of $\fu_p(\cH)$ is invariant under 
$\Isom(\cH)$. 

If $\cH = L^2(X,\fS,\mu)$, then each $\fS$-measurable map 
$\phi \: X \to X$ with $\phi_*\mu = \mu$ defines an isometry on 
$\cH$. Hence this monoid acts in particular by the adjoint 
action on $\fu_p(\cH)$ (cf.\ Subsection~\ref{subsec:5.3} 
and Remark~\ref{rem:5.8} below). If $X$ is finite and non-atomic, 
\cite[Thm.~5]{Br66} implies that the strong closure of the 
group of invertible measure preserving 
transformations is  the monoid of measure preserving transformations. 

\nin (c) On the level of diagonal matrices, we have 
the action of the Weyl  group $\cW \cong S_{(J)}$ 
on $\ell^p(J,\R)$, and this action extends to a continuous 
action of the monoid $\Inj(J)$ that has the same invariant 
weak-$*$-closed convex subsets in $\ell^p(J,\R)' \cong \ell^q(J,\R)$, 
where $p$ and $q$ are related by $p^{-1} + q^{-1} = 1$. 
In particular, $S_{(J)}$ and the much larger monoid $\Inj(J)$ 
have the same open invariant cones and the same 
invariant lower semicontinuous positively homogeneous convex functions 
on $\fu_p(\cH)$. 
Note that $\Inj(J)$ corresponds to the measure preserving maps 
on $J$ for $\fS = 2^J$ and the counting measure $\mu = \sum_{j \in J} \delta_j$.
\end{Remark}

\begin{Problem} \label{prob:2} Suppose that $\mu$ is $\sigma$-finite, so that 
$L^\infty(X,\fS,\mu)$ is the dual space of 
$L^1(X,\fS,\mu)$. Is it possible to describe the weak-$*$-closures 
of orbits of the group $\Gamma = \Meas(X,\fS,\mu)$ of (measurably invertible) 
measure preserving transformations 
in  $L^\infty(X,\fS,\mu)$? Note that, for every $E \in \fS$ of finite measure, 
we obtain a lower semicontinuous invariant convex function by \\[-4mm]
\[ s_E(f) := \sup_{\gamma \in \Gamma} \int_{\gamma(E)} f\, d\mu.\] 
For the special case where $X$ is countable, $\fS = 2^X$ 
and $\mu$ is the counting measure, 
this specializes to the situation of A.~Neumann's Theorem 
as in \eqref{neu:facts2}. 
\end{Problem}

\begin{Remark} \label{rem:5.8}
In \cite{Ry65, Ry67} J.V.~Ryff studies the 
case $X = (0,1)$, endowed with Lebesgue measure~$\mu$. We 
write $\cM$ for the monoid of measure preserving transformations of~$(0,1)$.
For every function $f \in L^1(0,1)$, the probability measure $f_*\mu$ 
defines a right-continuous decreasing function 
\[ f^* \: (0,1) \to \R, \quad f^*(s) := \sup\{ y \in \R \: \mu(\{ f > y\}) > s \} \] 
with the property that $f^* \in L^1(0,1)$ with 
\begin{equation}
  \label{eq:intineq}
\int_0^s f \leq \int_0^s f^* \quad \mbox{ for } \quad 0 < s < 1 \quad \mbox{ and }\quad
\int_0^1 f =  \int_0^1 f^*.
\end{equation} 
For $f,g \in L^1(0,1)$, we define $g \prec f$ ($f$ {\it majorizes} $g$) if 
\begin{equation}
  \label{eq:intineq2}
\int_0^s g^* \leq \int_0^s f^* \quad \mbox{ for } \quad 0 < s < 1 \quad \mbox{ and }\quad
\int_0^1 f =  \int_0^1 f^*.
\end{equation} 
For each $f \in L^1(0,1)$, the set $\Omega(f) 
:= \{ g \in L^1(0,1)  \: g \prec f\}$ of all functions 
majorized by $f$ then has the following properties: 
\begin{itemize}
\item It is a weakly compact convex subset of $L^1(0,1)$ 
(\cite[Thms.~2,3]{Ry65}). 
\item $\Ext(f) = \cM.f^* = \{ f \circ \sigma \: \sigma \in  \cM\}$ 
is the set of all functions equimeasurable with $f$ (i.e., $f^* = g^*$),  
in particular $\Omega(f) = \Omega(f^*)$ (\cite[Thm.~1]{Ry65}, \cite[Thm.~1]{Ry67}). 
\end{itemize} 
We thus obtain 
\[ \Omega(f) = \oline\conv(\cM.f), \] 
and that this set is characterized by the generalizations of the 
Hardy--Littlewood--P\'olya inequalities given by  \eqref{eq:intineq2}.
\end{Remark}

{\bf Operator algebras.}  
In the context of operator algebras, results resembling those 
for the full unitary group have been obtained 
recently for type~II factors and projections of unitary orbits onto maximal 
abelian subalgebras \cite{AM07, AM13}. See also \cite{HN91} for a general 
analysis of closed convex hulls of unitary orbits in von Neumann algebras. 
In this context one encounters maximal abelian 
subalgebras 
\[ \ft \subeq \g := \fu(\cA) := \{ x \in \cA \: x^* = - x\} \]
of the type $L^\infty(X,\fS,\mu)$ and ``Weyl groups'' acting by measure 
preserving transformations on this space 
(see Problem~\ref{prob:2}). 

In a $C^*$-algebra $\cA$ one often considers abelian subalgebras 
$\cT \subeq \cA$ and a 
continuous contraction $E \: \cA \to \cT$ (a {\it conditional expectation}), i.e., 
map preserving positivity and satisfying the equivariance condition 
\[ E(D_1 A D_2) = D_1 E(A) D_2 \quad \mbox{ for } \quad A \in \cA, D_1, D_2 \in \cT.\] 
Let $\cO_X = \{ UXU^* \: U \in \U(\cA)\}$ denote the adjoint orbit 
of $X = - X^* \in \fu(\cA)$. 
If $E \: \cA \to \cT$ is a conditional expectation, then we expect a 
relation of the type 
\[  E(\oline\conv(\cO_X)) = \oline\conv(\cW X) \quad \mbox{ for } \quad 
X \in \cT,\] 
where $\cW \subeq \GL(\cT)$ is an analog of the Weyl group. 

\begin{Example}
A particularly important example is the surjective conditional expectation 
$E \: \cM \to Z(\cM)$ of a $W^*$-algebra $\cM$ with values in its center 
(cf.\ \cite{Sa71}). 
In this case the action of $\cW$ is trivial and 
\[ \{E(X)\} = E(\oline\conv(\cO_X)) = \oline\conv(\cO_X) \cap Z(\cM), \] 
so that the convexity theorem holds in a trivial way. 
The situation becomes more 
interesting if we project onto a maximal abelian subalgebra which is 
not central; see 
\cite{AM07, AM13} for type II algebras. 
In this context one obtains interesting generalizations of the convex functionals 
$s_{\pi_k}$ from \eqref{eq:sk-funct}. Using a finite trace $\tau$, one obtains 
the $\U(\cA)$-invariant continuous convex/concave functionals 
(\cite[\S 4]{AM13}): 
\begin{align*}
U_t(a) &:= \sup \{ \tau(ap) \: p^2 = p = p^*, \tau(p) = t \}, \\  
L_t(a) &:= \inf\{ \tau(ap) \: p^2 = p = p^*, \tau(p) = t \}, \qquad t \geq 0. 
\end{align*}
\end{Example}

\begin{Remark} If $\cM$ is a von Neumann algebra 
of operators on a separable Hilbert space 
and $\cT \subeq \cM$  maximal abelian, then $\cT \cong L^\infty(X,\fS,\mu)$ for a finite 
measure~$\mu$, with predual $\cT_* = L^1(X,\fS,\mu)$. 
The majorization order on real functions in $\cT_*$ is expected to have 
natural connections with the theory of decreasing rearrangements 
(cf.~Remark~\ref{rem:5.8}). 
\end{Remark}

\pagebreak 
\subsection{Projective limits} 

Let $\g = \prolim \g_j$ be a projective limit of 
locally convex Lie algebras  $\g_j$ and $W \subeq \g$ an open invariant cone 
(\cite{HM07}). 
If $x \in W$, then $W - x$ is a $0$-neighborhood in $\g$, hence 
contains the kernel of a projection $p_j \: \g \to \g_j$. 
This implies that $x + \ker p_j \subeq W$, so that 
$\ker p_j\subeq H(W)$ (Lemma~\ref{lem:limcone}(iii)). 
Therefore $W$ is the inverse image of an 
open invariant cone in one of the Lie algebras $\g_j$. 
In particular, proper projective limits never contain pointed 
open invariant cones. 

\section{Direct limits}  
\label{sec:6} 

In this section we describe some special features of 
direct limit Lie algebras, first in the context of 
root decompositions and then for more general situations.

\subsection{Direct limits of finite-dimensional Lie algebras} 

First we consider the case, where $\ft \subeq \g$ is an elliptic 
Cartan subalgebra and $\g = \g^{\rm alg}$ is locally 
finite, i.e., every finite subset of $\g$ generates a finite-dimensional 
Lie subalgebra. Then $\g$ is the union of subalgebras 
$\g_F$ with 
\[ \g_{F,\C} = \ft_\C + \sum_{\alpha \in \Delta_F} \g_{F,\C}^\alpha, \qquad 
\Delta_F \subeq \Delta \ \mbox{ finite }, \qquad 
\dim \g_{F,\C}^\alpha < \infty.\] 
As $\Delta_F^\bot \subeq \ft_\C$ is of finite codimension, these subalgebras 
can be analyzed with the methods from \cite{Ne99}, outlined in Section~\ref{sec:4}.

\begin{Examples} \label{ex:3.15} 
(a) A particularly interesting case arises for $\Delta = \Delta_k$, 
i.e., when all roots are compact. 
One can show that the simple locally finite Lie algebras of this 
type are $\su(J,\C)$ and $\fu(J,\K)$, $\K = \R, \H$, for an infinite 
set $J$ (see \cite{St01}, \cite{Ne12} for details). 

Let $\lambda \in \ft'$. Applying Kostant's Convexity Theorem 
(Theorem~\ref{thm:compact-lie}(ii)) to the subalgebras $\g_F$, it follows that 
\[ p_{\ft'}(\cO_\lambda) = \conv(\cW\lambda).\] 
If $\cO_\lambda$ is semi-equicontinuous, the relation 
$\fg = \Ad(G)\ft$ implies that $B(\cO_\lambda)^0$ intersects~$\ft$, 
so that $\cW\lambda \subeq \ft'$ is also semi-equicontinuous 
(Remark~\ref{rem:maxcone}(b)). 

If $\g$ is one of the three  infinite-dimensional Lie algebras 
$\fu(J,\K)$, $\K = \R,\C,\H$, 
Lemma~\ref{lem:uj-seq} below implies that 
$\cW\lambda$ is semi-equicontinuous if and only if 
$\lambda(\Delta^\vee)$ is bounded.

\nin (b) If not all roots are compact, then $\g = \fk \oplus \fp$ 
with 
\[ \fk_\C = \ft_\C + \sum_{\alpha \in \Delta_k} \g_\C^\alpha
\quad \mbox{ and } \quad \fp_\C = \sum_{\alpha \in \Delta_p} \g_\C^\alpha.\] 
The semi-equicontinuity 
of $\cO_\lambda$ for $\lambda \in \ft'$ implies in particular that 
$\cO_\lambda^K := \Ad'(K)\lambda$ 
is semi-equicontinuous. 
If $\fk$ is a direct sum of a center and finitely many simple 
summands, then Lemma~\ref{lem:uj-seq} easily implies that 
$\cO_\lambda^K$ is semi-equicontinuous if and only if 
$\lambda(\Delta_k^\vee)$ is bounded. Here a key point is that 
the $\cW$-action on $\Delta_k$ has only finitely many orbits. 

A typical examples is the real form  
$\g = \fu(J_1, J_2,\C) \subeq \g_\C = \gl(J,\C)$, 
where $J = J_1 \dot\cup J_2$. 
Then $\ft = i \R^{(J)}$ and 
the Weyl group is 
$\cW \cong S_{(J_1)} \times S_{(J_2)} \subeq S_{(J)}$, acting 
by permutations. Up to sign, the only $\cW$-invariant positive system of 
non-compact roots is 
\[ \Delta_p^+ = \{ \eps_j - \eps_k \: j \in J_1, k \in J_2\},
\quad \mbox{ so that } \quad C_{\rm max}(\Delta_p^+) 
= \{ i x \:   x\res_{J_2} \leq x\res_{J_1} \}.\] 
If $J_1$ and $J_2$ are both infinite, we obtain the cone 
\[ C_{\rm max}(\Delta_p^+) 
= \{ - i x \:   x\res_{J_2} \leq 0 \leq x\res_{J_1} \}, \] 
which has no interior points. 
\end{Examples}

\begin{Lemma} \label{lem:uj-seq}
For $\g = \fu(J,\K)$, $\K = \R, \C, \H$ 
and an elliptic Cartan subalgebra $\ft$ obtained from a root decomposition of 
$\g_\C$, the following are equivalent for $\lambda \in \ft'$: 
\begin{itemize}
\setlength\itemsep{0em}
\item[\rm(i)] $\cO_\lambda \subeq \g'$ is semi-equicontinuous. 
\item[\rm(ii)] $\cO_\lambda \subeq \g'$ is equicontinuous. 
\item[\rm(iii)] $\cW\lambda \subeq \ft'$ is semi-equicontinuous. 
\item[\rm(iv)] $\cW\lambda \subeq \ft'$ is equicontinuous. 
\item[\rm(v)] $\lambda(\Delta^\vee)$ is bounded. 
\end{itemize}
\end{Lemma}

\begin{proof} The equivalences of (i) and (iii) and of 
(ii) and (iv) follow easily from Kostant's Convexity Theorem because
 $\g = \Ad(G)\ft$ (cf.~Examples~\ref{ex:3.15}(a) 
and Theorem~\ref{thm:compact-lie}). 

For a suitably chosen set $J'$, we can identify $\ft$ with $i \R^{(J')}$, 
and then $\cW$ contains the subgroup $S_{(J')}$ of finite permutations on $J'$ 
and the elements $e_j - e_k, j\not= k \in J'$, are coroots. 

If $\cW\lambda$ is semi-equicontinuous, then 
Proposition~\ref{prop:3.2}(c) implies that $\lambda \in \ft' \cong i \R^{J'}$ 
is bounded. This implies (v) and from the concrete description of the 
Weyl group in all cases, we immediately see that $\cW\lambda$ is 
equicontinuous. That (v) implies (iv) follows easily from the fact that 
by evaluating $\lambda$ on the coroots $e_{j} - e_k$, it follows that 
$\lambda$ is bounded. 
\end{proof}

\subsection{Direct limits of compact Lie algebras} 

In this subsection we assume that $\Delta =\Delta_k$. 
We have already seen in Example~\ref{ex:3.15}(a) 
that, for $\g = \fu(J,\K)$, any $\lambda \in \ft' \cong \R^J$ corresponding 
to a  semi-equicontinuous coadjoint orbit is bounded. 
This is reflected by the following proposition, which represents the 
information we have on $\g$ itself. 
The convex subset 
\[ B_\ft := \conv(i\Delta^\vee \cup -i\Delta^\vee) \subeq \ft \] 
is symmetric and invariant under the Weyl 
group~$\cW$.  
Since it is also generating, 
Kostant's Convexity 
Theorem (Theorem~\ref{thm:compact-lie}(ii)), 
applied to finite-dimensional subalgebras, 
implies that $B:= \Ad(K)B_\ft$ is a symmetric convex invariant
generating subset with $B \cap \ft = B_\ft$. Therefore 
\[  \|x\|_{\rm max} := \inf \{ t > 0 \: x \in t B \} \]
is an invariant norm on $\fg$ with $\|i\alpha^\vee\|_{\rm max} \leq 1$ for each 
$\alpha \in \Delta$. 
The following lemma shows that $\|\cdot\|_{\rm max}$ is a maximal 
invariant norm. 

\begin{Proposition} \label{prop:6.1} 
If $\fg$ is $\fsl(J,\C)$ or $\fu(J,\K)$ for $\K = \R,\H$, 
and $\|\cdot\|$ is an invariant norm on $\fg$, then 
\[  c := \sup_{\alpha \in \Delta} \|i\alpha^\vee\| < \infty 
\quad \mbox{ and  } \quad 
\|x\| \leq c \|x\|_{\rm max} \quad \mbox{ for all } x \in \fg. \] 
\end{Proposition}

\begin{proof} Since $\fg$ is simple, there exists at most two 
${\cal W}$-orbits in $\Delta$ 
and $\alpha$ and $\beta$ lie in the
same orbit if and only if $\|\alpha\| = \|\beta\|$ 
(\cite[\S 10.4, Lemma~C]{Hu80}).  
This shows that $c$ is finite. 

As every element $x \in \fg$ is conjugate under
$G$ to an element of $\ft$ (Theorem~\ref{thm:compact-lie}), 
it suffices to consider elements $x \in \ft$. Then 
$\|y\|_{\rm max} \leq 1$ implies $\|y\| \leq c$. We conclude
that $\|\cdot\| \leq c \|\cdot\|_{\rm max}$. 
\end{proof}

\begin{Example} \label{ex:6.2} To determine 
the maximal norm $\|\cdot\|_{\rm max}$ more explicitly, we 
consider a set $J$ and the simple Lie algebra 
$\g_\C = \fsl(J, \C)$ endowed with its canonical involution 
$x^* = \oline x^\top$, so that $\fg = \su(J,\C) 
:= \{ x \in \fsl(J, \C) \: x^* = - x\}$. 

Let $E_{ij} = (\delta_{ik} \delta_{j\ell})_{k,\ell \in J} 
\in \C^{J \times J}$ be the matrix with one entry $1$ in the 
position $(i,j)$. 
For $\alpha = \eps_i - \eps_j$ we have $\alpha^\vee = E_{ii} -
E_{jj}$ and $\|i\alpha^\vee\|_1 = 2$, where $\|\cdot\|_1$ is the
trace norm. We claim that $\|x\|_1 = 2\|x\|_{\rm max}$ on all of $\fg$. 
Since all coroots are conjugate under the Weyl group 
$\cW \cong S_{(J)}$, Proposition~\ref{prop:6.1} shows that 
$\|x\|_1 \leq 2 \|x\|_{\rm max}$ for all $x \in \fg$. 

If $J$ is an infinite set, then one readily verifies that the extreme points of the set 
\begin{equation}
  \label{eq:S}
 S := \Big\{ x \in \R^{(J)} \: \sum_j x_j = 0, \sum_j |x_j| \leq 2 \Big\}
\end{equation}
consist of the points of the form $e_i - e_j$, $i \not= j$. In fact, 
if $x = (x_j)_{j \in J} \in S$ has two positive entries 
$x_i$ and $x_j$, then $0 < x_i, x_j < 2$ implies the existence of an 
$\eps > 0$ such that 
$\|x\|_1 = \| x + t(e_i - e_j)\|_1 = 2$ for $|t| < \eps$, so that it cannot be
an extreme point. Therefore at most one entry $x_i$ is positive, and
likewise one entry $x_j$ is negative. From $x_j = - x_i$ we conclude
that $x_i = 1 = - x_j$, i.e., $x = e_i - e_j$. Therefore 
$$ S= \conv \{ e_i - e_j \: i \not=j\}. $$
Letting $J$ be arbitrary, this argument still shows that 
$2 \|x\|_{\rm max} = \|x\|_1$ for $x \in \ft$, hence for $x \in \fg$. 
\end{Example}

\begin{Corollary} \label{cor:6.5} 
If $\fg$ is $\su(J,\C)$ or $\fu(J,\K)$ for $\K \in \{\R,\C,\H\}$ and 
an infinite set~$J$, then the following assertions hold: \\[-7mm]
\begin{itemize}
\setlength\itemsep{0em}
\item[\rm(i)] For every invariant norm $\|\cdot\|$ on $\g$, there exists a 
$C > 0$ with $\|x\| \leq C\|x\|_1$, where $\|x\|_1 = \tr(|x|)$ 
denotes the trace norm.  
\item[\rm(ii)] For $\cH := \ell^2(J,\K)$, every norm-continuous representation 
$\rho \: \g \to \cL(E)$ on a Banach space $E$  extends 
to a representation of the Banach--Lie algebra 
completion $\su_1(\cH)$ for $\g = \su(J,\C)$, and to 
$\fu_1(\cH)$ for $\g = \fu(J,\K)$. 
\item[\rm(iii)] 
Every $\Ad(G)$-invariant convex subset of $\fg$ 
which is open with respect to the finest locally convex topology 
is also open with respect to $\|\cdot\|_1$. 
\end{itemize}
\end{Corollary}

\begin{proof} (i) follows from Proposition~\ref{prop:6.1} for the case 
where $\g$ is simple. For 
$\g = \fu(J,\C)$ we write $\g \cong \su(J,\C) \rtimes \R$ 
and observe that the open $\|\cdot\|$-unit ball $B$ intersects 
$\su(J,\C)$ in a $\|\cdot\|_1$-open subset. Since $\su(J,\C)$ has 
finite codimension, $B$~contains a $\|\cdot\|_1$-ball, and this implies~(i).

\nin (ii) follows from (i). 

\nin (iii) First we assume that $\g \not= \fu(J,\C)$, so that $\g$ is simple. 
Let $\Omega \subeq \g$ be an invariant convex subset 
which is open for the finest locally convex topology on $\g$, i.e., 
$\Omega - x_0$ is absorbing for each $x_0\in\Omega$.
Since $\g$ is a union of simple compact Lie algebras $\g_F$, 
$F \subeq J$ finite, we have $\g_F \cap \Omega \not=\eset$ for such a 
subset~$F$. Then $\eset \not= \Omega \cap \fz(\g_F) \subeq \{0\}$ 
(Remark~\ref{rem:center-compact}) shows that $0 \in \Omega$. 
Hence $\Omega \cap - \Omega$ is the open unit ball for an invariant norm, 
hence contains an interior point with respect to $\|\cdot\|_1$ 
and is therefore $\|\cdot\|_1$-open by Lemma~\ref{lem:new}. 

Now we consider the case $\g = \fu(J,\C)$. Using Konstant's Convexity 
Theorem, we have to show that every $S_{(J)}$-invariant 
open convex subset of $\R^{(J)}$ in $\ell^1$-open, 
resp., contains an $\ell^1$-inner point (Lemma~\ref{lem:new}). 
As in the proof of Proposition~\ref{prop:3.2}(b), we find a 
finite subset $F \subeq J$ and $a \not= 0$ with $x_0 := a e_F \in \Omega$. 
Then $\Omega - x_0$ is an open $0$-neighborhood invariant 
under the group $S_{(J \setminus F)}$. This implies that 
$(\Omega - x_0) \cap \R^{(J \setminus F)}$ is $\|\cdot\|_1$-open 
because it contains the $S_{(J\setminus F)}$-orbit of some 
element of the form $e_j - e_k$, $j \not=k \in J \setminus F$ 
(Lemma~\ref{lem:new}). 
As $\R^{(J\setminus F)}$ has finite codimension in $\R^{(J)}$, it 
follows that $\Omega$ is also open with respect to $\|\cdot\|_1$. 
%
\end{proof}

The maximal invariant norm $\|\cdot\|_{\rm max}$ 
on $\fg$ leads to a to a minimal Banach 
completion of the Lie algebra $\fg$. 
In this sense, Example~\ref{ex:6.2} shows that, 
for any infinite-dimensional Hilbert space over $\K = \R,\C,\H$, 
the subspace 
$\fu_1(\cH)$ of the Schatten ideal $B_1(\cH)$ 
on the Hilbert space $\cH = \ell^2(J,\K)$, is the minimal Banach
completion of the Lie algebra $\fu(J,\K)$ with respect to 
any $\U(J,\K)$-invariant norm.

\subsection{Euclidean Lie algebras} 
\label{subsec:6.3}

An extremely important tool for studying invariant convex subsets of 
a Lie algebra $\g$ are invariant symmetric bilinear forms. 
There are many interesting infinite-dimensional Lie algebras 
$\g$ which carry an  invariant positive definite form $\kappa$. 
We then call the pair $(\g,\kappa)$ a {\it euclidean Lie algebra}. 

\begin{Examples} (a) Compact Lie algebras (cf.\ Subsection~\ref{subsec:4.1}). 

\nin (b) Hilbert--Lie algebras (Definition~\ref{def:hilbertliealg}), such as 
the  Lie algebra $\fu_2(\cH)$ of skew-hermitian Hilbert--Schmidt operators on 
a Hilbert space $\cH$ over $\K = \R,\C,\H$, with $\kappa(x,y) = -\tr_\R(xy)$. 

\nin (c) Mapping Lie algebras, such as 
$C^\infty(M,\fk)$ for a compact Lie algebra $\fk$ and a compact 
smooth manifold $M$. Here every positive measure $\mu$ 
on $M$ leads to an invariant form on $C^\infty(M,\fk)$ by 
\[ \kappa(\xi,\eta) := \int_M \kappa_\fk(\xi(m),\eta(m))\, d\mu(m).\] 
Every vector field $X\in \cV(M)$ whose flow preserves $\mu$ 
defines a $\kappa$-skew derivation. 

Here the case $M = \bS^1 \cong \R/2\pi\Z$ is  of particular importance. 
In this case 
\[ \kappa(\xi,\eta) := \frac{1}{2\pi}\int_0^{2\pi}
 \kappa_\fk(\xi(t),\eta(t))\, dt\] 
and a vector field on $\bS^1$ preserves the measure if and only if it 
is constant, i.e., generates rigid rotations. 

\nin (d) Poisson bracket algebras, such as 
$(C^\infty(M),\{\cdot,\cdot\})$ 
for a compact symplectic manifold $(M,\omega)$, where we put 
$\kappa(F,H) := \int_M FH \omega^n.$ 
Here every symplectic vector field defines a 
$\kappa$-skew derivation. 
\end{Examples}

The following proposition applies to the increasing sequence of the 
Lie algebras $\su_n(\K)$ for $\K = \R,\C,\H$, but also to more 
complicated sequences, such as $\su_{2^n}(\K)$. In particular, 
we do not require the existence of an elliptic Cartan subalgebra. 

\begin{Proposition} \label{prop:6.5} 
Let $\g$ be a locally finite Lie algebra which is the 
directed union of compact simple Lie algebras. 
Then $\g$ has an invariant scalar product  which 
is unique up to a positive constant. 
\end{Proposition} 

\begin{proof} The existence of an invariant symmetric bilinear 
form, which is invariant under all derivations, is shown 
in \cite[Prop.~II.1]{Ne05}. Let $0 \not= x \in \g$ and 
$\g_0 \subeq \g$ be a finite-dimensional compact simple 
Lie algebra containing $x$ and an element $y$ with 
$\kappa(x,y) \not=0$. Then $\kappa$ restricts to a non-zero 
invariant symmetric bilinear form on $\fg_0$, hence is either 
positive or negative definite. Replacing $\kappa$ by $-\kappa$ if 
necessary, we may assume that $\kappa$ is positive definite on 
$\fg_0$. Then it is also positive definite on all simple compact 
subalgebras $\fg_1$ containing $\fg_0$, and since $\g$ is exhausted 
by such subalgebras, $\kappa$ is positive definite. 

The uniqueness assertion follows immediately from the corresponding 
uniqueness for simple compact Lie algebras. 
\end{proof}

\begin{Examples} (a) For the Lie algebra $\g = \su(\N,\C)$, the increasing 
union of the simple compact Lie algebras $\g_n = \su_n(\C)$, 
a natural invariant scalar product is given by 
$\kappa(x,y) =-\tr(xy)$. The corresponding norm is the Hilbert--Schmidt norm, 
and by completion we obtain the simple Hilbert--Lie algebra 
$\fu_2(\ell^2(\N))$. In particular, the bracket extends to the completion. 

\nin (b) For the union $\g := \su_{2^\infty}(\C)$ of the Lie algebras 
$\g_n = \su_{2^n}(\C)$, defined by the embeddings 
$a \mapsto \pmat{ a & 0 \\ 0  & a}$, 
a compatible sequence of scalar products is given by 
$\kappa_n(x,y) = - 2^{-n} \tr(xy)$. 
In this case the Lie bracket is not continuous with respect to the norm 
defined by the corresponding scalar product on $\g$. 
\end{Examples}

\begin{Problem} (a) Does Proposition~\ref{prop:6.1} generalize 
in a suitable way to Lie algebras $\fg$ 
which are direct limits of simple compact ones but which do not 
contain an elliptic Cartan subalgebra? 

\nin (b) Is it true that every direct limit of compact Lie algebras 
has an invariant scalar product? 
By Proposition~\ref{prop:6.5} this is true for direct limits 
of simple  compact Lie algebras. 
\end{Problem}

\begin{Definition} \label{def:hilbertliealg} 
A {\it Hilbert--Lie algebra} is a euclidean 
Lie algebra which is complete with respect to the scalar product. 
\end{Definition}

\begin{Proposition} \label{prop:2.12} If $\g$ is a Hilbert--Lie algebra, then each 
non-empty open invariant convex subset $\Omega \subeq \g$ intersects the center. 
\end{Proposition}

\begin{proof} If $\Omega = \g$, 
there is nothing to show because $0 \in \Omega$. 
We may therefore assume that $\Omega \not=\g$. 
Then $\Omega$ is a proper open convex subset, and since the group 
$\Gamma := \la e^{\ad \g} \ra$ 
of inner automorphisms of $\g$ acts isometrically on 
$\g$ preserving $\Omega$, the closed convex subsets 
$$ \Omega_c := \{ x \in \Omega \: d_\Omega(x) \leq c \} $$
are $\Gamma$-invariant (Lemma~\ref{lem:1.3}). 
Each set $\Omega_c$ is a Bruhat--Tits space 
with respect to the induced Hilbert metric from $\g$ and 
since $\Gamma$ acts on this space isometrically with bounded 
orbits, the Bruhat--Tits Theorem~\ref{thm:BT} 
implies the existence of a fixed point $z \in \Omega_c\subeq \Omega$. 
Now it only remain to observe that 
$\z(\g)$ is the set of fixed points for the action of 
$\Gamma$ on $\g$. 
\end{proof}

\begin{Corollary} \label{cor:2.13} A Hilbert--Lie algebra contains a proper 
open invariant cone if and only if its center is non-trivial. 
\end{Corollary}

\begin{proof} That the existence of 
an open invariant cone implies that $\z(\g)\not=\{0\}$ 
follows from Proposition~\ref{prop:2.12}. 

If, conversely, there exist a non-zero $z_0 \in \fz(\g)$, and for 
$0 < r < \|z_0\}$, the open ball $B_r(z_0)$ 
is invariant with $0\not\in \oline{B_r(z_0)}$. 
Then Lemma~\ref{lem:a.6}(b) implies that 
$\R^\times_+ B_r(z_0)$ is a pointed open invariant cone. 
\end{proof}

\section{Double extensions} 
\label{sec:7} 

In this final section we take a closer look at 
double extensions (Definition~\ref{def:doubleext}). 
This is motivated in particular by their natural 
emergence from covariant positive energy representations 
(cf.\ \cite{JN16} and \cite[\S 8]{JN19}). 

\subsection{Lorentzian double extensions} 
\label{subsec:7.1} 

In addition to invariant scalar products, 
invariant {\it Lorentzian forms} on Lie algebras are also very 
useful to analyze convex hulls of adjoint and coadjoint orbits. 
They arise naturally from double  extensions of euclidean Lie algebras. 
If $(\fg,\kappa)$ is a euclidean Lie algebra 
(Subsection~\ref{subsec:6.3}), and 
$\hat\g = (\g,\kappa,D)$ is a double extension defined by the 
cocycle $\omega_D(x,y) = \kappa(Dx,y)$, where 
\[ D \in \der(\g,\kappa) := \{ D \in \der(\g) \: (\forall x,y \in \g)\ 
\kappa(Dx,y) + \kappa(x,Dy) = 0 \} \] 
is a $\kappa$-skew-symmetric derivation on $\g$, 
then 
\begin{equation}
  \label{eq:invbil}
\beta((z,x,t), (z',x',t')) := zt' + tz' - \kappa(x,x') 
\end{equation}
defines an invariant Lorentzian form on $\hat\g$ 
(\cite{MR85}). We call such double extensions {\it Lorentzian}. 

On $\hat\g$ we have an action $\Ad$ of the Lie group 
\[  G^\sharp := G \rtimes_{\alpha^G} \R, \] 
where $\alpha^G \: \R \to \Aut(G)$ is a smooth action for which the 
corresponding action on $\g$ is generated by the derivation~$D$. 

\begin{Definition} 
In the following we write 
$\bc := (1,0,0)$, $\bd := (0,0,1)$ and, accordingly,  
$\bc^* := (1,0,0)$ and $\bd^* := (0,0,1)$ for the corresponding 
elements in $\hat\g'$. 
Note that $\bc \in \hat\g$ and $\bd^*\in \hat\g'$ are $G$-invariant. 
For $x \in \g$ we put $\|x\| := \|x\|_\kappa:=\sqrt{\kappa(x,x)}$. 

The {\it $\beta$-dual of $\hat\g$} is the subset $\hat\g^\beta$ of 
$\hat\g'$, consisting 
of all linear functionals whose restriction to $\g$ is continuous with 
respect to the $\kappa$-topology, 
i.e., contained in the Hilbert completion of the euclidean space 
$(\g,\kappa)$. 
\end{Definition}

Since $\beta$ is a Lorentzian, the open subset 
\[  W := \{ \bx = (z,x,t) \in \g \: t > 0, \beta(\bx,\bx) > 0 \} 
\subeq (0,\infty) \times \g \times (0,\infty) \] 
is an open convex invariant cone and 
$\chi(\bx) := \beta(\bx,\bx)^{-1}$
defines on $W$ a smooth strictly convex invariant function 
which tends to infinity at the boundary 
(Proposition~\ref{prop:c.2} and Example~\ref{ex:2.8}). 

\begin{Proposition} \label{prop:1.2} The following assertions hold: \\[-7mm]
\begin{itemize}
\setlength\itemsep{0em}
\item[\rm(a)] $\Omega := W - \R_+ \bc =  \{ (z,x,t) \: t > 0 \}$ 
is an open invariant half space bounded by the hyperplane ideal 
$\tilde\g = \R \bc + \g$. 
\item[\rm(b)] For each $x \in \Omega$, we have 
$x \in \Ext(\oline\conv(\cO_x)).$
\item[\rm(c)] If $\lambda = (z^*, \alpha, t^*) \in \hat\g^\beta$ with 
$z^* > 0$ and $\bx = (z,x,t)$ with $t > 0$, then 
$\lambda(\cO_x) = \cO_\lambda(x)$ is bounded from below. 
\item[\rm(d)] $W^\star$ is contained in the $\beta$-dual $\hat\g^\beta$ 
of $\hat\g$. 
\item[\rm(e)] If $\lambda \in \hat\g^\beta$ satisfies $\lambda(\bc) = z^*\not=0$, 
then its orbit $\cO_\lambda$ is semi-equicontinuous. 
\end{itemize}
\end{Proposition}

\begin{proof} (a) For $\bx = (z,x,t)$ with $t > 0$ we  have 
$\beta(\bx,\bx) = 2zt - \kappa(x, x).$
For $z_0 := \frac{1}{2t} \kappa(x, x)$ and 
$\bx_0 := (z_0, x,t)$ we thus obtain 
\[  \bx = \bx_0 + (z - z_0) \bc\quad \mbox{ and } \quad 
\beta(\bx_0,\bx_0) =2z_0 t - \kappa(x,x) =0, \] 
so that $\bx_0 \in \partial W$, and for $z_1 > z_0$ we get $\bx_1 := (z_1, x,t) \in W$. 

\nin (b) Since $\bc$ is fixed by the adjoint action, 
$\cO_\bx = \cO_{\bx_0} + (z- z_0)\bc.$ 
The orbit $\cO_{\bx_0}$ lies in the set 
$\big\{ \big(\frac{\kappa(y,y)}{2t}, y,t\big) \: y \in \g \big\},$
which is the graph of the strictly convex function 
$f(y) := \frac{\kappa(y,y)}{2t} = \frac{\|y\|^2}{2t}$ on $\g$. Therefore 
\[  \oline\conv(\cO_{\bx_0}) \subeq 
\{ \big(z, y,t) \: y \in \g, z \geq f(y) \} \] 
and thus 
\[  \cO_{\bx_0} 
\subeq \{ \big(f(y), y,t\big) \: y \in \g \} 
= \Ext(\{ \big(z, y,t) \: y \in \g, z \geq f(y) \}) \]
implies $\cO_{\bx_0} \subeq \Ext(\oline\conv(\cO_{\bx_0})).$

\nin (c) In view of (a), we may w.l.o.g.\ assume that 
$\bx = \bx_0 \in \partial W$. Then 
\[  \lambda(f(y),y,t) 
= tt^* + \alpha(y) + f(y)z^* 
= tt^* + \alpha(y) + z^*\frac{\|y\|^2}{2t} \]
is bounded below on $\g$ because 
$|\alpha(y)| \leq \|\alpha\| \cdot \|y\|$ and 
$z^*/t > 0$. 

\nin (d) If $\lambda = (z^*,\alpha,t^*) \in W^\star$, then 
$\lambda$ is bounded on the set 
$\big((\bc + \bd) - W\big) \cap (W-(\bc+\bd)).$ 
This set consists of all elements 
$\bx_0 = (z,x,t)$, satisfying 
$(1\pm  z,\pm x,1\pm t)  \in W,$
which means that $|z| < 1$, $|t| < 1$ and 
$\|x\|^2 < 2(1-|z|)(1-|t|).$
Therefore $\alpha$ is bounded on an open ball in $\g$, 
hence continuous with respect to the $\kappa$-topology, i.e., 
$\lambda \in \hat\g^\beta$. 

\nin (e) Write $\lambda = (z^*, x^*, t^*)$ and 
observe that 
\[  \beta^*((z^*, x^*,t^*), (\tilde z^*, \tilde x^*, \tilde t^*)) 
:= z^* \tilde t^* + \tilde z^* t^* - \kappa(x^*, \tilde x^*) \] 
defines an invariant Lorentzian form on the $\beta$-dual. 
If $z^* = \lambda(\bc) \not=0$, then 
\[  \beta^*(\lambda + t \bd^*, \lambda + t \bd^*) 
= \beta^*(\lambda,\lambda) + 2 t \beta^*(\lambda, \bd^*) 
= \beta^*(\lambda,\lambda) + 2 t z^*, \]
so that $t := -\frac{1}{2z^*} \beta^*(\lambda,\lambda)$ leads for 
$\lambda_0 := \lambda + t \bd^*$ to $\beta^*(\lambda_0, \lambda_0) = 0$. Now 
$\cO_\lambda = \cO_{\lambda_0} - t \bd^*$
and since $\lambda_0 \in W^\star$ follows from (d), 
$\cO_\lambda$ is semi-equicontinuous because $W^\star$ is the dual of an open cone 
(Example~\ref{ex:2.2}(b)). 
\end{proof}

We now turn to global aspects of double extensions. 
Assume that $\g = \L(G)$ for a $1$-connected Lie group 
$G$ and that 
$\alpha^G  \: \R \to \Aut(G)$ is a smooth action 
for which the corresponding action on the 
Lie algebra $\alpha_t := \L(\alpha^G_t)$ 
satisfies $\alpha'(0) = D$. 
We may thus form the 
semidirect product group $G^\sharp := G \rtimes_{\alpha^G} \R$ for which 
there exists a smooth action on the central extension 
$\hat\g$ of $\g^\sharp := \g \rtimes_D \R$ by $\R$. 

The action of $\R$ is simply given by 
\[  \Ad(\exp s\bd)(z,x,t) = (z, \alpha_s(x),t). \] 
The $G$-action on $\g^\sharp = \g \rtimes \R \bd$ has the form 
\[  \Ad(g)(x,t) := (\Ad(g)x + t \gamma(g), t)  \] 
for the cocycle 
\begin{equation}
  \label{eq:drel}
 \gamma \: G \to \g, \quad \gamma(g) = \Ad(g)\bd - \bd.
\end{equation}
For the central extension $\tilde\g = \R \bc \oplus \g$, the pairing 
\[ \tilde\g \times \g^\sharp \to \R, \quad 
\big((z,x), (x',t')\big) \mapsto zt' \] 
is $\Ad(G)$-invariant, so that the $G$-action on $\tilde\g$ 
takes the form 
\[ \Ad(g)(z,x) 
 = \Big(z - \kappa(x,\gamma(g^{-1})), \Ad(g)x\Big)
 = \Big(z + \kappa(\Ad(g)x,\gamma(g)), \Ad(g)x\Big),\] 
where we use that $\gamma(g^{-1}) = - \Ad(g)^{-1}\gamma(g)$.  
Using the fact that the function $\beta(\bx,\bx)$ on $\hat\g$ is 
$\Ad(G)$-invariant, we obtain: 

\begin{Proposition}\label{prop:act-form} 
For $g \in G$, $z,t \in \R$ and $x \in \g$, we have 
\[  \Ad_{\g}(g)(z,x,t) = 
\Big(z - \kappa(\gamma(g^{-1}),x) + \frac{t}{2}\kappa(\gamma(g),\gamma(g)), 
\Ad(g)x + t \gamma(g), t\Big), \]
and in particular, 
\[  \Ad_{\g}(g)\bd = 
\Big( \frac{1}{2}\kappa(\gamma(g),\gamma(g)), \gamma(g), 1\Big) 
= \Big(\frac{1}{2}\|\gamma(g)\|^2, \gamma(g), 1\Big). \] 
The coadjoint action is given by 
\begin{align*}
& \Ad_{\g}^*(g^{-1})(z^*,x^*,t^*) 
=  (z^*,x^*,t^*) \circ \Ad_{\g}(g) \\
&= \Big(z^*, x^* \circ \Ad_\g(g) - z^* \kappa(\gamma(g)^{-1}, \cdot), 
 \frac{z^*}{2}\kappa(\gamma(g),\gamma(g)) 
+ x^*(\gamma(g)) + t^*\Big). 
\end{align*}
\end{Proposition}

\begin{Remark}
As $\bd$ is fixed by $\Ad(\exp \R \bd)$, 
the cocycle 
$\gamma \: G \to \g, \gamma(g) = \Ad(g)\bd - \bd$ 
is equivariant with respect to the action of $\exp(\R \bd)$ on 
$G$ and $\g$, respectively. 
It follows in particular that the smooth function 
$G \to \R,\quad  g \mapsto \kappa(\gamma(g), \gamma(g))$ 
is $\exp(\R \bd)$-invariant. 
We also see that $\Ad(G)\bd = \Ad(G^\sharp)\bd$ is an adjoint orbit in~$\hat\g$. 
\end{Remark}

\begin{Examples} (a) Important examples arise from 
Hilbert--Lie algebras  
$\fk$ (Definition~\ref{def:hilbertliealg}) 
with the positive definite invariant scalar product $\kappa_\fk$, 
\[  \g = C^\infty(\bS^1, \fk), \quad 
\kappa(\xi,\eta) := \frac{1}{2\pi} \int_0^{2\pi} \kappa_\fk(\xi(t),\eta(t))\, dt, \quad 
\bS^1 \cong \R/2\pi\Z. \] 
Then $D\xi := \xi'$ is a derivation leaving $\kappa$ invariant. 
If $K$ is a $1$-connected Lie group with Lie algebra $\fk$, 
then $G := C^\infty(\bS^1,K)$ is a $1$-connected Lie group 
with Lie algebra~$\g$ (\cite{NW09}). 
The action of $\R$ on $\g$ is given by 
\[  \alpha_t(\xi)(s) := \xi(t + s) \quad \mbox { and } \quad 
 \alpha^G_t(f)(s) := f(t + s). \]
This implies that 
\[  \gamma(g)(s) = - \derat0 g(t+s)g(s)^{-1} 
=  -\delta^r(g)_s, \]
so that $\gamma = \delta^r$ is the right logarithmic derivative. 

\nin (b) A more general class of examples arises from principal $K$-bundles 
$(P,M,q,K)$, where $K$ is a compact group, and a volume form 
$\mu$ on the compact manifold $M$. Then 
$$\g := \gau(P) \cong \{ \xi \in C^\infty(P,\fk) \: 
(\forall p \in P)(\forall k \in K)\, \xi(pk) = \Ad(k)^{-1} \xi(p)\} $$
carries an invariant scalar product given by 
$$ \kappa(\xi,\eta) := \int_M \tilde\kappa_\fk(\xi,\eta) \cdot \mu, $$
where we use that the function $\kappa_\fk(\xi,\eta)$ on $P$ 
is constant on the $K$-orbits, hence factors through a smooth function 
$\tilde\kappa_\fk(\xi,\eta)$ on $M$. 

If $X \in \cV(M,\mu)$ is a vector field with $\cL_X \mu = 0$, 
then any lift $\tilde X \in \aut(P) = \cV(P)^K$ defines a 
skew-symmetric derivation of $(\g, \kappa)$ and we may form 
the corresponding double extension. 

\nin (c) An interesting example is the following: 
We consider a left invariant vector field $X$ on $M = \T^2$ whose orbits 
are dense. Then $\ker D = \{0\}$ on 
$\g^0 := C^\infty(\T^2, \R)$. Similar examples arise from (irrational) 
quantum tori. 
\end{Examples}

\subsection{Double extensions of Hilbert--Lie algebras} 
\label{subsec:7.2} 

\begin{Example} \label{ex:hs} Let $\cH$ be a complex Hilbert space and 
$\g = \fu_2(\cH)$ be the Lie algebra of skew-hermitian Hilbert--Schmidt 
operators. 

Each continuous cocycle $\omega \in Z^2(\g,\R)$ can be written as 
$$ \omega(x,y) = \tr([d,x]y) $$
for some bounded skew hermitian operator $d \in \fu(\cH)$ 
(\cite{dH72}). Then 
$Dx := [d,x]$ defines a continuous derivation of  $\g$ 
preserving the invariant scalar product defined by $(x,y) := -\tr(xy)$. 
From that we easily derive that the cocycle $\gamma \: G = \U_2(\cH) \to \g$ 
is given by 
\[  \gamma(g) = gdg^{-1}- d.  \] 
If $d \in \R i \1 + \fu_2(\cH)$, then the cocycle $\omega$ is trivial 
and $\gamma$ is obviously  bounded. 

If, conversely, $\gamma$ is bounded, then the Bruhat--Tits 
Theorem~\ref{thm:BT}  implies that the affine isometric action of $G$ on 
$\g$ defined by $g * x := \Ad(g)x + \gamma(g)$  has a fixed 
point $x \in \g$. This means that 
$$ g(d + x)g^{-1} = d + x\quad \mbox{ for each } \quad g \in \U_2(\cH), $$
and hence that $d + x \in \C \1$, which leads to 
$d \in \fu_2(\cH) + \R i \1$. We conclude that each non-trivial 
$2$-cocycle $\omega$ leads to an unbounded cocycle $\gamma$. 

If $\gamma(G)$ is semi-equicontinuous, then 
$B(\im(\gamma))^0$ is an open invariant cone in 
$\fu_2(\cH)$, hence equal to the whole space 
(Proposition~\ref{prop:2.12}) and 
Remark~\ref{rem:2.3old} further entails that 
$\gamma$ is bounded. 
\end{Example}

\begin{Remark} Consider the Lorentzian double extension $\hat\g$ 
of $\g = \fu_2(\cH)$, defined by the skew-hermitian bounded operator $d$ 
(Example~\ref{ex:hs}). We want to classify semi-equicontinuous 
coadjoint orbits in $\hat\g'$ and open invariant cones in $\hat\g$. 

Beyond the general observations that hold for all Lorentzian 
double extensions (Proposition~\ref{prop:1.2}), this seems 
rather difficult. 
However, according to \cite[Satz~1]{NeJ35}, 
there exists a skew-hermitian Hilbert--Schmidt operator $s$ such that 
$d' := d + s$ is diagonalizable. Then the double extensions defined by 
$d$ and $d'$ are isomorphic, so that we may assume that 
$d$ is diagonalizable. Then the centralizer 
$\g^d := \{ x \in \g \: [d,x] = 0\}$ contains a maximal 
abelian subalgebra $\ft$ of $\g$ and 
\[ \hat\ft = \R \bc \oplus \ft \oplus \R \bd \] 
is an elliptic Cartan subalgebra of $\hat\g$. In the 
corresponding root system all roots are compact and $\g^{\rm alg}$ 
is a direct limit of compact Lie algebras. 
We refer to \cite{MN16} for a concrete analysis of this 
situation. 
\end{Remark}

\subsection{Twisted loop algebras} 
\label{subsec:twistloop}

Let $K$ be a connected Lie group whose Lie algebra 
$\fk$ is a simple Hilbert--Lie algebra. 
We fix $\Phi\in \Aut(K)$ with $\Phi^N = \id_K$, 
and write $\phi = \Lie(\Phi) \in \Aut(\fk)$. 
For $\tau := \frac{2\pi}{N}$, the 
corresponding {\it twisted loop group (of period $\tau$)} is 
\begin{equation}
  \label{eq:twistloop1}
\cL_\Phi^\tau(K) := \{ \xi \in C^\infty(\R,K)\: 
(\forall x \in \R)\ \xi(x + \tau) = \Phi^{-1}(\xi(x))\}
\end{equation}
with Lie algebra 
\begin{equation}
  \label{eq:twistloop1b}
\cL_{\phi}^\tau(\fk) := \{ \xi \in C^\infty(\R,\fk)\: 
(\forall x \in \R)\ \xi(x + \tau) = \phi^{-1}(\xi(x))\}, 
\end{equation}
where $\phi = \L(\Phi) 
\in  \Aut(\fk)$ is the automorphism of $\fk$ induced by~$\Phi$. 
We have a natural action of $\R$ on $\cL_\Phi(K)$ by
\begin{equation}
\label{eq:translat}
\alpha_t(\xi)(x) =\xi(x+t) \quad \mbox{ and } \quad 
D\xi = \frac{d}{dt}\Big|_{t = 0} \alpha_t(\xi) = \xi'.
\end{equation}
It satisfies $\alpha_{2\pi} = \id_{\cL_\Phi(K)}$, hence factors through an action 
of $\R/2\pi \Z \cong \T$. 

We write $\kappa$ for an invariant scalar product on $\fk$ 
normalized in such a way that 
$\kappa(i\check \alpha, i\check \alpha) = 2$ for long roots~$\alpha$. 
Accordingly, we obtain the double extension 
\[  \hat\cL^{\tau}_\phi(\fk) := 
(\R \oplus_{\omega} \cL^{\tau}_\phi(\fk)) \rtimes_{D} \R
\quad \mbox{ with } \quad D\xi = \xi', \] 
where 
\[ \omega(\xi,\eta) 
:=  \frac{c}{2\pi \tau}\int_0^\tau \kappa(\xi'(t),\eta(t))\, dt \] 
for some $c \in \Z$ (the {\it central charge}).  

Let  $\ft^\circ \subeq \fk^\phi$ be a maximal abelian 
subalgebra. Then $\fz_\fk(\ft^\circ)$ 
is an elliptic Cartan subalgebra of $\fk$ 
by \cite[Lemma~D.2]{Ne14}, and 
\begin{equation}
  \label{eq:ft}
\ft = \R \bc\oplus \ft^\circ \oplus \R \bd   
\end{equation}
is an elliptic Cartan subalgebra of $\hat\cL^{\tau}_\phi(\fk)$. 
The corresponding set of roots 
$\Delta\subeq i \ft'$ can be identified with the set of pairs 
$(\alpha,n)$, where 
\[ (\alpha,n)(z,h,s) := (0,\alpha,n)(z,h,s) = \alpha(h) + i s n, \quad 
n \in \Z, \alpha \in \Delta_n. \] 
Here $\Delta_n \subeq i(\ft^\circ)^*$ is the set of 
$\ft^\circ$-weights in the $\phi$-eigenspace 
\[ \fk_\C^n  = \{ x \in \fk_\C \: \phi^{-1}(x) = e^{in\tau} x\}.\] 

In this case all root vectors are of type (N) or (CS). 
If $\fk$ is finite-dimensional, then $\hat\g_\C^{\rm alg}$ is an affine Kac--Moody 
algebra, and in general it is a direct limit of affine Kac--Moody Lie algebras. 
For these Lie algebras, we have the following convexity theorem 
for coadjoint orbits. It is obtained by applying direct limit 
arguments to the topological version of the Kac--Peterson Convexity Theorem 
in Subsection~\ref{subsec:kacmoo}. We refer to P.~Helmreich's thesis \cite{He19} 
for details. 

\begin{Theorem} {\rm(Convexity Theorem for twisted loop groups)} 
For $\lambda \in \hat\ft'$ with $\lambda(\bc) \not=0$, we have 
\[ p_{\ft'}(\oline\conv\cO_\lambda) =  \oline\conv(\cW \lambda).\] 
\end{Theorem}

\begin{Problem} Classify the adjoint orbits of elements of the 
form $(z,x,t) \in \g$ with $t \not=0$. 
\end{Problem}

\subsection{Non-Lorentzian double extensions} 

To deal with double extensions of Hilbert--Lie algebras 
and of twisted loop algebras, one has to determine to which 
extent the rich geometric information available for 
Lorentzian double extensions is still available 
for non-Lorentzian ones. Here a promising approach is to 
consider a Lie algebra $\g$ with a positive definite form $\kappa_\g$ 
and two $\kappa_\g$-skew derivations $D_1, D_2 \in \der(\g,\kappa_\g)$. 
We further assume that the derivation $[D_1, D_2]$ is inner. One can show that 
in this case the Lorentzian double extension 
$\hat\g_1$ defined by $(\omega_{D_1}, D_1)$ carries a 
$\kappa_1$-skew symmetric derivation $\tilde D_2$ obtained by lifting 
$D_2$, so that we can form a double extension 
$\hats\g$ of $\hat \g_1$ defined by $(\omega_{\tilde D_2}, \tilde D_2)$. 
This results in a so-called {\it bidouble extension} which 
carries an invariant symmetric bilinear form which is non-Lorentzian 
of index $2$ but which contains two hyperplane ideals that are 
central extensions of Lorentzian double extensions. 
This geometric structure should be explored with respect to 
its impact on closed convex hulls of adjoint and coadjoint orbits.

\appendix 
\section{Appendices} 
\subsection{Some facts on convex sets} 
\label{app:a}

The following observation shows 
that semi-equicontinuous convex sets share many  
important properties with compact ones (cf.\ \cite[Prop.~6.13]{Ne08}): 

\begin{proposition} \label{prop:locomp} 
Let $C \subeq E'$ be a 
non-empty weak-$*$-closed convex semi-equicontinuous subset and $v \in 
B(C)^0$. Then $C$ is weak-$*$-locally compact, 
the function 
$\eta_v \: C \to \R, \eta_v(\alpha) := \alpha(v)$ 
is proper, and there exists an extreme point $\alpha \in C$ 
with $\alpha(v) = \min \la C, v \ra$. 
\end{proposition}

The following lemma (\cite[Cor.~II.2.6.1]{Bou07}) is often useful: 
\begin{Lemma}
  \label{lem:bou} 
For a convex subset $C$ of a locally convex space $E$, 
the sets $C^0$ and $\oline C$ are convex, 
$\oline C$ and $C$ have the same interior, 
and if this is non-empty, then $\oline{C^0} = \oline C$. 
\end{Lemma}

From the Hahn--Banach separation theorem, we obtain immediately: 
\begin{Proposition} \label{prop:bidual}
If $D \subeq E$ is a convex cone, then 
$\oline D = (D^\star)^\star$, and if $C \subeq E'$ is a convex cone,
 then $(C^\star)^\star = \oline{C}$ is its weak-$*$-closure. 
\end{Proposition}

The following lemma is often useful to obtain results 
on open convex subsets of normed spaces from those on closed ones. 
If $\dim E < \infty$, then the characteristic function 
\begin{equation}
  \label{eq:charfct}
 \phi_C(x) := \int_{B(C)} e^{-\alpha(x)+ \inf \alpha(C)}\, d\mu_{E^*}(\alpha) 
\end{equation}
of an open convex subset $C \subeq E$ 
has similar properties, but it is invariant under all 
affine automorphisms of $C$ preserving Lebesgue 
measure $\mu_{E'}$ (\cite[Thm.~V.5.4]{Ne99}). 

\begin{Lemma} \label{lem:1.3} 
Let $p \: E \to \R_+$ be a seminorm on the real vector space 
$E$ and $U \subeq E$ be a proper $p$-open convex subset. 
We write 
$B^p_r(x) := \{ y \in E \: p(x-y) < r\}$ 
for the open $p$-ball of radius $r$ and 
$d_{\partial U}(x) := \sup \{ r \geq 0 \: B^p_r(x) \subeq U \}$ 
for the distance of $x$ from $\partial U$ with respect to~$p$. 
Then 
\[ d^p_U \: U \to \R, \quad d^p_U(x) := d_{\partial U}(x)^{-1} \] 
is a continuous convex function on $U$ with 
$\lim_{x_n \to x \in \partial U} d^p_U(x_n) = \infty.$ 
In particular, all the set 
$U_c := \{ x \in U \: d^p_U(x) \leq c \}$
are closed convex subsets of $E$. 
\end{Lemma} 

\begin{proof} To see that $d^p_U$ is convex, let $B := B^p_1(0)$ 
denote the open $p$-unit ball in $E$. For 
$x,y \in U$ we have $x + {1 \over d^p_U(x)} B ,
y + {1 \over d^p_U(y)} B \subeq U$. For $0 < \lambda <1$, 
$z = \lambda x + (1- \lambda)y$
and $0 < p(v)< s := {\lambda \over d^p_U(x)} 
+ { 1-\lambda \over  d^p_U(y)}$, let 
$v_0 := {v \over s}$. Then $v_0 \in B$, so that 
\begin{align*}
z + v 
&= \lambda x + (1- \lambda)y + \Big({\lambda \over d^p_U(x)}+ {1- \lambda
\over d^p_U(y)}\Big)v_0\\
&\in \lambda\Big(x + {1 \over d^p_U(x)}B \Big) + (1 - \lambda) \Big(y + {1
\over d^p_U(y)} B\Big) \subeq U. 
\end{align*}
Hence $d^p_{\partial U}(z) \geq {\lambda  \over d^p_U(x)} + {1 - \lambda
\over d^p_U(y)}$. We conclude from the convexity and the antitony 
of the function $r \mapsto {1 \over r}$ for $r > 0$ that 
$ d^p_U(z) \leq \lambda d^p_U(x) + (1- \lambda) d^p_U(y), $
i.e.\ that $d^p_U$ is a convex function. 

The continuity of $d^p_U$ on $U$ follows from the continuity of 
the distance function $d^p_{\partial U}$ with respect to the semi-metric 
defined by $p$. If $x_n \to x \in \partial U$, then 
$d^p_{\partial U}(x_n) \to 0$, so that $d^p_U(x_n) \to \infty$. 
\end{proof}

\begin{Lemma} \label{lem:a.6} 
{\rm(Constructing open cones)} 
Let $E$ be a locally convex space. \\[-7mm]
\begin{itemize}
\setlength\itemsep{0em}
\item[\rm(a)] Let $H \subeq E$ be a closed hyperplane not containing~$0$, 
$E_0:= H -H$, and let 
$\eset\not=\Omega \subeq H$ be a relatively open convex subset. 
Then $W := \R_+^\times \Omega$ is an open convex cone in~$E$ 
with 
\[  \oline W \cap E_0 = \lim(\Omega) 
\quad \mbox{ and } \quad H(W) = H(\Omega).\] 
In particular, $W$ is pointed if and only if $\Omega$ contains no affine lines. 
\item[\rm(b)]
If $\Omega$ is a bounded open convex subset of the 
real locally convex space $E$ with $0\not\in\oline\Omega$, 
then $\R^\times_+\Omega$ is a pointed open cone. 
\end{itemize}
\end{Lemma} 

\begin{proof} (a) The convexity of $W$ follows immediately from the convexity 
of $\Omega$. 
Let $\lambda \: E \to \R$ be the unique continuous linear 
functional with $\lambda(H) = \{1\}$. Then $E_0 = \ker \lambda$, 
the map 
\[  \R_+^\times \times H \to  \lambda^{-1}(\R^\times_+), \quad 
(t,v) \mapsto tv  \] 
is a homeomorphism, and $W$ is the image of the open subset 
$\R_+^\times \times \Omega$, hence open. 

Let $w \in \oline W \cap E_0$ 
and pick a net $t_j w_j \to w$ with $w_j \in \Omega$. 
Then we have  ${t_j = \lambda(t_j w_j) \to \lambda(w) = 0}$, 
so that Lemma~\ref{lem:limcone}(ii) 
shows that $w \in \lim(\Omega)$. 
If, conversely, $0 \not=w \in \lim(\Omega)$ and 
$x \in \Omega$, then $x + \R_+ w \subeq  \Omega$, and therefore 
$\frac{1}{n}(x + nw) \to w$ implies that $w \in \oline W \cap E_0$. 
This proves that $\oline W \cap E_0 = \lim(\Omega)$. 

From $\lambda \in W^\star$ we  derive $H(W) \subeq \ker \lambda = E_0$, 
so that 
\[  H(W) = H(\oline W) = \oline W \cap (-\oline W) \cap E_0 
 = \lim(\Omega) \cap - \lim(\Omega) = H(\Omega).  \] 

\nin (b) To see that $W:= \R^\times_+\Omega$ is a convex cone, 
we observe that, for $x,y \in \Omega$ and $t,s > 0$, we have 
\[  tx + s y = (t+s)\Big(\frac{t}{t+s} x + \frac{s}{t+s}y\Big) 
\in (t+s) \Omega \subeq W.\] 

To see that $W$ is pointed, we use the Hahn--Banach--Separation Theorem 
and $0 \not\in\oline\Omega$ 
to find a continuous linear functional $\lambda \in E'$ with 
$\inf \lambda(\Omega)= 1$. 
Let $w\in\oline W$ and $t_j x_j \to w$ with $x_j \in \Omega$ and 
$t_j > 0$. Then $t_j \leq t_j \lambda(x_j) \to \lambda(w)$. 
By passing to a suitable subnet, we may thus 
assume that $t_j \to t_0 \geq 0$. If $t_0 = 0$, then 
the boundedness of $\Omega$ implies that $w = 0$. 
If $t_0 > 0$, then $x_j \to t_0^{-1} w \in \oline\Omega$ 
and $\lambda(w) = t_0 \lim \lambda(x_j) \geq t_0$. 
This shows that $\oline W \setminus \{0\} \subeq \lambda^{-1}(\R^\times_+)$, 
and hence that $\oline W \cap - \oline W = \{0\}$. 
\end{proof}

\begin{Lemma} \label{lem:a.2} 
{\rm(Enlarging semi-equicontinuous subsets 
by cones)} Let $E$ be a locally convex space, $\Omega \subeq E$ be a non-empty 
open convex cone and let $C \subeq E'$ be a weak-$*$-closed convex 
subset with $\Omega \subeq B(-C)$. 
Then $C - \Omega^\star$ is weak-$*$-closed semi-equicontinuous, and 
a convex subset $C_1 \subeq E'$ is contained in 
$C - \Omega^\star$ if and only if $s_{C_1}\res_\Omega \leq s_C\res_\Omega.$ 
\end{Lemma}

\begin{proof} The convex subset 
$C_2 := C - \Omega^\star \subeq E'$. It clearly satisfies 
$s_{C_2} = s_C$ on $\Omega$, so that any subset ${C_1} \subeq {C_2}$ satisfies 
$s_{C_1} \leq s_{C_2} = s_C$ on $\Omega$. 

Assume, conversely, that $s_{C_1} \leq s_C$ holds on $\Omega$. 
From Proposition~\ref{prop:locomp} we know that 
$C$ and $\Omega^\star$ are weak-$*$-locally compact, and 
for every $x \in \Omega \subeq B(-C)$, the evaluation maps
$\eta_x \: C \to \R$ and $\eta_x \: -\Omega^\star \to \R$ are 
proper and bounded from above. This implies that the function 
\[ f \: C \times \Omega^\star \to \R, \quad (\lambda, \beta) \mapsto 
(\lambda - \beta)(x) \] 
is also proper and bounded from above. It follows that 
the subset ${C_2} = C - \Omega^\star \subeq E'$ is weak-$*$-closed, 
and since $s_{C_2} = s_C$ holds on $\Omega$, we see that ${C_2}$ is semi-equicontinuous. 

For any $y \in B(-C_2)$, the evaluation map $\eta_y$ must be bounded 
from above on the convex cone $- \Omega^\star$, which implies that 
$y \in (\Omega^\star)^\star = \oline \Omega$. This proves that 
$\Omega \subeq B(-C_2)\subeq \oline \Omega$, and since 
$B(-C_2)^0$ is open and convex, we obtain $\Omega = B(-C_2)^0$ from 
Lemma~\ref{lem:bou}. 
Now \eqref{eq:recov} in the proof Theorem~\ref{thm:2.3} yields 
\begin{align*}
 C - \Omega^\star = {C_2} 
&=  \{ \lambda \in E' \: (\forall v \in B(-C_2)^0)\,  \lambda(v) 
\leq s_{C_2}(v)\} \\
&=  \{ \lambda \in E' \: (\forall v \in \Omega)\ \lambda(v) \leq s_C(v)\},
\end{align*}
and our assumption thus implies that ${C_1} \subeq C - \Omega^\star$. 
\end{proof}

\begin{Lemma}
  \label{lem:new} Let $(E, \tau)$ be a locally convex space 
and $\tau' \subeq \tau$ be a coarser  locally convex topology on $E$. 
If $\Omega \subeq E$ is $\tau$-open and convex and contains a 
$\tau'$-inner point, then $\Omega$ is also $\tau'$-open. 
\end{Lemma}

\begin{proof} Let $x_0 \in \Omega$ be $\tau'$-inner and 
$U\subeq \Omega$ be an open convex $\tau'$-neighborhood of~$x_0$. 
If $x_1 \in \Omega$ is any other point, there exists an $\eps > 0$ with 
$x_2 := x_1 + \eps(x_1 - x_0) = (1 + \eps) x_1 - \eps x_0 \in \Omega$ 
because $\Omega$ is $\tau$-open. Then 
\[ x_1 
= \frac{1}{1 + \eps}(\eps x_0 + x_2) \in 
 \frac{\eps}{1 + \eps}U + \frac{1}{1+ \eps} x_2 \subeq \Omega \] 
shows that $\Omega$ is also a $\tau'$-neighborhood of $x_1$. 
\end{proof}

\nin {\bf Some Lorentzian geometry.} 
Let $E$ be a locally convex space and 
$\beta \: E \times E \to \R$ be a {\it Lorentzian form}, i.e., 
$\beta$ is symmetric bilinear, and there exists 
$v_0 \in E$ with $\beta(v_0, v_0) > 0$ 
such that $\beta$ is negative definite 
on $v_0^\bot$. Note that this implies that the same holds for any 
$v \in E$ with $\beta(v,v) > 0$ 
(cf.\ \cite[Cor.~IV.7.5]{Bog74}). 
Then $\{ v \: \beta(v,v) > 0\}$ is the union of two 
open convex cones, and 
$$ W := \{ v \in E \: \beta(v,v) > 0, \beta(v_0, v) > 0\} $$
is one of them. 

\begin{Lemma} \label{lem:ics} {\rm(Inverse CS inequality)} 
For $\beta(v,v) \geq 0$ and $x \in E$ we have 
$$ \beta(v,v) \beta(x,x) \leq \beta(x,v)^2. $$
\end{Lemma}

\begin{proof} If $\beta(v,v) = 0$, there is nothing to show. So we assume 
that $\beta(v,v) > 0$ and we may w.l.o.g.\ assume that 
$\beta(v,v) =1$. Writing $E = \R v \oplus v^\bot$ and, 
accordingly, $x = \lambda v + w$ with $\lambda = \beta(v,x)$, we have 
$$ \beta(v,v)\beta(x,x) = \beta(x,x) = \lambda^2 + \beta(w,w) 
\leq \lambda^2 = \beta(x,v)^2 $$
because $\beta$ is negative definite on $v^\bot$. 
\end{proof}

\begin{Proposition} \label{prop:c.2} The smooth function 
$\chi \: W \to \R, \chi(v) := \beta(v,v)^{-1}$
has the following properties:\\[-7mm] 
\begin{description}
\setlength\itemsep{0em}
\item[\rm(i)] $\chi(v+ w) \leq \chi(v)$ for $v,w \in W$. 
\item[\rm(ii)] If $v_n \to v \in \partial W$, then $\chi(v_n) \to \infty$. 
\item[\rm(iii)] $\chi$ is strictly convex.  
\end{description}
\end{Proposition}

\begin{proof} (i) We have to show that,  for $v,w \in W$, we have 
$$ \beta(v + w, v + w) 
= \beta(v,v) + \beta(w,w) + 2\beta(v,w)
\geq \beta(v,v). $$
Since $\beta(w,w) > 0$, it suffices to see that 
$\beta(v,w) > 0$. To verify this inequality, 
pick $v_0 \in W$ with $\beta(v_0, v_0) = 1$ and 
write $E = \R v_0 \oplus v_0^\bot$. Accordingly, we have 
$v = \lambda v_0 + v'$ and 
$w = \mu v_0 + w'$ with $\lambda, \mu > 0$, so that 
$$\beta(v,w) = \lambda \mu + \beta(v',w')
\quad \mbox{ and } \quad  0 < \beta(v,v) = \lambda^2 + \beta(v',v'). $$
Applying the Cauchy--Schwarz inequality to the restriction of $-\beta$ 
to $v_0^\bot$, we get 
$$ \beta(v',w')^2 \leq \beta(v',v')\beta(w',w') 
< \lambda^2 \mu^2, $$
which leads to $\beta(v,w) > 0$. 

\nin (ii) is a trivial consequence of the definition. 

\nin (iii) For $v \in W$ we have 
\[  (\partial_x \chi)(v) = - \frac{2}{\beta(v,v)^2} \beta(v,x), 
\quad \mbox{ and } \quad  (\partial_x^2 \chi)(v) 
= \frac{2}{\beta(v,v)^3}
\Big(2\beta(v,x)^2 - \beta(v,v) \beta(x,x)\Big). \]
From the inverse Cauchy--Schwarz inequality (Lemma~\ref{lem:ics}) we derive 
$$ \beta(x,x)\beta(v,v) \leq \beta(x,v)^2, $$
so that $\partial_x^2 \chi(v)\geq 0$ for every $x \in E$. 
The relation $(\partial_x^2\chi)(v) = 0$ implies 
$$ \beta(v,x)^2 + \big(\beta(v,x)^2 - \beta(v,v)\beta(x,x)\big) = 0, $$
and since both summands are non-negative, we obtain 
$\beta(v,x) = 0$. Then the relation 
$-\beta(v,v)\beta(x,x) = 0$ further  leads to 
$\beta(x,x) = 0$. As $\beta$ is negative definite on $v^\bot$, 
we finally obtain $x = 0$. 
\end{proof}

\subsection{The Bruhat--Tits Theorem} 

\begin{Definition} Let $(X,d)$ be a metric space. We say that 
$(X,d)$ satisfies the {\it semi parallelogram law} if, 
for $x_1, x_2 \in X$, there exists a point $z \in X$ such that, for 
each $x \in X$, we have 
\begin{align}
  \label{eq:spl}
d(x_1, x_2)^2 + 4 d(x,z)^2 \leq 2 d(x, x_1)^2 + 2 d(x, x_2)^2.  \tag{SPL}
\end{align}
\end{Definition}

Every euclidean space satisfies the 
parallelogram law if $z$ is chosen to be the midpoint of $x_1$ and $x_2$, 
and this trivially implies the semi parallelogram law. 
Accordingly, the point $z$ occurring in the preceding definition 
plays the role of a ``midpoint'' of $x_1$ and $x_2$. 

\begin{Remark}  Suppose that $(E, \|\cdot\|)$ is a normed space endowed
with the metric $d(x,y) := \|x - y\|$. When does $(E,d)$ satisfy the
semi parallelogram law? This is the case if 
the norm is defined by a positive definite symmetric bilinear form
$\la \cdot, \cdot \ra$ via $\|x\| := \sqrt{ \la x, x \ra}$. 
If, conversely, the semi parallelogram law is satisfied, then we
obtain with $x_1 = - x_2 = a$ and $x = b$ for all $a,b \in E$ the relation 
$$ 4 \|a\|^2 + 4\|b\|^2 \leq 2 \|a + b\|^2 + 2 \|a - b \|^2, $$
which in turn leads to 
$$ 2 \|a + b\|^2 + 2 \|a - b \|^2 
\leq \|(a+b) + (a-b)\|^2 + \|(a+b)-(a-b)\|^2
= 4 \|a\|^2 + 4\|b\|^2, $$
and therefore to the parallelogram law 
\[  2 \|a\|^2 + 2\|b\|^2 =  \|a + b\|^2 + \|a - b \|^2 
\quad \mbox{ for all } \quad a,b \in E.\]
 By a theorem of P.~Jordan and J.~von Neumann 
(\cite{JvN35}), 
each normed satisfying the parallelogram law is euclidean. In fact, 
$$ \la a,b \ra := {\|a + b\|^2 - \|a - b \|^2 \over 4} $$
is a positive definite symmetric bilinear form with 
$\|a\| = \sqrt{\la a,a\ra}$. Therefore the euclidean spaces are 
precisely the normed spaces satisfying the semi parallelogram law. 
\end{Remark}

\begin{Definition}
 A {\it Bruhat--Tits space} is a complete metric space 
$(X,d)$ satisfying the semi parallelogram law. 
\end{Definition}

We write $\Aut(X,d)$ for the group of {\it automorphisms of the metric
space $(X,d)$}. For an elementary proof of the following theorem 
we refer to \cite{La99}. 

\begin{Theorem}\label{thm:BT}  {\rm(Bruhat--Tits Fixed Point Theorem)} Let $X$ be a
Bruhat--Tits space and $G \subeq \Aut(X,d)$ a subgroup with a bounded
orbit. Then $G$ has a fixed point. 
\end{Theorem}

\begin{Example} \label{ex:BTexample}
 If $E$ is a real Hilbert space and $C \subeq E$ a closed convex 
subset, then $C$ is a Bruhat--Tits space with respect to the euclidean metric. 
\end{Example}

\subsection{A general lemma} 

The following lemma (\cite[Lemma~6.15]{Ne10}) captures the spirit of various 
constructions of invariant convex cones in Lie algebras in a quite natural way. 
F.i., it applies to finite-dimensional simple 
algebras as well as $\cV(\bS^1)$ (\cite[\S 8]{Ne10}). 

\begin{Lemma}\label{lem:doubcone} 
Suppose that the element $\bd \in \g$ has the following 
properties: 
\begin{itemize}
\setlength\itemsep{0em}
\item[\rm(a)] The interior $W_{\rm min}$ of the invariant convex cone 
generated by $\cO_d = \Ad(G)\bd$ is non-empty and different from $\g$. 
\item[\rm(b)] There exists a continuous linear  projection 
$p_\bd \: \g \to \R \bd$ which preserves every open and closed convex invariant 
subset. 
\item[\rm(c)] There exists an element $x \in \g$ for which 
$p_\bd(\cO_x)$ is unbounded. 
\end{itemize}
Then each non-empty open invariant cone $W$ contains $W_{\rm min}$ or 
$-W_{\rm min}$, and for $\lambda \in \g'$ the following are equivalent
\begin{description}
\item[\rm(i)] $\lambda \in \g'_{\rm seq}$, i.e., 
$\cO_\lambda$ is semi-equicontinuous. 
\item[\rm(ii)] $\cO_\lambda(d)$ is bounded from below or above. 
\item[\rm(iii)] $\lambda \in W_{\rm min}^\star \cup - W_{\rm min}^\star$. 
\end{description}
\end{Lemma}

Typical classes of Lie algebras to which 
this lemma applies are finite-dimensional 
hermitian simple Lie algebras with $\fz(\fk) = \R \bd$ 
(\cite{Ne99}) and 
$\cV(\bS^1)$ with the generator $\bd$ of rigid rotations of the circle 
(cf.~Subsection~\ref{subsec:5.3}). For further 
applications we refer to the restricted metaplectic group 
discussed in \cite[\S 9.2]{Ne10}.

\end{document}